%% file: main.tex
\documentclass[10pt]{article}
\usepackage[a4paper,margin=1in]{geometry}

\usepackage[T1]{fontenc}
\usepackage[utf8]{inputenc}
\usepackage{amsmath,amssymb,mathrsfs,bm}
\usepackage{mathtools} 
\usepackage{amsthm}
\usepackage{dsfont}
\usepackage{stmaryrd} 
\usepackage[mathscr]{eucal}

\usepackage{graphicx}
\usepackage[table,xcdraw]{xcolor}
\usepackage{booktabs}
\usepackage{multirow,pbox}
\usepackage{caption,subcaption}
\usepackage{float}
\usepackage{algorithm}
\usepackage{algpseudocode}

\usepackage{enumitem}
\setlist{leftmargin=9mm}

\usepackage[square, authoryear]{natbib}
\usepackage[colorlinks,linkcolor=blue,citecolor=blue,urlcolor=blue,pagebackref]{hyperref}

\renewcommand*{\backrefalt}[4]{%
    \ifcase #1 \footnotesize{(not cited)}%
    \or \footnotesize{(cited on page~#2)}%
    \else \footnotesize{(cited on pages~#2)}%
    \fi}

\allowdisplaybreaks[4]
\numberwithin{equation}{section}

\newcommand*{\dif}{\mathop{}\!\mathrm{d}}
\newcommand{\smallop}{\mathfrak{o}_{\mathbb{P}}}
\newcommand{\bigop}{\mathcal{O}_{\mathbb{P}}}
\newcommand{\smallo}{\mathfrak{o}}
\newcommand{\bigo}{\mathcal{O}}
\newcommand{\iu}{\mathrm{i}} 
\renewcommand{\epsilon}{\varepsilon}

\newcommand{\pto}{\xrightarrow{\mathtt{p}}}
\newcommand{\dto}{\xrightarrow{\mathtt{d}}}
\DeclarePairedDelimiter{\intset}{\llbracket}{\rrbracket}
\makeatletter
\renewcommand*{\top}{{\mathpalette\@transpose{}}}
\newcommand*{\@transpose}[2]{\raisebox{\depth}{$\m@th#1\intercal$}}
\makeatother
\DeclareMathOperator{\Tr}{Tr}

\DeclareMathOperator{\E}{\mathbb{E}}
\DeclareMathOperator{\Prob}{\mathbb{P}}

\DeclareMathOperator{\var}{Var}
\DeclareMathOperator{\cov}{Cov}

\DeclareMathOperator{\diag}{diag}

\let\limsup\relax
\DeclareMathOperator*\limsup{\overline{lim}}

\DeclarePairedDelimiter{\floor}{\lfloor}{\rfloor}

\renewcommand{\hat}{\widehat}
\renewcommand{\tilde}{\widetilde}
\newcommand{\one}{\mathds{1}} 

\newcommand{\bbR}{\mathbb{R}}

\newcommand{\calF}{\mathcal{F}}
\newcommand{\calG}{\mathcal{G}}
\newcommand{\deq}{\coloneqq}
\newcommand{\eqd}{\eqqcolon}

\newcommand{\mycomment}[1]{}

\theoremstyle{plain}
\newtheorem{theorem}{Theorem}[section] 
\newtheorem{lemma}[theorem]{Lemma}
\newtheorem{proposition}[theorem]{Proposition}

\theoremstyle{definition}

\theoremstyle{remark}
\newtheorem{remark}[theorem]{Remark}

\theoremstyle{plain}
\newtheorem{assumption}{\textbf{Assumption}}
\renewcommand*{\theassumption}{\Alph{assumption}}

\AtBeginDocument{%
	\def\MR#1{}
}

\title{\textbf{A spectral approach for online covariance change point detection}}
\author{Zhigang Bao\thanks{University of Hong Kong, \texttt{zgbao@hku.hk}},
Kha Man Cheong\thanks{Hong Kong University of Science and Technology, \texttt{kmcheong@connect.ust.hk}},
Yuji Li\thanks{University of Hong Kong, \texttt{u3011732@connect.hku.hk}},
Jiaxin Qiu\thanks{University of Hong Kong, \texttt{jxqiu@hku.hk}}
}

\begin{document}
\maketitle

\begin{abstract}

Change point detection in covariance structures is a fundamental and crucial problem for sequential data. Under the high-dimensional setting, most of the existing research has focused on identifying change points in historical data. However, there is a significant lack of studies on the practically relevant online change point problem, which means promptly detecting change points as they occur.
In this paper, applying the limiting theory of linear spectral statistics for random matrices, we propose a class of spectrum based CUSUM-type statistic. We first construct a martingale from the difference of linear spectral statistics of sequential sample Fisher matrices, which converges to a  Brownian motion. Our CUSUM-type statistic is then defined as the maximum of a variant of this process. Finally, we develop our detection procedure based on the invariance principle. Simulation results show that our detection method is highly sensitive to the occurrence of change point and is able to identify it shortly after they arise, outperforming the existing approaches.

\end{abstract}

\setcounter{tocdepth}{2}
\tableofcontents

\section{Introduction}

Online change point detection plays a pivotal role in statistical process control and time series analysis. 
The field has a rich history, rooted in the classical works on Page's cumulative sum (CUSUM) procedure \citep{Page1954continuous} and Shiryaev-Roberts procedure \citep{Shiryaev1963optimum, Roberts1966comparison}. 
For a detailed exposition of the methodology and a survey of recent developments, we refer readers to \citet{Aue2024state} and the references therein. 
While classical theories are well-established, monitoring the covariance structure in high-dimensional settings remains a critical yet challenging task, as it captures the dynamic dependencies among a large number of variables.
Motivated by this, we consider the problem of online change point detection for the covariance structure of high-dimensional random vectors. 
Specifically, we observe a sequence of $p$-dimensional random vectors $\{y_i\}_{i=1}^{\infty}$, modeled as $y_i = \Sigma_i^{1/2} x_i$. 
Here, $\Sigma_i$ is a deterministic $p\times p$ symmetric positive definite matrix with uniformly bounded spectral norms, and the entries of $x_i$ are independent and identically distributed (i.i.d.) with zero mean and unit variance. More generally, we can potentially consider the model $y_i = T_i x_i$, where $T_i$ satisfies $\Sigma_i = T_i T_i^*$. But for the sake of technical conciseness, we work with the mentioned setting.
Our objective is to sequentially monitor the stream and detect whether the covariance structure $\Sigma_i$ changes at some unknown time instant. 
Formally, we test the following hypothesis
\begin{align*}
    H_0: \Sigma_i = \Sigma_0\quad \forall i \geq 1 
    \quad \text{vs} \quad 
    H_1: \exists k^{\star} \text{ s.t. } 
    \begin{cases}
        \Sigma_i = \Sigma_0, & 1\leq i \leq k^{\star},\\
        \Sigma_i = \Sigma_1 \neq \Sigma_0, & i > k^{\star}.
    \end{cases} 
\end{align*}
Here, $k^{\star}$ denotes the change point location.
Our goal is to detect the change as promptly after the change point $k^{\star}$ while controlling the false alarm rate under the null hypothesis $H_0$. 

The problem of detecting structural changes in covariance matrices has been extensively investigated, primarily in the offline setting where the complete data sequence is available for analysis. In the classical fixed-dimensional regime, a multitude of methods have been developed based on likelihood ratios and CUSUM; see, for example, \citet{Chen2004statistical}, \citet{Lavielle2006detection}, \citet{Galeano2007covariance}, \citet{Aue2009break}, \citet{Xie2013sequential}, \citet{Dette2015detecting},
For a comprehensive overview of these classical approaches, we refer readers to \citet{csorgo1997limit}, \citet{chen2000parametric}, \citet{Aue2012structural} and the references therein. 
More recently, there has been a surge of literature addressing high-dimensional offline covariance change-point detection from various perspectives, including projection based technique \citep{steland2020testing}, binary segmentation algorithm \citep{Wang2021optimal,li2023detection}, $U$-statistics based method \citep{Liu2020unified, Cui2025simultaneous}, graph-based method \citep{Chen2015graph,Chu2019asymptotic}, dimension reduction technique \citep{Dette2022estimating}, self-normalization principle \citep{Bours2021large, Wang2022inference}, and spectral method \citep{Ryan2023detecting, dornemann2024detecting, Dornemann2024detecting_spike}.

In contrast, online monitoring for high-dimensional covariance structures is still relatively limited in the literature.
Unlike offline methods that utilize the full data set to estimate the change point, online procedures must make decisions in real-time with limited historical data and no knowledge of the future.
This constraint requires statistics that are recursive, computationally efficient, and able to control false alarm rates. 
\citet{Avanesov2018change} and \citet{Avanesov2019structural} investigated the problem of change point detection for high-dimensional precision and covariance matrices within a multiple-testing framework. 
Their procedure constructs a collection of local window-based statistics based on $\ell_{\infty}$-norm, with critical values determined via high-dimensional Gaussian approximation and bootstrap techniques.
Addressing a specific structural change, \citet{Xie2020sequential} focused on sequentially detecting deviations from an identity matrix toward an unknown low-rank (spiked) covariance model. 
To this end, they proposed two online detection procedures: a Largest-Eigenvalue Shewhart chart and a Subspace-CUSUM procedure.
\citet{Li2023online} studied the online change point detection of high-dimensional covariance matrices utilizing a window-based $U$-statistic coupled with a sequential stopping rule. 
A key advantage of their method is that it admits explicit control of the Average Run Length (ARL) while accommodating non-Gaussian data with spatial and temporal dependence.
Most recently, \citet{Gao2025online} proposed a framework for the online detection of changes in the correlation structure of high-dimensional streaming data. 
Their approach leverages $\ell_2$- and $\ell_{\infty}$-norm based statistics, employing a sign-flip permutation method to determine adaptive detection thresholds.

In this work, we take a different perspective rooted in Random Matrix Theory. Specifically, we construct a sequential monitoring procedure based on the linear spectral statistics of the Fisher matrix, which offers a powerful tool for distinguishing the null and alternative hypotheses in high-dimensional settings.
Suppose that at the current time $k$, we have observed the sequence $y_1, \ldots, y_k$. 
To detect potential changes, we compare the covariance structures of two adjacent segments partitioned by an index $k_1$, where $1 < k_1 < k$. 
We define the reference sample covariance matrix $S_1$ and the monitoring sample covariance matrix $S_{2,k}$ as follows:
\[
    S_1 = \frac{1}{k_1} \sum_{i=1}^{k_1} y_i y_i^{\top}, 
    \qquad 
    S_{2,k} = \frac{1}{k - k_1} \sum_{i=k_1 + 1}^{k} y_i y_i^{\top}.
\]
We then construct the multivariate \emph{$F$-matrix} (or \emph{Fisher matrix}) defined as $F_k = S_1^{-1} S_{2,k}$,
and consider its \emph{linear spectral statistic} (LSS)
\[
    \Tr f(F_k) \deq \sum_{i=1}^{p} f\bigl(\lambda_i(F_k)\bigr),
\]
where $\lambda_1(F_k) \geq \cdots \geq \lambda_p(F_k)$ denote the ordered eigenvalues of $F_k$ and $f:\mathbb{R}\to\mathbb{R}$ is a suitable test function. 
The spectral properties of the $F$-matrix have been extensively investigated within the framework of random matrix theory. 
The Limiting Spectral Distribution (LSD) was established in the foundational works of \citet{bai1988on} and \citet{Bai1986on}, and subsequently extended by \citet{Silverstein1995strong} and \citet{Zhang2022limiting}.
Furthermore, the fluctuations of its LSS have been derived under various settings; see, for instance, \citet{Zheng2012central} and \citet{Zheng2017CLT}. 
Moreover, the asymptotic behavior of the largest eigenvalue has attracted considerable attention. Relevant results include \citet{Han2016tracy,Han2018unified} and \citet{Wang2017extreme}, with more recent advancements discussed in \citet{Jiang2021limits, Xie2021limiting, Hou2023spiked, Wang2025inference} and \citet{Jiang2026invariance}.

The core of our methodology relies on monitoring the incremental differences of the LSS. To this end, we define the standardized increment
\[
    \widetilde{L}_k(f) \deq 
    \frac{ \Tr f(F_k) - \Tr f(F_{k-1}) - \mu_k}{\sigma_k}, \label{onstepdiff}
\]
where $\mu_k$ and $\sigma_k$ represent the mean and standard deviation of the one-step difference $\Tr f(F_k) - \Tr f(F_{k-1})$ under the null hypothesis $H_0$. The explicit forms of $\mu_k$ and $\sigma_k$ will be provided in Proposition \ref{pro:stepone_diff_mean_var}.
Under $H_0$ of no change point, the expectation of $\widetilde{L}_k(f)$ is asymptotically negligible.
Once a structural change occurs at some unknown time $k^{\star}$,  the expectation of $\widetilde{L}_k(f)$ deviates systematically from zero.
To effectively accumulate the signal from such potentially subtle drifts, we further construct a normalized CUSUM statistic $T_p(n,i)$ (see definition in \eqref{eq:online_CUSUM}) based on the sequence $\{\widetilde{L}_k(f)\}$.
Under $H_0$, the process $\{T_p(n,i)\}$ converges in distribution to a  Brownian motion type process. Such convergence is based on the standardization in (\ref{onstepdiff}), which makes the increment behave as a normalized martingale difference; thus, the CUSUM statistic behaves as a random walk. 
When a change point exists ($H_1$), the mean shift in the LSS increments will cause $T_p(n,i)$ to deviate significantly from zero, thus triggering detection.
Accordingly, we define the \emph{stopping rule}
\[
\hat{k}^{\star} = \inf\{i: |T_p(n,i)| > c_{\alpha}\},
\]
where the threshold $c_{\alpha}$ is chosen to control the false alarm rate at a prescribed level $\alpha$. The value of $c_{\alpha}$ is determined by the quantile of the limiting null distribution of the process $\{T_p(n,i)\}$.

The rest of this paper is organized as follows. 
Section \ref{sec:main_results} presents our main theoretical results, including the explicit expressions for the mean and variance of the one-step LSS difference under $H_0$, as well as the weak convergence of the associated Brownian motion type process. 
Based on these results, we propose our online change point detection procedure. 
In Section \ref{sec:simulation}, we conduct extensive simulation studies to evaluate the performance of our proposed method in various scenarios. 
In Section \ref{sec:real_data}, we apply our proposed method to the S\&P 500 dataset to demonstrate its practical utility.
To facilitate reproducibility, the source code for the simulation studies and real data analysis is available online\footnote{Code available at: \url{https://github.com/jxqiu77/OnlineCovCPD}.}.
All technical proofs are provided in Sections \ref{sec:proof_thm_BM} -- \ref{sec:proof_lem_trace_limit_Sigma} and the Appendix.

Before moving forward, let us introduce some notations that will be used throughout this paper.
The symbol $X \deq Y$ (or $Y \eqd X$) indicates that $X$ is defined as $Y$. 
We use double brackets to denote index sets, i.e., $\intset{n_1, n_2} : = \mathbb{Z} \cap [n_1, n_2]$ for $n_1, n_2 \in \mathbb{R}$. 
For brevity, we write $\intset{n}$ for $\intset{1, n}$ when $n$ is a positive integer.
We use $\one\{E\}$ or $\one_{E}$ to denote the indicator function of an event $E$.
For any matrix $A$, we denote its $(i,j)$-th entry by $A_{ij}$, its transpose by $A^{\top}$, its Hermitian transpose by $A^*$, its trace by $\Tr(A)$, its $j$-th largest eigenvalue by $\lambda_j(A)$, its spectral norm by $\|A\|_2$. 
For two sequences of positive numbers $\{a_n\}_{n=1}^{\infty}$ and $\{b_n\}_{n=1}^{\infty}$, we write $a_n\ll b_n$ (or $b_n\gg a_n$) if $a_n/b_n\to 0$. 
We write $a_n\gtrsim b_n$ when $a_n\geq cb_n$ for some absolute constant $c>0$, and $a_n\lesssim b_n$ when $b_n\gtrsim a_n$. 
We write $a_n\asymp b_n$ if $c_0\leq a_n/b_n \leq C_0$ for some absolute constants $c_0, C_0 >0$. 
For random variables $\{X_n\}_{n=1}^{\infty}$ and positive real numbers $\{a_n\}_{n=1}^{\infty}$, we write $X_n=\bigop(a_n)$ if $X_n/a_n$ is bounded in probability, and $X_n=\smallop(a_n)$ if $X_n/a_n\pto 0$.
We use $X_n \pto X$ and $X_n \dto X$ to denote convergence in probability and convergence in distribution, respectively.

\section{Main Results}\label{sec:main_results}

\subsection{One-step LSS difference}

In this section, we present our main theoretical results concerning the one-step LSS difference of the Fisher matrix.
These results form the foundation for our online change point detection procedure.
To facilitate the analysis, we first impose the following assumptions.

\begin{assumption}[Data generating model]\label{asm:data}
    For any $i\geq 1$, let $y_i = \Sigma_i^{1/2} x_i$, where $\{\Sigma_i\}_{i=1}^{\infty}$ are deterministic $p\times p$  positive-definite symmetric matrix with uniformly bounded spectral norms. The vectors $\{x_i\}_{i=1}^{\infty}$ are i.i.d. with $x_i=(x_{ij})_{j=1}^p$, satisfying $\E x_{11} = 0$, $\E x_{11}^2 = 1$, and $\E x_{11}^{4+\delta} < \infty$ for some $\delta >0$.
\end{assumption}

\begin{assumption}[High-dimensional scaling]\label{asm:sample_dimension}
    The sample sizes $k_1, k$ and the dimension $p$ satisfy
    \[ 
        \frac{p}{k_1} \to c_1 \in (0, 1) \quad \text{and} \quad \frac{p}{k-k_1} \to c_2 \in (0, + \infty),
    \]
    as $p, k, k_1 \to \infty$.
\end{assumption}

Let $\{y_i\}_{i\geq 1}$ be a sequence of $p$-dimensional random vectors satisfying Assumptions \ref{asm:data} and \ref{asm:sample_dimension}, with $\Sigma_i = I_p$ for all $i\geq 1$ under $H_0$. 
We define the two sample covariance matrices as
\[
    S_1 \deq \frac{1}{k_1} \sum_{i=1}^{k_1} y_i y_i^{\top}, 
    \quad 
    S_{2,k} \deq \frac{1}{k-k_1} \sum_{i=k_1+1}^k y_i y_i^{\top}.
\]
The corresponding $F$-matrix is given by
\[
    F_k \deq S_1^{-1}S_{2,k}.
\]
From \citep{bai1988on,Bai1986on}, the LSD of $F_k$ is the distribution with the density function given by
\[
    f_{c_1,c_2}(x) = \frac{(1-c_1)\sqrt{(b-x)(x-a)}}{2\pi x(c_1x+c_2)}, \quad a \leq x \leq b,
\]
where
\[
    a \deq \frac{(1-h)^2}{(1-c_1)^2}, \quad 
    b \deq \frac{(1+h)^2}{(1-c_1)^2}, \quad
    h \deq \sqrt{c_1+c_2-c_1c_2}.
\]
The Stieltjes transform of the LSD is
\begin{align*}
    m(z) 
    & = \frac{1-c_2}{zc_2} -\frac{ c_2 \{z(1-c_1)+1-c_2\} + 2zc_1 - c_2(1-c_1)\sqrt{ (z-b)(z-a) } }{2zc_2(c_2 + zc_1)},
\end{align*}
where $z$ is any complex number with positive imaginary part.
The non-zero eigenvalues of $S_1^{-1}S_2$ and its companion matrix $\mathscr{F}_k \deq\frac{1}{k-k_1}X_2^{\top} S_1^{-1}X_2$ are the same, where $X_2 = (x_{k_1+1}, \ldots, x_k)$ is the $p\times (k-k_1)$ data matrix formed by the second sample. The Stieltjes transform of $\mathscr{F}_k$ is given by
\begin{align}
    \underline{m}(z)  & = -\frac{1-c_2}{z} + c_2 m(z)\notag\\
    & = -\frac{ c_2 \{z(1-c_1)+1-c_2\} + 2zc_1 - c_2(1-c_1)\sqrt{ (z-b)(z-a) } }{2z(c_2 + zc_1)}. \label{eq:mbar}
\end{align}
For any test function $f$, the \emph{linear spectral statistic} (LSS) of $F_k$ is defined as
\[
    \Tr f(F_k) \deq \sum_{i=1}^p f(\lambda_i(F_k)). 
\]
Let $L_k(f) \deq \Tr f(F_k) - \Tr f(F_{k-1})$ be the \emph{one-step difference} of LSS.
By using the Cauchy integral formula, we have for analytic $f$, 
\begin{align}\label{eq:Lk_Cauchy}
    L_k(f) = -\frac{1}{2\pi \iu} \oint_{\Gamma} f(z) \underbrace{\bigl\{\Tr G_k(z) -\Tr G_{k-1}(z) \bigr\}}_{M_k(z)} \dif z,
\end{align}
where $G_k(z) \deq (F_k - z)^{-1}$ is the \emph{resolvent} (or \emph{Green function}) of $F_k$, and the contour $\Gamma$  is taken in the anticlockwise direction enclosing the support $[a,b]$ of the LSD of $F_k$. 

For simplicity, we denote $k_2 \deq k - 1 - k_1$ in the following. The following proposition provides the explicit expressions of the mean and variance of the one-step difference $L_k(f)$ under the null hypothesis $H_0$.
\begin{proposition}\label{pro:stepone_diff_mean_var} 
Suppose that Assumptions \ref{asm:data} and \ref{asm:sample_dimension} hold.
Under $H_0$, for any analytic function $f$, the mean and variance functions of $L_k(f)$ are given by
\begin{align}
    \mu_k(f) \deq \E (L_k(f))
    & = \frac{1}{2 \pi \iu} \oint_{\Gamma} \biggl(f'(z) z \underline{m} - f(z) \frac{\underline{m}'}{\underline{m}} \biggr) \dif z  + \frac{\nu_4-3}{4p \pi \iu} \oint_{\Gamma} f'(z) (1+z\underline{m})^2 \dif z \notag\\
    & \qquad + \frac{1}{2 k_2 \pi \iu} \oint_{\Gamma} \biggl(-\frac{\{zf(z)\}'' z \underline{m}}{2} + z f'(z) \frac{\underline{m}'}{\underline{m}} + \frac{f''(z)}{\underline{m}}\biggr) \dif z + \smallo\biggl(\frac{1}{p}\biggr), \label{eq:LSS_mean}\\
    \sigma_k^2(f) \deq \var (L_k(f)) & = \frac{\nu_4-3}{k_2 c_2} \biggl(\frac{1}{2\pi \iu}\oint_{\Gamma} f'(z) z \underline{m}\dif z\biggr)^2
    + \frac{1}{k_2 \pi \iu} \oint_{\Gamma} \frac{\{f'(z)\}^2}{\underline{m}} \dif z + \smallo\biggl(\frac{1}{p}\biggr).\label{eq:LSS_var}
\end{align}
\end{proposition}
The proof of Proposition \ref{pro:stepone_diff_mean_var} is given in Section \ref{sec:proof_stepone_diff_mean_var}. 
We refer to Section \ref{sec:simulation} for explicit expressions of the mean and variance of $L_k(f)$ for some commonly used test functions.
Based on Proposition \ref{pro:stepone_diff_mean_var}, we first define a random walk type process. For a fixed $T>0$, 
\begin{align*}
    W_{n,t} \deq \frac{1}{\sqrt{n}} \sum_{k=n+1}^{n+\floor{nt}} \tilde{L}_k(f), 
    \quad 0 \leq t \leq T,
\end{align*}
where 
\begin{align}\label{eq:Lk_tilde}  
    \tilde{L}_k(f) \deq \frac{L_k(f) - \mu_k(f)}{\sigma_k(f)},
\end{align}
and $\mu_k(f), \sigma_k(f)$ are defined in Proposition \ref{pro:stepone_diff_mean_var}.
This process is used in defining our monitoring statistics in Section \ref{sec:detect_procedure}. The following theorem establish the weak convergence of the above process under the null hypothesis.
\begin{theorem}[Weak Convergence of $W_{n,t}$ under $H_0$]\label{thm:BM}
    Suppose that Assumptions \ref{asm:data} and \ref{asm:sample_dimension} hold.
    Under the null hypothesis, $\{W_{n,t}, t\in [0,T]\}$ converges in distribution to a standard Brownian motion $\{W(t), t\in [0,T]\}$, i.e., with covariance kernel 
    \[ 
        \cov(W(t), W(s)) = t \wedge s.
    \]
\end{theorem}
The proof of Theorem \ref{thm:BM} is given in Section \ref{sec:proof_thm_BM}.

While the result above characterizes the dynamic behavior of the LLS-based monitoring process, the asymptotic theory of LSS was primarily developed for static hypothesis testing.
In the offline setting, the LSS of large-dimensional sample covariance matrices has served as a cornerstone for hypothesis testing \citep{Bai2004CLT}. 
Based on this theory, an extensive literature has emerged, offering rigorous procedures for testing covariance structures.
These methodologies encompass one-sample tests for identity or sphericity \citep{Bai2009corrections, Wang2013sphericity, Zheng2015substitution, Li2016testing, Liu2023clt,Qiu2023asymptotic, Yin2024liberating, Liu2025asymptotic}, as well as two-sample tests assessing the equality of covariance matrices \citep{Zheng2017CLT, Yang2017weighted, Zhang2019invariant, Zou2021two, Ding2024two}.
However, despite these results in static regimes, extending LSS-based inference to the sequential monitoring domain presents distinct challenges.
We remark here that a stochastic process constructed from the sequential linear spectral statistics of the sequential sample covariance matrices has been considered in \cite{Drnemann2024linear}, for which the limiting process is a non-standard Gaussian process. Here, instead, we construct the process as a random walk with the normalized increments, so that the limiting process is the standard Brownian motion. A major advantage of  such construction is that the local  envelope of the Brownian motion is well-understood, which will be important for the construction of our testing statistics in the sequel.

\subsection{Online change point detection procedure}\label{sec:detect_procedure}

In this section, we introduce our online change point detection procedure, which leverages the asymptotic results of the one-step LSS difference established in Proposition \ref{pro:stepone_diff_mean_var}. 
We assume an initial historical dataset of size $n=k_1+k_2^{\circ}$ is available. 
Specifically, the first $k_1$ observations are utilized to construct the baseline covariance matrix $S_1$, while the subsequent $k_2^{\circ}$ observations form the initial covariance matrix $S_{2,n}$.
The initial sample size $k_2^{\circ}$ is predetermined to ensure that Assumption \ref{asm:sample_dimension} is satisfied at $k=n$. 
Following this initialization, the procedure proceeds to sequentially monitor the incoming data stream $\{y_{n+i}\}_{i \geq 1}$.

Motivated by the online CUSUM procedure proposed in \citet{Chu1996monitoring}, we construct a monitoring statistic at time $n+i$ ($i\geq 1$) based on the cumulative sum of the standardized one-step LSS differences defined in \eqref{eq:Lk_tilde}:
\begin{equation}\label{eq:online_CUSUM}
    T_p(n,i) =  w(n,i) |\Psi(n,i)|,
\end{equation}
where $w(n,i)$ is the weight function satisfying Assumption \ref{asm:weight} below, and $\Psi(n,i)$ is a CUSUM-type statistic defined as
\begin{align}\label{eq:Psi_def}
    \Psi(n,i) = \frac{1}{\sqrt{n}} \sum_{k= n + 1}^{n + i} \tilde{L}_k(f), 
    \quad n = k_1+ k_2^{\circ}, \quad i \geq 1,
\end{align} 
and $\tilde{L}_k(f)$ is defined in \eqref{eq:Lk_tilde}.
At each time point $n+i$ ($i\geq 1$), we update $T_p(n,i)$ by including the new sample $y_{n+i}$.
An alarm is triggered immediately if the monitoring statistic $T_p(n,i)$ exceeds the critical threshold $c_\alpha$; otherwise, the monitoring process continues. 
Hence, the stopping time $\hat{k}^{\star}$ of the procedure is defined by
\begin{align}\label{eq:kstar_hat}
    \hat{k}^{\star} = \inf \bigl\{ i\geq 1: T_p(n,i) > c_{\alpha} \bigr\}
\end{align}
with the convention that $\inf(\varnothing) = \infty$.
The implementation details of the proposed method are summarized in Algorithm~\ref{alg:online_cov_cpt}.

In our online detection procedure, we aim to control the asymptotic false alarm rate at a nominal level $\alpha \in (0,1)$ while ensuring that the detection power converges to one.
Formally, this entails selecting a critical value $c_{\alpha}$ such that
\begin{align*}
    \Prob_{H_0}(\hat{k}^{\star} < \infty) \to \alpha,
    \quad \text{and} \quad
    \Prob_{H_1}(\hat{k}^{\star} < \infty) \to 1,
\end{align*}
where $\Prob_{H_0}$ and $\Prob_{H_1}$ denote the probability under the null and alternative hypothesis, respectively.
To achieve this goal, we need to impose the following regularity conditions on the weight function $w(n,i)$:
\begin{assumption}[Weight function]\label{asm:weight}
    Let the weight function satisfy 
    \begin{enumerate}[label=(\roman*), ref=\theassumption(\roman*)]
        \item \label{asm:w_form} $w(n,i)=\rho(i/n)\one\{i > \ell_n\}$ with $\ell_n/n \to0$, where $\rho: (0,\infty)\to[0,\infty)$ is nonnegative and continuous if restricted to $(0,e_{\rho}]\cap \mathbb{R}$, $e_{\rho}=\sup\{t>0: \rho(t)>0\}$;
        \item \label{asm:w_lim_zero} $\lim_{t\to 0}t^{\gamma} \rho(t)<\infty$ for some $0 \leq \gamma < \frac{1}{2}$;
        \item \label{asm:w_lim_inf} $\lim_{t\to \infty} t \rho(t)<\infty$.
    \end{enumerate}
\end{assumption}
Condition \ref{asm:w_form} mandates a burn-in period (determined by $\ell_n$) prior to the start of monitoring.
This delay is designed to mitigate false alarms during the early stages, 
where the monitoring statistic is prone to instability due to the limited sample size.
Conditions \ref{asm:w_lim_zero} and \ref{asm:w_lim_inf} characterize the weight function's asymptotic behavior to balance sensitivity and robustness: the former restricts the growth rate at $t \to 0$ to mitigate excessive false alarms during the initial monitoring phase, while the latter prevents slow decay at infinity to ensure effective long-term monitoring.

\begin{figure}[htbp]
    \centering
    \begin{subfigure}[b]{0.8\linewidth}
        \centering
        \includegraphics[width=\linewidth]{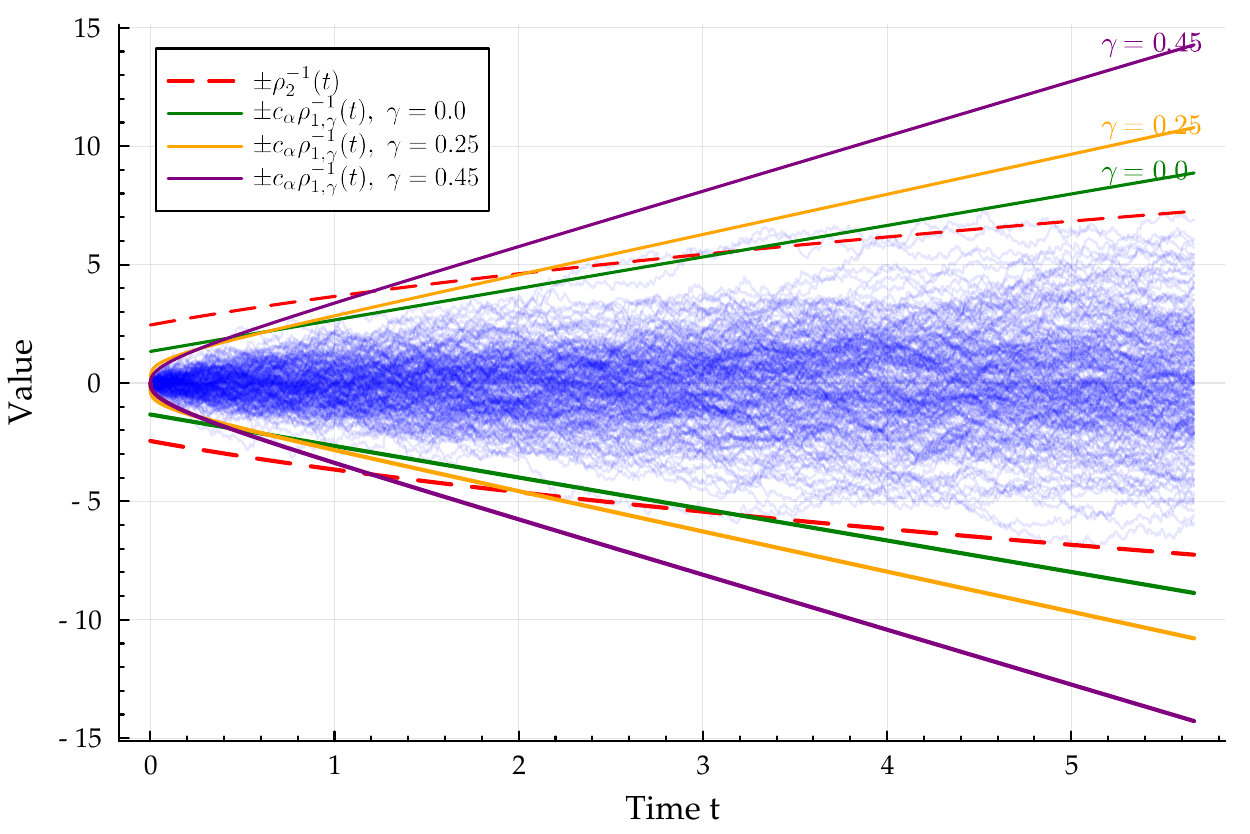}
        \caption{Under $H_0$: $200$ sample paths where the data $\{y_j\}_{j\geq 1}$ are i.i.d. $N(0,I_p)$.}
        \label{fig:trace_diff_BM_H0}
    \end{subfigure}

    \vspace{1em}

    \begin{subfigure}[b]{0.8\linewidth}
        \centering
        \includegraphics[width=\linewidth]{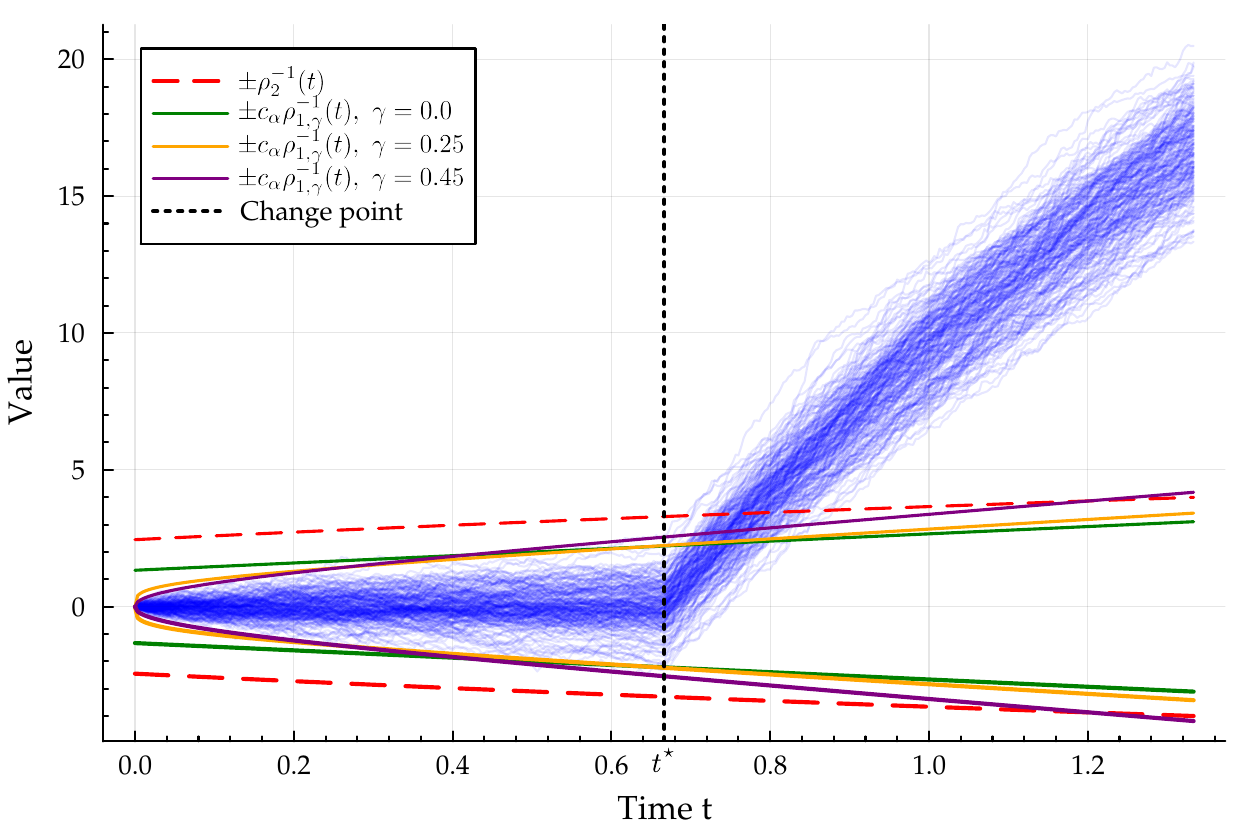}
        \caption{Under $H_1$: $200$ sample paths with a change point at $t^{\star}=2/3$ ($k^{\star}=500$). Data are i.i.d. $N(0,I_p)$ before and $N(0,1.2 I_p)$ after the change.}
        \label{fig:trace_diff_BM_H1}
    \end{subfigure}

    \caption{
        Sample paths (\textcolor{blue}{blue}) of the standardized LSS difference process $\{\frac{1}{\sqrt{n}} \sum_{k=n+1}^{n+\floor{nt}} \tilde{L}_k(f), t\geq 0\}$ with test function $f(x)=x$, together with boundary curves.
        The simulation setting is $p=100$, $n=300$ ($k_1=k_2^{\circ}=150$).
        The \textcolor{red}{red} dashed lines represent the boundary $\pm \rho_2^{-1}(t)$, while the \textcolor{green}{green}, \textcolor{orange}{orange}, and \textcolor{purple}{purple} lines depict the boundaries $\pm c_{\alpha} \rho_{1, \gamma}^{-1}(t)$ for $\gamma = 0.0, 0.25, 0.45$ at significance level $\alpha=0.05$.
        Reject $H_0$ whenever the sample path crosses the selected boundaries.
    }
    \label{fig:trace_diff_BM_combined}
\end{figure}

In this paper, we set the burn-in period as $\ell_n=\log(n)$.
To evaluate the detection performance under different weighting schemes, we employ the following two weight functions:
\begin{align}
    &\rho_{1,\gamma}(t) = (1+t)^{\gamma-1}t^{-\gamma} , \quad 0\leq \gamma < \frac{1}{2}, \label{eq:weight_1}\\
    &\rho_{2}(t) = (1+t)^{-1/2}\bigl\{ -2 \log \alpha + \log(1+t)\bigr\}^{-1/2}.\label{eq:weight_2}
\end{align}
Both the weight family $\rho_{1,\gamma}$ and the logarithmic weight $\rho_2$ are adopted from \citet{Chu1996monitoring}.
To visually demonstrate the behavior of the standardized LSS difference process $\Psi(n,\floor{nt}) = \frac{1}{\sqrt{n}} \sum_{k=n+1}^{n+\floor{nt}} \tilde{L}_k(f)$ and the corresponding detection boundaries implied by these weight functions, we present a simulation example in Figure \ref{fig:trace_diff_BM_combined}.
Under $H_0$, the sample paths fluctuate around zero and remain largely within the prescribed boundaries, while under $H_1$ the process exhibits a clear upward drift after the change point and crosses the boundaries shortly thereafter.

Note that the stopping condition implies $\Prob_{H_0}(\hat{k}^{\star} < \infty) = \Prob_{H_0}(\sup_{i\geq 1} T_p(n,i) > c_{\alpha})$. 
Consequently, to determine the critical value $c_{\alpha}$, it is necessary to derive the limiting null distribution of the process $\{T_p(n,i), i\geq 1\}$. 
The following theorem establishes this asymptotic result.
\begin{theorem}[Limiting null distribution]\label{thm:limit_dist}
    Suppose that Assumptions \ref{asm:data} -- \ref{asm:weight} hold.
    Under the null hypothesis, we have
    \[
        \sup_{i\geq 1} T_p(n,i) \dto  \sup_{t > 0} \rho(t)|W(t)|,
    \]
    where $W(\cdot)$ is a standard Brownian motion.  
\end{theorem}

The proof of Theorem \ref{thm:limit_dist} is given in Section \ref{sec:proof_thm_lim_dist}.

\begin{remark}
    As $t\to 0$, $W(t)=\bigop(\sqrt{t\log\log t^{-1}})$ by the law of iterated logarithm at time $0$. By Assumption \ref{asm:w_lim_zero}, we have $\rho(t)W(t)=\smallop(1)$ for $0\leq \gamma < 1/2$. 
    As $t\to\infty$, the law of iterated logarithm implies $W(t)=\bigop(\sqrt{t\log\log t})$. From Assumption \ref{asm:w_lim_inf}, we have $\rho(t)W(t)=\bigop(\sqrt{t^{-1}\log\log t})=\smallop(1)$. Hence, the supremum in Theorem \ref{thm:limit_dist} can be taken ove $(0,\infty)$. A major reason of adding a weight function to the original process is the non-stationarity of Brownian motion. This non-stationarity introduces bias in the location of the process maxima, even under the null hypothesis, making early change points unlikely to be detected. Applying a weight function to amplify the small-$t$ regime therefore increases sensitivity for detecting early changes. Choosing suitable weight function is possible only if we have a good understanding on the local envelope of the process, which motivated us to consider the normalized one-step LSS difference, whose partial sum process converges to Brownian motion. 
\end{remark}
Based on the asymptotic results in Theorem \ref{thm:limit_dist}, we determine the critical values $c_{\alpha}$ for the weight function $\rho_{1,\gamma}(t)$ via Monte Carlo simulations. The obtained values are summarized in Table \ref{tab:critical_value}.
\begin{table}[htbp]
    \centering
    \caption{Critical value $c_{\alpha}$ for the weighted function $\rho_{1,\gamma}(t)$.}
    \label{tab:critical_value}
    \begin{tabular}{cccccc}
    \toprule
    $\alpha$\verb|\|$\gamma$& 0.00 & 0.15 & 0.25 & 0.35 & 0.45 \\
    \midrule
    0.01 & 1.56949 & 1.81747 & 1.97581 & 2.29276 & 2.78885 \\
    0.05 & 1.33027 & 1.5131 & 1.68472 & 1.93445 & 2.30402 \\
    0.10 & 1.19574 & 1.35757 & 1.50264 & 1.73564 & 2.11163 \\
    \bottomrule
    \end{tabular}
\end{table}

To investigate the asymptotic power and detection delay of our online detection procedure, we  impose the following assumption on the change point $k^{\star}$:
\begin{assumption}[Change point]\label{asm:change_point}
    The change point $k^{\star}$ satisfies that $k^{\star}/n = \bigo(1)$.
\end{assumption}

In practice, if the proposed change point detection procedure indicates that no change has occurred, while the current sample size $i$ becomes large relative to $n$ (i.e., $i/n \gg 1$), we may enlarge $n$ adaptively so that the ratio $i/n$ remains of order $\bigo(1)$.

Recall that the detection time $\hat{k}^{\star}$ is determined by the stopping rule \eqref{eq:kstar_hat}.
The detection delay is formally defined as $d^{\star} \deq \hat{k}^{\star} - k^{\star}$, which quantifies the latency between the true change points and the detection time.
The following theorem characterizes the asymptotic behavior of this delay for our online procedure, along with the asymptotic power.

\begin{theorem}[Detection delay time]\label{thm:power}
    Suppose that Assumptions \ref{asm:data} -- \ref{asm:change_point} hold.
    Let
    \begin{align*} 
        & \tau_1 \deq \lim_{p\to\infty} \frac{1}{p} \Tr(\Sigma_0^{-1}\Sigma_1),
        \quad 
        \tau_2 \deq \lim_{p\to\infty} \frac{1}{p} \Tr(\Sigma_0^{-1}\Sigma_1)^2,\\   
        & I_1(f, \tau_1) \deq \frac{1}{2\pi \iu} \oint_{\Gamma} f'(z) \log \bigl\{\tau_1 - (1-\tau_1)z\underline{m}(z)\bigr\} \dif z,\\
        &I_2(f) \deq \frac{1}{2\pi \iu} \oint_{\Gamma} f'(z) \bigl\{ \log(z\underline{m}) + z \underline{m}\bigr\} \dif z.
    \end{align*}
    We further assume that one of the following conditions holds:
    \begin{enumerate}[label=(\Roman*)]
        \item\label{cond:tau1_non_one} $\tau_1 \neq 1$ and $I_1(f,\tau_1) \neq 0$;
        \item\label{cond:tau1_one} $\tau_1 = 1$, $\tau_2 \neq 1$, and $I_2(f) \neq 0$.
    \end{enumerate}
    The detection delay $d^{\star}$ satisfies one of the following: 
    \begin{enumerate}
        \item (Early-change) If $k^{\star}-n=\smallo(n)$, then the detection delay $d^{\star}$ satisfies $d^{\star} \asymp \log n$ under condition \ref{cond:tau1_non_one}, and $d^{\star} \asymp n^{1/2-\delta}$ for some small $0<\delta<1/2$ under condition \ref{cond:tau1_one}.
        \item (Late-change) If $k^{\star}-n \asymp n$, then the detection delay $d^{\star}$ satisfies $d^{\star} \asymp \log n$ under condition \ref{cond:tau1_non_one}, and $d^{\star} \asymp \sqrt{n}$ under condition \ref{cond:tau1_one}. 
    \end{enumerate}
    
    Furthermore, it holds under $H_1$ that
    \[
        \sup_{i\geq 1} T_p(n,i) \pto +\infty,
    \]
    such that the online test based on the monitoring statistic $T_p(n,i)$ has asymptotic power one.
\end{theorem}
The proof of Theorem \ref{thm:power} is given in Section \ref{sec:proof_thm_power}.

\begin{remark}
    We examine whether the contour integrals $I_1(f,\tau_1)$ and $I_2(f)$
    vanish or not in general. 
    Let 
    \[ 
        L_1(z) \deq \log \bigl\{\tau_1 - (1-\tau_1)z\underline{m}(z)\bigr\} 
        \quad \text{and} \quad
        L_2(z) \deq \log(z\underline{m}) + z \underline{m}.
    \] 
    According to the Residue Theorem, these integrals are non-zero if the sum of the residues of the integrands at all poles in the region exterior to $\Gamma$ is non-zero.    
    We consider two kinds of test functions based on the analytic structure of its derivative $f'(z)$ in the domain exterior to $\Gamma$:
    \begin{itemize}
        \item \emph{Case A: polynomial functions.} Assume $f(z)$ is a polynomial of degree $k$: $f(z) = \sum_{j=0}^{k} c_j z^j$, where $k$ is a finite integer. 
        Its derivative is a polynomial of degree $k-1$:
        $
            f'(z) = \sum_{j=1}^{k} j c_j z^{j-1}
        $
        Since $f'(z)$ is entire, there are no finite poles in the exterior region. The integral is determined solely by the residue at infinity. To compute this, we consider the Laurent expansion of $L_1(z)$ at infinity:
        \[ 
            L_1(z) = \sum_{m=0}^{\infty} \frac{A_m}{z^m} = A_0 + \frac{A_1}{z} + \frac{A_2}{z^2} + \dots + \frac{A_k}{z^k} + \bigo(z^{-(k+1)}).
        \]
        The integrand is the product:
        \[ 
            f'(z) L_1(z) = \left( \sum_{j=1}^k j c_j z^{j-1} \right) \left( \sum_{m=0}^{\infty} \frac{A_m}{z^m} \right).
        \]
        According to the Residue Theorem, we have
        \[ 
            I_1(f,\tau_1) = - \mathrm{Res}(f'(z)L(z), \infty) = \sum_{j=1}^k j c_j A_j.
        \]
        The integral is non-zero provided that the linear combination 
        \[ 
            \sum_{j=1}^k j c_j A_j \neq 0.
        \]
        Specifically, if $f(z) = z$, we need $A_1\neq 0$.
        If $f(z) = z^2$, we need $2 A_2 \neq 0$. 
        Note that 
        $A_1 = \frac{(1-\tau_1)c_2}{1-c_1}$ and
        $A_2 = \frac{c_2(1-\tau_1)\{2+c_2(1-c_1)(1+\tau_1)\}}{2(1-c_1)^3}$.
        Both cases hold under the condition $\tau_1\neq 1$. 
        Hence, $f(x) = x$ and $x^2$ are valid choice to detect the change point when $\tau_1 \neq 1$.
  
        Similarly, for the case $\tau_1 = 1$, we have the Laurent expansion of $L_2(z)$ at infinity:
        \[ 
            L_2(z) = \sum_{m=0}^{\infty} \frac{B_m}{z^m} = B_0 + \frac{B_1}{z} + \frac{B_2}{z^2} + \dots + \frac{B_k}{z^k} + \bigo(z^{-(k+1)}).
        \]
        The integral $I_2(f)$ is non-zero provided that
        \[ 
            \sum_{j=1}^k j c_j B_j \neq 0.
        \]
        Note that $B_1 = 0$ and $B_2 = - \frac{c_2^2}{2(1-c_1)^2}$.
        Thus, $f(x)$ is not a valid choice to detect the change point when $\tau_1 = 1$.
        However, $f(x) = x^2$ is a valid choice since $2 B_2 \neq 0$.

        \item \emph{Case B: function with singularities.} Assume $f(z)$ possesses singularities (e.g., logarithmic functions like $\log(1+z)$). In this scenario, the derivative $f'(z)$ may exhibit isolated poles $z_k$ in the region exterior to $\Gamma$. The integral is non-zero if the sum of residues at these poles is non-zero.
        
        For both integrals, we consider the specific example $f(z) = \log(1+z)$. Using the Residue Theorem, we have
        \[ 
            I_1(f,\tau_1) 
            = - \log \{\tau_1 + (1-\tau_1)m(-1) \}.
        \]
        This integral is non-zero as long as $\underline{m}(-1) \neq 1$. It follows from that $\underline{m}(-1) = \int_a^b \frac{1}{x + 1} \dif F(\lambda)< 1$, where $F(x)$ is the LSD of the companion matrix  $\mathscr{F}_k$. Hence, $f(z) = \log(1+z)$ is a valid choice to detect the change point when $\tau_1 \neq 1$.
        
        Similarly, for the case $\tau_1 = 1$, we have
        \[ 
            I_2(f) 
            = - \bigl\{\log \{ \underline{m}(-1) \} - m(-1) + 1\bigr\}.
        \]
        This integral is non-zero as long as $\underline{m}(-1) \neq 1$, which is verified above. 
        Hence, $f(z) = \log(1+z)$ is still a valid choice to detect the change point when $\tau_1 = 1$.

    \end{itemize}
\end{remark}

\renewcommand{\algorithmicrequire}{\textbf{Input:}}
\renewcommand{\algorithmicensure}{\textbf{Output:}}
\begin{algorithm}[htbp]
\caption{Online covariance change point detection via one-step $F$-matrix LSS differences}
\label{alg:online_cov_cpt}
\begin{algorithmic}[1]
\Require Streaming observations $\{y_k\}_{k\ge 1}$; historical sample sizes $k_1$ and $k_2^{\circ}$; test function $f$; weight function $\rho(\cdot)$; significance level $\alpha$.
\Ensure Detection time $\hat{k}^{\star}$ (or ``no detection'' if the monitoring is stopped externally).

\vspace{0.3em}
\Statex \textbf{Step 1: Critical threshold} 
\State Use Monte Carlo simulations to determine the critical value $c_\alpha$ 
such that $\Prob_{H_0}\bigl(\sup_{t>0} \rho(t)|W(t)| > c_\alpha\bigr) = \alpha$, where $W(t)$ is standard Brownian motion; see Theorem~\ref{thm:limit_dist}.

\vspace{0.3em}
\Statex \textbf{Step 2: Initialization}
\State Set $n \gets k_1 + k_2^{\circ}, \, i \gets 0$.
\State Compute initial sample covariance matrices and $F$-matrix:
\[ 
    S_1 \gets \frac{1}{k_1}\sum_{i=1}^{k_1} y_i y_i^{\top}, 
    \quad
    S_{2,n} \gets \frac{1}{n-k_1}\sum_{i=k_1+1}^{n} y_i y_i^{\top},
    \quad 
    F_n \gets S_1^{-1}S_{2,n}.
\] 
\State Set current time index $k \gets n$.

\vspace{0.3em}
\Statex \textbf{Step 3: Online Monitoring}
\While{new observation $y_{k+1}$ arrives}
    \State Update indices:  $i \gets i +1, \, k \gets k+1$.
    
    \vspace{0.3em}
    \State \textbf{Step 3.1: one-step LSS difference}
    \State Update the second sample covariance matrices $S_{2,k}$ (using recursive update or full sum).
    \State Compute current $F$-matrix: $F_k \gets S_1^{-1}S_{2,k}$.
    \State Compute the one-step LSS difference
    \[
        L_k(f) \gets \Tr f(F_k) - \Tr f(F_{k-1}).
    \]

    \vspace{0.3em}
    \State \textbf{Step3.2: standardization}
    \State Compute the theoretical mean $\mu_k(f)$ and variance $\sigma_k^2(f)$ of $L_k(f)$ using \eqref{eq:LSS_mean} and \eqref{eq:LSS_var}.
    \State Standardize:
    \[
        \tilde{L}_k \gets \frac{L_k(f)-\mu_k(f)}{\sigma_k(f)}.
    \]

    \vspace{0.3em}
    \State \textbf{Step 3.3: monitoring statistic}
    \State Update CUSUM-type statistic (\eqref{eq:online_CUSUM} and \eqref{eq:Psi_def}):
    \[
        \Psi(n,i) \gets \frac{1}{\sqrt{n}}\sum_{t=n+1}^{n+i} \tilde{L}_t,
        \qquad
        T_p(n,i) \gets \rho(i/n)\,|\Psi(n,i)|.
    \]
    \If{$T_p(n,i) > c_\alpha$}
        \State \Return Detection time $\hat{k}^{\star} \gets k$.
    \EndIf
\EndWhile

\State \Return ``no detection''.
\end{algorithmic}
\end{algorithm}

\section{Numerical Simulations}\label{sec:simulation}

\subsection{Simulation Setup}\label{sec:simulation_setting}

We conduct Monte Carlo experiments to examine the finite-sample performance of the proposed sequential detection procedure under both the null ($H_0$) and alternative ($H_1$) hypotheses. 
All reported results are based on $M=2000$ independent replications for each configuration.
The details of the simulation setup are as follows. 

\paragraph{Common Parameters.}
Throughout the simulations, the nominal significance level is set to $\alpha=0.05$.
The weighting functions are chosen as those defined in \eqref{eq:weight_1} and \eqref{eq:weight_2}, where the parameter $\gamma$ in $\rho_{1,\gamma}(t)$ takes values in $\{0.0, 0.25, 0.45\}$.
The critical value $c_{\alpha}$ for $\rho_{1,\gamma}(t)$ with different $\alpha$ and $\gamma$ is obtained through Monte Carlo simulation of the limiting distribution in Theorem \ref{thm:limit_dist}, as summarized in Table \ref{tab:critical_value}.
We set the data dimension to $p=100$. 
The historical data sizes are $k_1=k_2^{\circ}=150$, such that the sequential monitoring starts at time $k_1+k_2^{\circ}+1=301$. 
The data are generated in the form $y_k=\Sigma^{1/2}_k x_k$, where $x_k$ and $\Sigma_k$ satisfy Assumption \ref{asm:data}.
For the estimation of the fourth moment $\nu_4$ of underlying distributions, we adopt the estimator proposed in \citet{Lopes2019bootstrapping}, given by
\[ 
    \hat{\nu}_4 = \max \biggl( 3 + \frac{\hat{\gamma}_n-2\hat{\tau}_n}{\hat{\omega}_n}, \, 1\biggr), 
\]
where 
\begin{align*}
    & \hat{\tau}_n  = \Tr(S_n^2) - \frac{1}{n}\bigl\{\Tr(S_n)\bigr\}^2, 
    \qquad S_n = \frac{1}{n} \sum_{j=1}^n y_j y_j^{\top},\\
    & \hat{\gamma}_n  = \frac{1}{n-1} \sum_{j=1}^n \biggl(\|y_j\|_2^2 - \frac{1}{n}\sum_{j'=1}^n \|y_{j'}\|_2^2\biggr)^2, 
    \qquad
    \hat{\omega}_n^2  = \sum_{i=1}^p \biggl( \frac{1}{n} \sum_{j=1}^n y_{ij}^2 \biggr)^2.
\end{align*}

\paragraph{Empirical Size under $H_0$.}

To evaluate the empirical size, we simulate data under the null hypothesis $H_0$, where no structural change occurs. 
Without loss of generality, we set the true covariance matrix as $\Sigma_k = I_p$ for all $k$.
The standardized observations $x_k$ are generated from three distinct distributions with varying fourth moments $\nu_4$, allowing us to evaluate robustness to tail behavior:
\begin{itemize}
    \item Gaussian distribution $N(0,1)$, with $\nu_4=3$;
    \item Uniform distribution $\mathrm{Unif}(-\sqrt{3},\sqrt{3})$, with $\nu_4=1.8$;
    \item Student's $t$ distribution $t(10)/\sqrt{1.25}$, with $\nu_4=4$.
\end{itemize}
Additionally, we examine four test functions:
\begin{itemize}
    \item Linear: $f(x)=x$;
    \item Logarithmic: $f(x)=\log(1+x)$;
    \item Mixed: $f(x)=x+\log(1+x)$;
    \item Square: $f(x)=x^2$.
\end{itemize}

\paragraph{Empirical Power under $H_1$.}

To assess detection power and detection delay time, we simulate data under $H_1$ with a single change point occurring at $k^{\star}$. 
To examine the impact of the change point location, we vary $k^{\star} \in \{350,450,500\}$.
The underlying distributions and test functions remain consistent with the $H_0$ setting, excluding the \texttt{square} function.
The covariance structure is defined as $\Sigma_k = I_p$ for the pre-change samples ($k \leq k^{\star}$), and switches to $\Sigma_k = \Sigma \neq I_p$ for the post-change samples ($k > k^{\star}$).
Specifically, we investigate three distinct structural changes:
\begin{itemize}
    \item \emph{Homogeneous variance inflation}: $\Sigma=\sigma^2 I$. 
    This scenario simulates a uniform variance change across all dimensions. 
    We vary $\sigma^2\in\{1.1, 1.2, 1.3, 1.4, 1.5\}$;
    \item \emph{Correlation structure change}: 
    The covariance matrix follows a Toeplitz structure with inflated variance, where the diagonal entries are $\Sigma_{jj}=2$ and off-diagonal entries are $\Sigma_{jk}=\rho^{|j-k|}$ for $j\neq k$. 
    This scenario introduces correlations among variables. 
    We vary $\rho\in\{0.1,0.3,0.5,0.7,0.9\}$;
    \item \emph{Heterogeneous variance inflation}: $\Sigma=1.5 I+\delta \sum_{i=1}^5 e_i e_i^\top$, where $e_i$ is the $i$-th standard basis vector in $\mathbb{R}^p$.
    This scenario combines a global variance inflation (from $1$ to $1.5$) with a localized low-rank perturbation, and is designed to test the capability to detect localized large shifts masked by a global variance increase.
    We vary  $\delta\in\{2, 2.5, 3, 3.5, 4\}$. 
\end{itemize}

\paragraph{Explicit forms of $\mu_k(f)$ and $\sigma_k^2(f)$ for selected $f$.}

Now, we present the explicit forms of the mean $\mu_k(f)$ and variance $\sigma_k^2(f)$ of $L_k(f)$ under $H_0$ for several commonly used test functions $f$. 
These explicit expressions facilitate efficient computation of the monitoring statistic $T_p(n,i)$ in practical applications. 
We first define several auxiliary quantities that appear in the expressions:
\begin{align*}
    M_1 & = \frac{c_2}{1-c_1},
    \qquad 
    M_2 = \frac{c_2\bigl(1 + c_2 - c_1 c_2\bigr)}{(1-c_1)^3},\\ 
    M_3 & = \frac{c_2\bigl(c_1^{2}c_2^{2} - 2c_1c_2^{2} - 3c_1c_2 + c_1 + c_2^{2} + 3c_2 + 1\bigr)}{(1-c_1)^5},\\
    M_4 &= \frac{c_2}{(1-c_1)^7}\Big( -c_1^{3}c_2^{3} + 3c_1^{2}c_2^{3} + 6c_1^{2}c_2^{2} - 4c_1^{2}c_2 + c_1^{2} - 3c_1c_2^{3} \\ 
    &\qquad \qquad \qquad \qquad - 12c_1c_2^{2} - 2c_1c_2 + 3c_1  + c_2^{3} + 6c_2^{2} + 6c_2 + 1 \Big), \\
    C_3 & = M_1^4 - 3M_1^2 M_2 + 2M_1 M_3 + M_2^2 - M_4.
\end{align*}

For the selected test functions, the explicit forms of the mean $\mu_k(f)$ of $L_k(f)$ under $H_0$ are as follows:
\begin{align*}
    \mu_k(x) = 0,
    \qquad
    \mu_k(x^2) & = - M_1^2 + \frac{\nu_4-3}{p} M_1^2 + \frac{1}{k_2} M_2,
\end{align*}
\begin{align*}
    &\; \mu_k(\log(1+x)) = \mu_k(x + \log(1+x)) \\ 
    =&\; \{\underline{m}(-1) - 1 - \ln \underline{m}(-1)\} - \frac{\nu_4-3}{2p} \{1 - \underline{m}(-1)\}^2 
    + \frac{1}{k_2} \biggl\{ \frac{1}{2} - \underline{m}'(-1) \left( \frac{1}{2} - \frac{1}{\underline{m}(-1)} + \frac{1}{\underline{m}^2(-1)} \right) \biggr\}.
\end{align*}
The variance $\sigma_k^2(f)$ of $L_k(f)$ under $H_0$ are given by:
\begin{align*}
    \sigma_k^2(x) = \frac{\nu_4-3}{k_2 c_2} M_1^2 - \frac{2}{k_2} (M_1^2 - M_2), 
    \qquad 
    \sigma_k^2(x^2) = \frac{4 M_2^2 (\nu_4-3)}{k_2 c_2} - \frac{8 C_3}{k_2},
\end{align*}
\begin{align*}
    \sigma_k^2(\log(1+z)) = \frac{\nu_4-3}{k_2 c_2} \{\underline{m}(-1)-1\}^2 + \frac{2}{k_2} \left( \frac{\underline{m}'(-1)}{\underline{m}(-1)^2} - 1 \right),
\end{align*}
\begin{align*}
    \sigma_k^2 (z+\log(1+z)) & = \frac{\nu_4-3}{k_2 c_2} \{M_1 + 1 - \underline{m}(-1)\}^2 
    + \frac{2}{k_2} \biggl( M_2 - (M_1-1)^2 +2 - \frac{2}{\underline{m}(-1)} + \frac{\underline{m}'(-1)}{\underline{m}(-1)^2} \biggr).
\end{align*}

The explicit expressions presented above follow from Proposition \ref{pro:stepone_diff_mean_var} and the residue theorem. The derivation details are provided in Section \ref{sec:proof_explicit_mean_var}.
\begin{remark}\label{rmk:integral}
For general test functions $f$, one may convert the corresponding contour integrals in Proposition \ref{pro:stepone_diff_mean_var} into equivalent real integrals, allowing for accurate numerical evaluation of the mean and variance.
For any analytic test function $f$, the mean and variance of $L_k(f)$ under $H_0$ can be computed via the following integrals:
\begin{align*}
    \mu_k(f) & = -\frac{1}{\pi}\int_a^b\biggl( x B f'
    + \frac{A'B-AB'}{A^2+B^2} f\biggr) \dif x - \frac{\nu_4-3}{p\pi} \int_a^b (1+xA) x B f' \dif x \\
    &\qquad - \frac{1}{k_2 \pi} \int_a^b \biggl\{- \frac{(x f'' + 2f') x B}{2} 
    -\frac{(A'B-AB')x f' + B f''}{A^2+B^2}
    \biggr\} \dif x + \smallo\biggl(\frac{1}{p}\biggr),\\
    \sigma_k^2(f) & = \frac{\nu_4-3}{k_2c_2\pi^2}   \biggl(\int_a^b x f' B  \dif x \biggr)^2 
    +\frac{2}{k_2\pi} \int_a^b \frac{B f'^2}{A^2+B^2}\dif x + \smallo\biggl(\frac{1}{p}\biggr),
\end{align*}
where
$A\equiv A(x) \deq -\frac{x(h^2+c_1)+c_2(1-c_2)}{2x(c_2+xc_1)}$, and 
$B\equiv B(x) \deq \frac{c_2(1-c_1)\sqrt{(x-b)(a-x)}}{2x(c_2+xc_1)}$.

The derivation of these integral forms is provided in Section \ref{sec:proof_rmk_integral}.
\end{remark}

\subsection{Simulation Results and Discussion}

With $M=2000$ replications of the dataset simulated under the null/alternative hypothesis, we calculate the empirical size/power and expected detection delay (EDD) as follows:
\begin{align*}
    \text{Empirical Size/Power} = \frac{\sum_{\ell=1}^M \one\{\hat{k}_{\ell}^{\star} < \infty\}}{M}, 
    \qquad 
    \text{EDD} 
    = \frac{\sum_{\ell=1}^M (\hat{k}_{\ell}^{\star} - k^{\star})_{+} }{ \sum_{\ell=1}^M  \one\{\hat{k}_{\ell}^{\star} \ge k^{\star}\}},
\end{align*}
where $\hat{k}_{\ell}^{\star}$ is the stopping time (see \eqref{eq:kstar_hat}) in the $\ell$-th replication.
The former measures the probability of successful detection within the monitoring horizon, 
while the latter quantifies the average number of observations required after the true change point to trigger an alarm. 

Table~\ref{tab:empirical_size} reports the empirical sizes under different combinations of data distributions, weight functions, and test functions at the nominal level $\alpha=0.05$.
The results indicate that the choice of weight and test function impacts the empirical size control.
Specifically, the weight functions $\rho_{1,\gamma}$ with larger values of $\gamma$ tend to produce slightly inflated sizes, particularly when paired with the \texttt{square} test function.
In contrast, the weight $\rho_2$ leads to conservative behavior with empirical sizes frequently below the nominal level. 
Regarding test functions, the \texttt{log}, \texttt{linear} and \texttt{mix} generally provide more reliable size control than the \texttt{square} function, which leads to over-rejection in many settings.
Overall, the combination of the $\rho_{1,0}$ weight with the \texttt{log}, \texttt{linear}, or \texttt{mix} test functions consistently achieves empirical sizes close to the nominal level across all distributions.
The empirical size behavior reported in Table~\ref{tab:empirical_size} can be partly explained by the boundary geometry illustrated in Figure~\ref{fig:trace_diff_BM_combined}.
For the $\rho_{1,\gamma}$ family, increasing $\gamma$ causes the boundary $c_{\alpha}\rho_{1,\gamma}^{-1}$ to narrow around the monitoring process trajectories in the initial phase, making false alarms more probable and leading to slight size inflation.
Conversely, the logarithmic weight $\rho_2$ generates a significantly wider boundary that stays well clear of the trajectories, which decreases spurious detections and yields conservative empirical sizes.

\input{table/table_H0_results}

\begin{figure}[htbp]
    \centering
    \includegraphics[width=\linewidth]{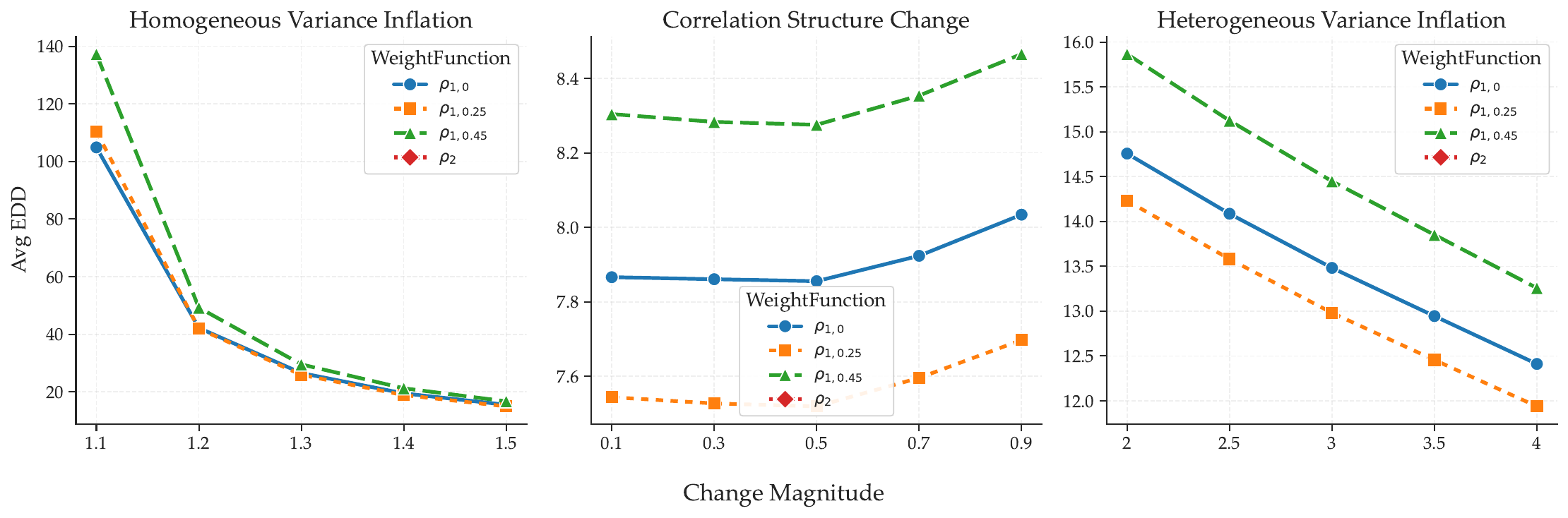}
    \caption{Performance comparison of different weight functions.}
    \label{fig:weight_function_comparison}
\end{figure}

\begin{figure}[htbp]
    \centering
    \includegraphics[width=\linewidth]{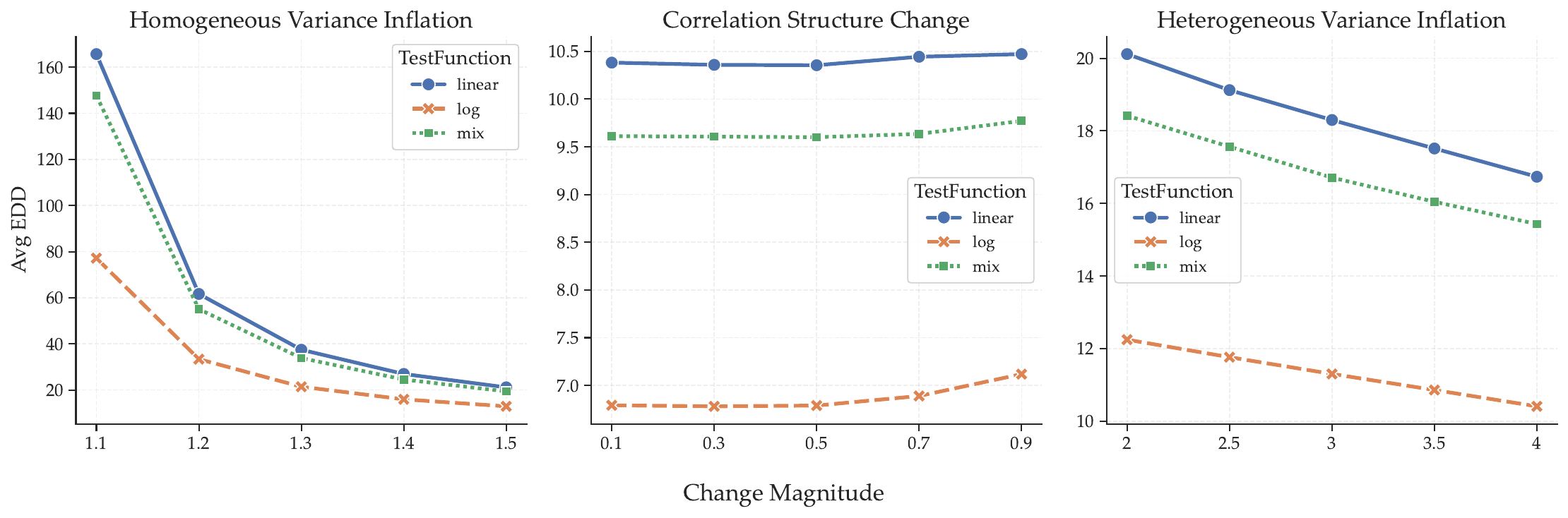}
    \caption{Performance comparison of different test functions.}
    \label{fig:function_comparison}
\end{figure}

\begin{figure}[htbp]
    \centering
    \includegraphics[width=\linewidth]{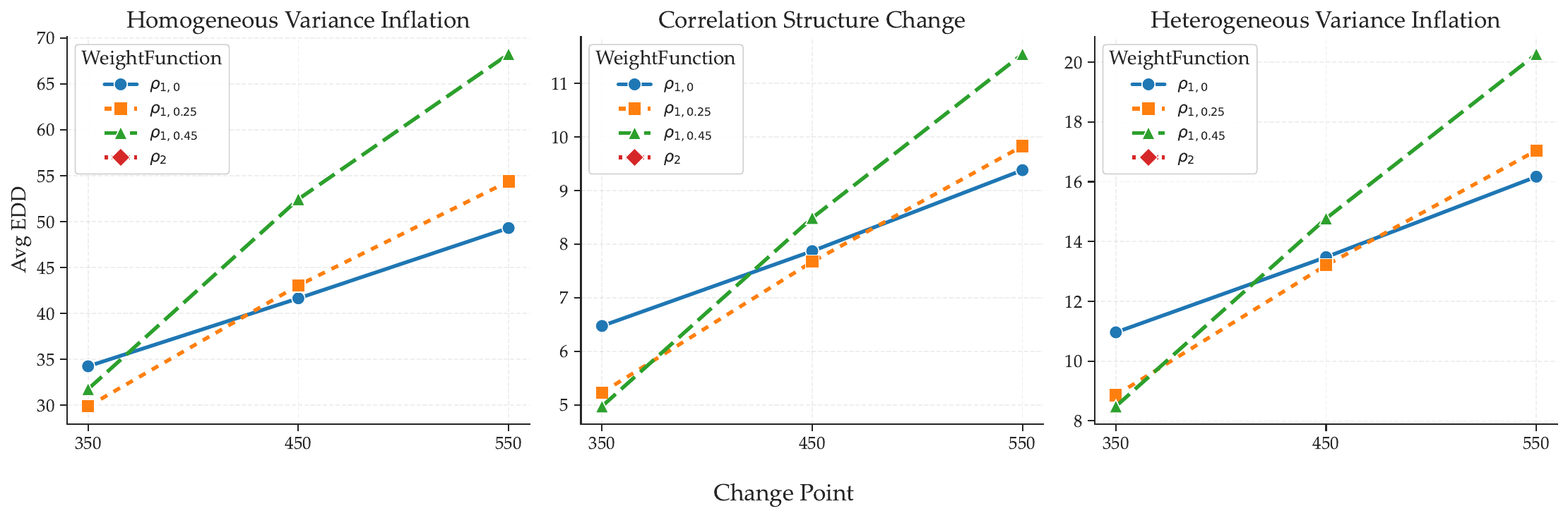}
    \caption{Impact of the change point $k^{\star}$ on detection delay.}
    \label{fig:kstar_sensitivity}
\end{figure}

\begin{figure}[htbp]
    \centering
    \includegraphics[width=\linewidth]{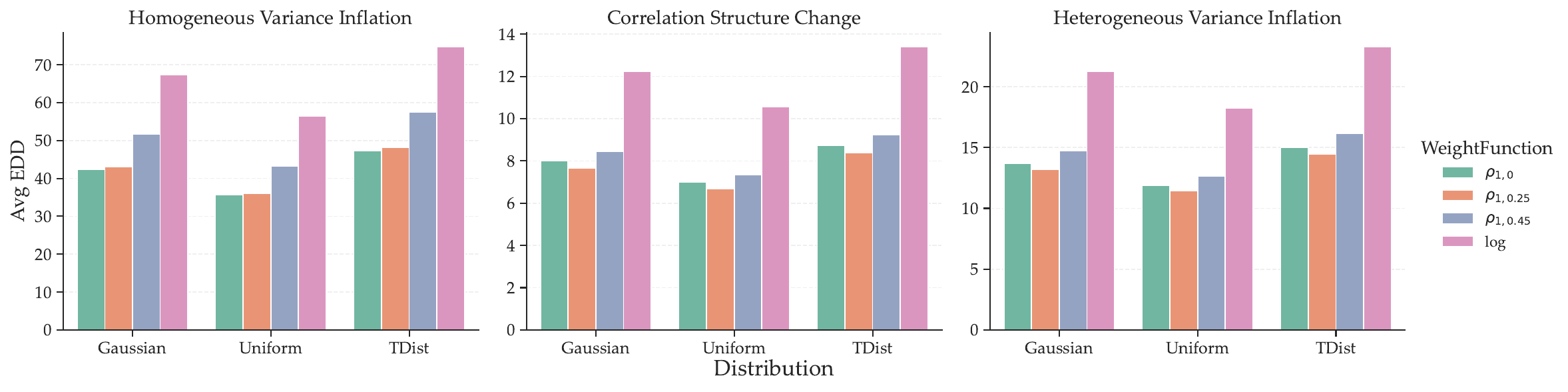}
    \caption{Impact of data distributions on detection delay.}
    \label{fig:distribution_robustness}
\end{figure}

Tables~\ref{tab:inflation_gamma_00} -- \ref{tab:spike_gamma_00} report the empirical power and EDD for weight function $\rho_{1,\gamma}$ with $\gamma=0$. 
To analyze the impact of weight functions, test functions, locations of change point, and data distributions on detection performance, we present Figures~\ref{fig:weight_function_comparison} -- \ref{fig:distribution_robustness}.
These results demonstrate the comprehensive effectiveness and robustness of our proposed detection procedure. We can draw several key conclusions:
\begin{itemize}
    \item (General power and detection delay.) 
    For weak to strong signal strengths, the proposed method exhibits exceptional sensitivity. As illustrated in Tables~\ref{tab:inflation_gamma_00} -- \ref{tab:spike_gamma_00} and Figures~\ref{fig:weight_function_comparison} -- \ref{fig:function_comparison}, it consistently achieves empirical power near unity (1.0 or $\approx 1.0$) with rapidly decreasing EDD across nearly all scenarios.
    Theoretical insights from Theorem \ref{thm:power} reveal that the detection delay primarily rely on the quantity $\tau_1$.
    In the scenario of correlation structure changes, $\tau_1$ remains invariant (constant at $2$) with respect to $\rho$, resulting in a relatively stable EDD across varying correlation magnitudes. 
    In contrast, for homogeneous/heterogeneous variance inflation scenarios, $\tau_1$ changes with signal strength, leading to a substantial reduction in EDD as the signal magnitude increases.  
    
    \item (Superiority of the \texttt{log} test function.) 
    As evidenced Figure~\ref{fig:function_comparison}, the \texttt{log} test function consistently achieves the lowest EDD across all change scenarios. 
    Notably, in the homogeneous variance inflation scenario (see Table \ref{tab:inflation_gamma_00}), it exhibits superior sensitivity, yielding higher empirical power than both the \texttt{linear} and \texttt{mixed} functions in the weak signal regime where $\sigma^2=1.1$.
    
    \item (Impact of change point location.) 
    Figure~\ref{fig:kstar_sensitivity} illustrates the impact of the change point location on detection performance.
    We observe consistent trends across diverse weight functions and change scenarios:  earlier changes lead to shorter detection delays.
    Comparing the weighting schemes, the $\rho_{1,\gamma}$ family generally yields shorter detection delays than the $\rho_2$ weight. 
    Within the $\rho_{1,\gamma}$ class, the choice of $\gamma$ involves a trade-off: larger $\gamma$ values enhance sensitivity to early changes, whereas smaller $\gamma$ values prove more effective for changes occurring later in the sequence. 
    Among the candidates considered, the $\rho_{1,0}$ weight is the most robust option and shows relatively little dependence on the change point location.

    \item (Robustness to data distribution.) 
    A primary strength of our method is its remarkable robustness to the underlying data generating process.
    As evidenced in Figure \ref{fig:distribution_robustness}, the detection performance remains virtually invariant across Gaussian, Uniform, and Student's $t$ distributions. 
    This strongly suggests that the efficacy of our method are insensitive to the underlying data distribution.
    This property constitutes a significant practical advantage, particularly in applications involving high-dimensional data that may exhibit heavy tails or deviate from normality.
\end{itemize}

\input{table/table_H1_results_gamma_0}

In summary, based on the comprehensive numerical simulations, we recommend the \texttt{log} test function paired with the $\rho_{1,\gamma}$ weight function with $\gamma = 0$ as the default configuration for practical applications. 
The \texttt{log} function demonstrates superior sensitivity across all scenarios, particularly in the challenging weak signal regime. 
Regarding the weighting scheme, $\rho_{1,0}$ offers the most favorable trade-off, ensuring rapid detection while maintaining stability invariant to the unknown change point location. 
This combination provides a robust and powerful solution that performs reliably across diverse data distributions.

\subsection{Comparison with Existing Methods}
We conduct simulation studies to benchmark our proposed method against two alternative procedures: \citet{Avanesov2019structural} (referred to as \textsf{A19}) and \citet{Li2023online} (referred to as \textsf{LL23}). 
The simulation setup fix the dimension at $p=50$ and the true change point at $k^{\star}=200$. The observations are generated from a standard normal distribution, with the covariance structure shifting from $\Sigma_0 = I_p$ to $\Sigma_1$ according to the three scenarios described in Section \ref{sec:simulation_setting}.

For our proposed method, we use the weight function $\rho_{1,0}$ and the test function $f(x)=\log(1+x)$. 
The historical data sizes are set to $k_1=60$ and $k_2^{\circ}=100$, with a nominal significance level of $\alpha=0.05$. 
For the \textsf{A19} method, we specific a window size of $20$ and a significance level of $\alpha=0.05$.
The detection threshold is determined via the bootstrap procedure described in Section 2.2 of \citet{Avanesov2019structural}, using $500$ bootstrap replications.
Regarding the \textsf{LL23} method, we set the window size parameter to $100$.
To ensure a fair comparison at the same false alarm rate, we calibrate the average run length (ARL) of \textsf{LL23} to $18492$. 
This ARL value corresponds to a threshold of $a=3.95$ \citep[see Section 4.4]{Li2023online} in \textsf{LL23}'s stopping rule and their Theorem 1.

\input{table/table_comparison}

All results, summarized in Table \ref{tab:comparison}, are based on $500$ independent Monte Carlo replications for each configuration.
Our method outperforms the competitors in both accuracy and speed, maintaining high empirical power ($>0.99$) across all change scenarios, even in the weakest signal scenarios.
Conversely, \textsf{LL23} suffers from a significant loss of power under small variance inflation (e.g., Power is only $0.332$ at $\sigma^2=1.1$), although its performance improves as the signal strength increases. \textsf{A19} also displays reduced power in this setting ($0.794$).
Regarding detection efficiency, our approach achieves shorter detection delays compared to both benchmark procedures.
This advantage is particularly pronounced in weak signal settings; for example, our method reduces the detection delay by approximately $35\%$ to $80\%$ compared to competitors across the smallest change magnitudes ($\sigma^2=1.1$, $\rho=0.1$, and $\delta=2$), verifying its superior sensitivity to minor structural changes.

\section{A Real Data Example}\label{sec:real_data}

In this section, we demonstrate the practical applicability of the proposed online covariance change point detection method through an analysis of the S\&P 500 stock data.
Specifically, the covariance structure characterizes the dynamic comovement and risk integration among assets; detecting its abrupt changes is crucial for identifying market regime shifts and adjusting risk management strategies accordingly.
The dataset consists of historical daily adjusted closing prices for the constituents of the S\&P 500 index, obtained from Yahoo Finance\footnote{\url{https://finance.yahoo.com/}}.
We utilized adjusted closing prices to account for corporate actions such as dividends and stock splits.
To ensure data quality, we excluded stocks with a missing data rate exceeding $5\%$ and imputed minor gaps in the price series using the forward-filling method. 
After preprocessing, we obtain daily adjusted closing prices for $496$ stocks spanning $756$ trading days.
Daily log-returns are computed as $y_k = \log(P_{k}/P_{k-1})$, where $P_k$ denotes the adjusted closing price on day $k$.
To mitigate the impact of extreme outliers or data errors, we applied a Winsorization technique to the computed log-returns, clamping values exceeding 5 standard deviations from the mean.
To capture the most significant market fluctuations, we select a subset of the top $30$ stocks with the highest volatility (defined as the standard deviation of daily log-returns) from the pool of $497$ stocks.

\begin{figure}[!htbp]
    \centering
    \includegraphics[width=0.8\linewidth]{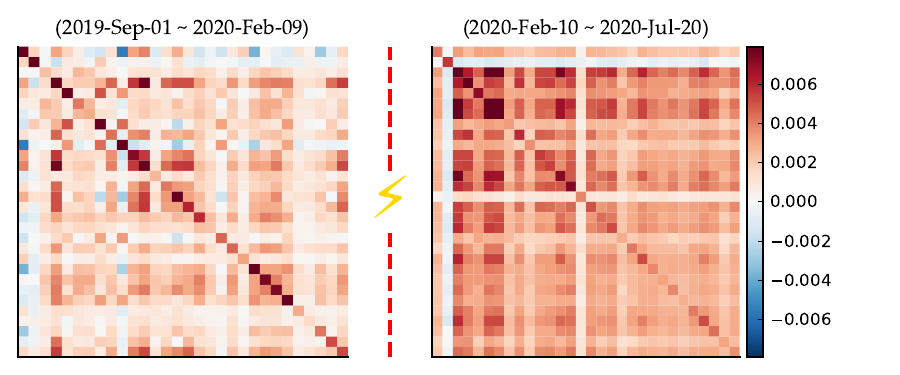}
    \caption{Heatmaps of the sample covariance matrices for the top 30 S\&P 500 stocks. The left panel represents the covariance structure before the detected change point ($\hat{k}^{\star} =$ February 10, 2020), while the right panel corresponds to the post-change period. The increased color intensity in the right panel illustrates the surge in volatility and correlation triggered by the onset of the COVID-19 market crash.}
    \label{fig:real_data_sp500}
\end{figure}

We apply the proposed online monitoring procedure (Algorithm \ref{alg:online_cov_cpt}) to the log-return series, utilizing the $\log(1+x)$ test function and the $\rho_{1,0}$ weight function. 
The historical window sizes are set to $k_1=k_2^{\circ}=40$, and the monitoring phase commences at $k=81$ (corresponding to December 26, 2019).
The nominal significance level is set to $\alpha=0.05$.
The procedure detects a change point on February 10, 2020. 
This detection serves as a early warning signal preceding the major COVID-19 market crash. 
It captures the initial structural change in the covariance matrix as investors began pricing in the systemic risk of the coronavirus outbreak, notably weeks before the widespread panic selling occurred in late February.
Figure \ref{fig:real_data_sp500} compares the sample covariance matrices before and after this detected change.
It is evident that the covariance structure undergoes a drastic shift, with the post-change period characterized by universally larger entries corresponding to the onset of widespread market stress.

\section{Proof of Theorem \ref{thm:BM}}\label{sec:proof_thm_BM}
To prove the weak convergence of the process $\{W_{n,t}, t\in [0,T]\}$, we only need to show to the \emph{finite dimensional distribution convergence} and the \emph{tightness} of $\{W_{n,t}, t\in [0,T]\}$ (see e.g. \citet{billingsley1968convergence}). Therefore, Theorem \ref{thm:BM} follows immediately from the following two lemmas.

\begin{lemma}[Finite-dimensional distribution]\label{lem:finite-dim-dist}
    Suppose that Assumptions \ref{asm:data} and \ref{asm:sample_dimension} hold.
    Under the null hypothesis, for any fixed $r\in \mathbb{N}$ and $t_1,t_2,\cdots,t_r \in [0, T]$, the random vector $(W_{n,t_1}, \ldots, W_{n,t_r})$ converges in distribution to a $r$-dimensional Gaussian distribution with mean zero and covariance matrix $\Sigma = (\sigma_{ij})_{1\leq i,j \leq r}$, where $\sigma_{ij} = t_i \wedge t_j$.
\end{lemma}

\begin{lemma}[Tightness]\label{lem:tightness}
    The sequence $\{W_{n,t}, t \in [0,T]\}$ is asymptotically tight in the space $\ell^{\infty}([0,T])$. 
\end{lemma}

The proof of Lemma \ref{lem:finite-dim-dist} and Lemma \ref{lem:tightness} are given in Section \ref{sec:lem_FDD} and Section \ref{sec:lem_tight} respectively.

\subsection{Proof of Lemma \ref{lem:finite-dim-dist}}\label{sec:lem_FDD}
Recall that 
\begin{align*}
    W_{n,t} \deq \frac{1}{\sqrt{n}} \sum_{k=n+1}^{n+\floor{nt}} \tilde{L}_k(f), 
    \quad 0 \leq t \leq T.
\end{align*}
By the Cram\'er-Wold device, it suffices to show that, for any fixed $a_1, a_2 \in \bbR$ and $t_1, t_2 \in [0,T]$,
\begin{align*}
	a_1W_{n,t_1} + a_2W_{n,t_2} \dto N(0, \tau^2_{1,2}),
\end{align*}
where $\tau_{1,2}^2 \deq a_1^2 t_1 + a_2^2 t_2 + 2 a_1 a_2 (t_1\wedge t_2)$, and the extension to more points is straightforward.
WLOG, we assume that $T\geq t_1 > t_2 \geq 0$. Note that
\begin{align*}
    a_1 W_{n,t_1} + a_2 W_{n,t_2} = \sum_{k=n+1}^{n+\floor{nt_1}} b_k (L_k-\mu_k), 
\end{align*}
where 
\[
b_k = \begin{cases}
    \frac{a_1+a_2}{\sqrt{n} \sigma_k}, & n+1 \leq k \leq n+\floor{nt_2}, \\
    \frac{a_1}{\sqrt{n} \sigma_k}, & n+\floor{nt_2}+1 \leq k \leq n+\floor{nt_1}.
\end{cases}
\]
Let $\calF_k$ be the $\sigma$-algebra generated by $x_{n+1}, \ldots, x_{n+k}$.
Denote $\E_k \deq \E(\cdot | \calF_k)$.
From the definition of $L_k$, we have
\begin{align*}
L_k - \mu_k & = \frac{1}{2\pi \iu} \oint_{\Gamma} f(z) (1-\E) Y_k(z) \dif z + \bigo(1/p),
\end{align*}
where 
\[ 
    Y_k(z) \deq (\E_k - \E_{k-1}) \frac{x_k^{\top}A_k(z) x_k}{k_2}
\] 
forms a martingale difference sequence, and $A_k(z)$ is defined as
\begin{align*}
    A_k(z) \deq \alpha(z) G^2 S_1^{-1} + \alpha'(z)  G S_1^{-1},
    \quad 
    \alpha(z) \deq \frac{ 1 }{1+\frac{1}{k_2} \Tr (GS_1^{-1}) }. 
\end{align*}
The detailed calculation about the mean structure $L_k$ can be found in Section \ref{sec:proof_lem_resolvent_mean_var}. Hence, we have
\begin{align*}
    a_1 W_{n,t_1} + a_2 W_{n,t_2} 
    = \frac{1-\E}{2\pi \iu} \oint_{\Gamma} f(z) S_{m,t_1,t_2}(z) \dif z, 
    \qquad 
    S_{m,t_1,t_2}(z) \deq \sum_{k=n+1}^{n+\floor{nt_1}} b_k Y_k(z). 
\end{align*}
By functional CLT and Lemma \ref{lem:resolvent_mean_var}, it suffices to show that $\{S_{m,t_1,t_2}(z)\}$ converges weakly to a Gaussian process $S(z)$ with covariance function
\begin{align}\label{eq:cov_Sz1_Sz2}
    \cov(S(z_1), S(z_2)) = \sum_{k=n+1}^{n+\floor{nt_1}} b_k^2 \biggl\{ \frac{\nu_4-3}{k_2c_2} \partial_{z_1} \bigl(z_1\underline{m}_1' \bigr) \partial_{z_2} \bigl(z_2\underline{m}_2' \bigr) - \frac{2}{k_2} \partial_{z_1}\partial_{z_2} \biggl( \frac{1/\underline{m}_1-1/\underline{m}_2}{z_1-z_2} \biggr) \biggr\}. 
\end{align}
Using this formula, we conclude that $a_1 W_{n,t_1} + a_2 W_{n,t_2}$ converges in distribution to a centered normal distribution with variance
\begin{align*} 
    \var(a_1 W_{n,t_1} + a_2 W_{n,t_2}) = \sum_{k=n+1}^{n+\floor{nt_2}} \frac{(a_1+a_2)^2}{m} + \sum_{k=n+\floor{nt_2}+1}^{n+\floor{nt_1}} \frac{a_1^2}{m} = \tau_{1,2}^2.
\end{align*}

Now, prove the weak convergence of $\{S_{m,t_1,t_2}(z)\}$.
First, we derive the finite-dimensional distribution of $S_{m,t_1,t_2}(z)$.
From the martingale CLT it suffices to verify that
\begin{align*}
    &\sum_{k=n+1}^{n+\floor{nt_1}} \E \bigl\{b_k^2 Y_k(z_1)Y_k(z_2) | \calF_{k-1}\bigr\} \pto \eqref{eq:cov_Sz1_Sz2}, \\
    &\sum_{k=n+1}^{n+\floor{nt_1}} \E \bigl\{b_k^2 Y_k^2(z) \one_{\{|b_k Y_k(z)| > \epsilon\}}\bigr\} \to 0, \quad \forall \epsilon > 0.
\end{align*}
The first condition can be verified by using Lemmas \ref{lem:resolvent_mean_var} and \ref{lem:trace_limits}.
We consider the second condition.
For any small $\delta > 0$, using Lemma B.26 in \cite{bai2010spectral}, we have the following bound
\begin{align*}
    &\;\sum_{k=n+1}^{n+\floor{nt_1}} \E \bigl\{b_k^2 Y_k^2(z) \one_{\{|b_k Y_k(z)| > \epsilon\}}\bigr\} \leq \frac{1}{\varepsilon^{\delta/2}}\sum_{k=n+1}^{n+\floor{nt_1}} b_k^{2+\delta/2} \E |Y_k|^{2+\delta/2}  \\
    \lesssim &\;  \frac{\max_{k}b_k^{2+\delta/2}}{k_2^{2+\delta/2} \varepsilon^{\delta/2}} \sum_{k=n+1}^{n+\floor{nt_1}} \Bigl[\E \bigl\{ \nu_4  \Tr \bigl(A_k(z)A_k(\bar{z})\bigr)\bigr\}^{1+\delta/4} + \nu_{4+\delta} \E \Tr \bigl\{A_k(z)A_k(\bar{z})\bigr\}^{1+\delta/4} \Bigr] \\ 
    \lesssim &\; \frac{1}{k^{2+\delta/2}} \sum_{k=n+1}^{n+\floor{nt_1}} p^{1+\delta/4} = \bigo(p^{-\delta/4}).
\end{align*}
where $\nu_{\ell}\deq \E x_{11}^{\ell}$, and in the last step we used the fact that $\max_k b_k = \bigo(1)$ and Lemma \ref{lem:trace_limits}. 

Then, we show the tightness (in $z$) of the process $\{S_{m,t_1,t_2}(z)\}$. 
It is sufficient to verify the moment condition (see \citet{billingsley1968convergence} (12.51) for instance), i.e.,
\begin{align*}
    \sup_{z_1,z_2\in\Gamma} \frac{\E|S_{m,t_1,t_2}(z_1)-S_{m,t_1,t_2}(z_2)|^2}{|z_1-z_2|^2} < \infty.
\end{align*}
Using the definition of $S_{m,t_1,t_2}(z)$, for any $z_1,z_2\in\Gamma$, we have
\begin{align*}
    \E|S_{m,t_1,t_2}(z_1)-S_{m,t_1,t_2}(z_2)|^2 = \sum_{k=n+1}^{n+\floor{nt_1}} b_k^2 \E |Y_k(z_1)-Y_k(z_2)|^2. 
\end{align*}
For any $k$, $Y_k(z)$ is analytic in $z$. Thus, we have
\[ 
    Y_k(z_1) - Y_k(z_2) = (z_1-z_2) \int_0^1 Y_k'\bigl(z_2 + t(z_1-z_2)\bigr) \dif t,
\]
then 
\[ 
    \E |Y_k(z_1)-Y_k(z_2)|^2 \le |z_1-z_2|^2 \sup_{z\in\Gamma} \E |Y_k'(z)|^2.
\]
Plugging into the required moment condition, we have
\[ 
    \sup_{z_1,z_2\in\Gamma} \frac{\E|S_{m,t_1,t_2}(z_1)-S_{m,t_1,t_2}(z_2)|^2}{|z_1-z_2|^2} \le \biggl(\sum_k b_k^2 \biggr) \sup_k\sup_{z\in\Gamma} \E |Y_k'(z)|^2.
\]
From the definition of $b_k$, it is easy to see that $\sum_k b_k^2 = \bigo(p)$. It remains to show that 
\begin{align}\label{eq:Yk_derivative}
    \sup_k\sup_{z\in\Gamma} \E |Y_k'(z)|^2=\bigo(1/p).
\end{align}
From the definition of $Y_k(z)$ , we have
\[ 
    \E |Y_k'(z)|^2 
    = \frac{1}{k_2^2}\E \Big| x_k^{\top} A_k'(z) x_k - \Tr A_k'(z) \Big|^2
    \lesssim \frac{1}{p^2} \E \Tr \bigl\{A_k'(z) A_k'(\bar{z})\bigr\} \leq \frac{1}{p} \E \|A'(z)\|_2^2,
\]
and the derivative of $A_k(z)$ is given by
\[ 
    A_k'(z) = \Bigl\{\alpha'(z) G_{k-1}^2(z)  + 2\alpha(z) G_{k-1}^3(z)  + \alpha''(z) G(z)  +\alpha'(z) G_{k-1}^2(z) \Bigr\}S_1^{-1}.
\]
Since $y_1 < 1$, $\|S_1^{-1}\|_2$ is bounded with high probability. By Lemma \ref{lem:trace_limits} and the fact that $\|G_{k-1}(z)\|_2 \le \eta^{-1}$ with $\eta = \Im z$, we have
\[ 
    \|A'(z)\|_2 = \bigop(1).
\]
The above estimates imply \eqref{eq:Yk_derivative}, which completes the proof of tightness of $\{S_{m,t_1,t_2}(z)\}$.

\subsection{Proof of Lemma \ref{lem:tightness}}\label{sec:lem_tight}
Without loss of generality, we further assume $t\in [0,1]$. 
Consider a partition $[0,1]=\cup_{j=1}^n K_j$, where
\begin{equation*}
    K_j=\left[ \frac{j-1}{n},\frac{j}{n} \right], \quad n\in \mathbb{N}^+
\end{equation*}
From \cite[Theorem 1.5.6]{van1996weak}, it suffices to show that for any $\lambda >0$,
\begin{equation}\label{eq:tightness_goal}
    \limsup _{n\to \infty}\mathbb{P}\left( \sup_{1\leq j \leq n}\sup_{t_1,t_2\in K_j}|W_{n,t_1}-W_{n,t_2}|>\lambda \right) = 0.
\end{equation}

Note that $\{\tilde{L}_k(f)\}_{k\in \intset{n+1, 2n}}$ forms a martingale difference sequence. According to Lemma \ref{lem:DFJ_inequality}, for some small $\varepsilon > 0$, we have the following bound on the increment of $W_{n,t}$:
\begin{equation*}
    \E |W_{n,t}-W_{n,s}|^{2+\varepsilon} \lesssim (mt-ms)^{1+\varepsilon/2} \frac{\E |\tilde{L}_k|^{2+\varepsilon}}{m^{1+\varepsilon/2}} = (t-s)^{1+\varepsilon/2} \E |\tilde{L}_k|^{2+\varepsilon}.
\end{equation*}
From the definition of $\tilde{L}_k$, the mean structure of $M_k(z)$ given in \eqref{eq:Mk}, and Lemma \ref{lem:resolvent_mean_var}, we have
\begin{align*}
    \E |\tilde{L}_k|^{2+\varepsilon} & \asymp p^{1+\varepsilon/2} \E \biggl| \oint_{\Gamma} f(z) (1-\E_k)M_k(z) \dif z\biggr|^{2+\varepsilon}\\
    & \lesssim p^{1+\varepsilon/2}  \int_{\Gamma} |f(z)|^{2+\varepsilon} \E\bigl|(1-\E_k)M_k(z)\bigr|^{2+\varepsilon} |\mathrm{d}z| \\ 
    & \lesssim p^{-1-\varepsilon/2} \E \bigl| x_k^{\top} A_k(z) x_k -\Tr A_k(z) \bigr|^{2+\varepsilon}.
\end{align*}
By Assumption \ref{asm:data}, we have
\begin{equation*}
    \E \bigl|x_k^{\top} A_k x_k - \Tr A_k\bigr|^{2+\varepsilon} \lesssim  \E |x_{11}|^{4+2\varepsilon} \Tr \{(A_k A_k^*)^{1+\varepsilon/2}\} \lesssim p.
\end{equation*}
Using Markov's inequality, it holds for all $\lambda>0$ that
\begin{equation}\label{09041432}
    \Prob(|W_{n,t}-W_{n,s}| >\lambda)\lesssim (t-s)^{1+\varepsilon/2} p^{-\varepsilon/2}.
\end{equation}
Similarly, we have for $0\leq r\leq s\leq t\leq 1$,
\begin{align}\label{09041433}
     \; \Prob\bigl(\min \{|W_{n,s}-W_{n,t}|,|W_{n,s}-W_{n,r}|\}>\lambda\bigr)
    &\leq \; \mathbb{P}(|W_{n,s}-W_{n,t}| >\lambda)+\mathbb{P}(|W_{n,s}-W_{n,r}| >\lambda) \notag\\
   & \lesssim \; (t-r)^{1+\varepsilon/2} p^{-\varepsilon/2}.
\end{align}
Combining  (\ref{09041432}) and (\ref{09041433}) with Lemma \ref{lem:increment_control}, we get for any interval $K_j$ and any $\lambda >0$,
\begin{equation*}
    \Prob \biggl( \sup_{t_1,t_2\in K_j}|W_{n,t_1}-W_{n,t_2}| >\lambda \biggr)
    \lesssim \frac{1}{n^{1+\varepsilon/2}},
\end{equation*}
which implies \eqref{eq:tightness_goal} and thus completes the proof.

\section{Proof of Theorem \ref{thm:limit_dist}}\label{sec:proof_thm_lim_dist}
We start with the proof for the closed-end situation with $\rho(t)=0$ for all $t > T$. 
From Theorem \ref{thm:BM}, for any $0 < \xi <T$, we have
\begin{align}
    \sup_{\xi n \leq i \leq nT} T_p(n,i) & = \sup_{\xi n < i \leq n T} \rho \biggl(\frac{i}{n}\biggr) \Biggl|\frac{1}{\sqrt{n}}\sum_{k=n+1}^{n+i} \tilde{L}_k \Biggr| \notag \\
    &= \sup_{\xi < t \leq T} \rho \biggl(\frac{\floor{nt}}{n}\biggr) \Biggl|\frac{1}{\sqrt{n}}\sum_{k=n+1}^{n+\floor{nt}} \tilde{L}_k \Biggr| 
    \dto \sup_{\xi < t \leq T} \rho(t) |W(t)|. \label{eq:limit_CUSUM_large_part}
\end{align}
Let $t = i/n\in(0,\xi)$, and $w(n,i)=\rho(i/n)=\rho(t)$. For any fixed $\tilde{\gamma}>\gamma$, we have 
\[
    \rho(t) = t^{-\tilde{\gamma}} \cdot t^{\tilde{\gamma}} \rho(t)
    \le t^{-\tilde{\gamma}} \sup_{0< x < \xi} x^{\tilde{\gamma}} \rho(x) = \frac{n^{\tilde{\gamma}}}{i^{\tilde{\gamma}}}\sup_{0<x<\xi} x^{\tilde{\gamma}} \rho(x).
\]
Hence, 
\begin{align*}
    \sup_{1\leq i < \xi n} T_p(n,i) & \leq \sup_{0 < x < \xi} x^{\tilde{\gamma}} \rho(x) \cdot \Biggl( \sup_{1\leq i < n} \frac{1}{n^{1/2-\tilde{\gamma}}i^{\tilde{\gamma}}} \Biggl|\sum_{k=n+1}^{n+i} \tilde{L}_k \Biggr| \Biggr).
\end{align*}
Under Assumption \ref{asm:w_lim_zero}, we have 
\[ 
    \sup_{0 < x <\xi} x^{\tilde{\gamma}} \rho(x) 
    = \sup_{0 < x <\xi} x^{\tilde{\gamma}-\gamma} x^{\gamma}\rho(x) 
    \lesssim \xi^{\tilde{\gamma}-\gamma} = \smallo(1)
    \quad \text{as } \xi \to 0.
\]
From the H\'{a}jek-R\'{e}nyi inequality (Lemma \ref{lem:hajek-renyi}), we have, for any $\varepsilon > 0$,
\[
    \Prob\Biggl(\sup_{1\leq i < n} \frac{1}{n^{1/2-\tilde{\gamma}}i^{\tilde{\gamma}}} \Biggl|\sum_{k=n+1}^{n+i} \tilde{L}_k \Biggr| \geq \varepsilon \Biggr) \leq \frac{1}{\varepsilon^2} \frac{1}{n^{1-2\tilde{\gamma}}} \sum_{i=1}^{n-1} i^{-2\tilde{\gamma}} \asymp \frac{1}{\varepsilon^2},
\]
and thus 
\[ 
    \sup_{1\leq i < n} \frac{1}{n^{1/2-\tilde{\gamma}}i^{\tilde{\gamma}}} \Biggl|\sum_{k=n+1}^{n+i} \tilde{L}_k \Biggr| = \bigop(1).
\]
From the above results, it holds for $\xi \to 0$ that
\begin{align}\label{eq:limit_CUSUM_small_part}
    \sup_{1\leq i < \xi n} T_p(n,i) = \smallop(1).
\end{align}
Recall that a standard Brownian motion satisfies $|W(t)| = \bigop(\sqrt{t})$ as $t \to 0$. By Assumption \ref{asm:w_lim_zero}, we have $\rho(t) = \bigo(t^{-\gamma})$ as $t\to 0$ for some $0 \le \gamma < 1/2$.
Hence, it holds that
\begin{align}\label{eq:limit_BM_small_part}
    \sup_{0 < t <\xi} \rho(t) |W(t)| =  \sup_{0 < t <\xi} \bigop(t^{1/2-\gamma}) = \smallop(1), 
    \quad \text{as } \xi \to 0. 
\end{align}
Combining \eqref{eq:limit_CUSUM_large_part}, \eqref{eq:limit_CUSUM_small_part} and \eqref{eq:limit_BM_small_part}, we conclude that for any $T > 0$,
\begin{align}\label{eq:limit_CUSUM_closed_end}
    \sup_{1 \leq i \leq nT} T_p(n,i) \dto  \sup_{0 < t \leq T} \rho(t) |W(t)|.
\end{align}

Now, we consider the open-end situation with $\rho(t) > 0$ for all $t > 0$.
Using the similar arguments above, for any $T > 0$, we have
\begin{align}\label{eq:limit_CUSUM_finite_N}
    \frac{1}{\sqrt{n}}
    \sup_{i\geq 1} \rho\biggl(\frac{i}{n}\biggr)\biggl| \sum_{k=n+1}^{n+\min(i,nT)} \tilde{L}_k \biggr| \dto  \sup_{t>0} \rho(t) \bigl|W(\min(t,T))\bigr|.
\end{align}
For any $i > nT$, we let $t=i/n \geq T$ and get $\rho(i/n)\leq \frac{m}{i}\sup_{x \geq N} x\rho(x)$.
This, together with the triangular inequality, yields that
\begin{align}
    &\;\Biggl|\sup_{i \geq 1} T_p(n,i) - \frac{1}{\sqrt{n}}
    \sup_{i\geq 1} \rho(i/n)\biggl| \sum_{k=n+1}^{n+\min(i,nT)} \tilde{L}_k \biggr| \Biggr| \notag\\
    \leq &\; \frac{1}{\sqrt{n}} \sup_{i \geq nT} \rho(i/n) \Biggl|\sum_{k=n+1}^{n+i} \tilde{L}_k - \sum_{k=n+1}^{n+\min(i,nT)} \tilde{L}_k \Biggr| \notag\\
    \leq &\; \sup_{t \geq T} t\rho(t) \cdot \sup_{i\geq nT} \frac{\sqrt{n}}{i} \Biggl|\sum_{k=n+nT+1}^{n+i} \tilde{L}_k \Biggr|\notag\\ 
    = &\; \bigop(T^{-1/2}) \cdot \sup_{t \geq T} t\rho(t) \notag\\ 
    = &\; \bigop(T^{-1/2}), \qquad (\text{Assumption \ref{asm:w_lim_inf}})
    \label{eq:infty_diff}
\end{align}
where the $\bigop(T^{-1/2})$ term follows from 
that, for any $\varepsilon > 0$,
\begin{align*}
    \Prob\biggl( \sup_{i\geq nT} \frac{\sqrt{n}}{i} \biggl|\sum_{k=n+nT+1}^{n+i} \tilde{L}_k \biggr| \geq \varepsilon \biggr) \leq \frac{1}{\varepsilon^2} \sum_{i=nT}^{\infty} \frac{n}{i^2} \leq \frac{1}{\varepsilon^2 T}.
\end{align*}
This bound is obtained by applying the H\'{a}jek-R\'{e}nyi inequality.
From Assumption \ref{asm:w_lim_inf} and the property of Brownian motion, we have
\begin{align}\label{eq:infty_BM_diff}
    \biggl| \, \sup_{t >0} \rho(t) \bigl|W(\min(t,T))\bigr| - \sup_{t>0} \rho(t) \bigl|W(t)\bigr|\, \biggr|
    \le \sup_{t > T} \rho(t) \bigl|W(T)-W(t)\bigr| = \bigop(T^{-1/2}).    
\end{align} 
From \eqref{eq:limit_CUSUM_closed_end}, \eqref{eq:limit_CUSUM_finite_N}, \eqref{eq:infty_diff} and \eqref{eq:infty_BM_diff}, we conclude that
\[ 
    \sup_{i\geq 1} T_p(n,i) \dto  \sup_{t > 0} \rho(t)|W(t)|.
\]

\section{Proof of Theorem \ref{thm:power}}\label{sec:proof_thm_power}
Let $\ell \deq n+i-k^{\star}$ be the number of post-change samples included in the current monitoring procedure at time $n+i$,
and assume that $\ell = \smallo(n)$. 
Define $d_k \deq L_k(f) - \mu_k$, and let
\[
    \Phi_n(\ell) := C_{\sigma}\Bigl|\rho_{1,\gamma}(i/n)\sum_{k=k^{\star}+1}^{k^{\star}+\ell}\E d_k\Bigr|,
    \qquad
    C_{\sigma} := \min_{k^{\star}+1\le k\le k^{\star}+\ell}\frac{1}{\sqrt{n}\,\sigma_k},
\]
\[
    R_n(\ell) := \Bigl|\frac{\rho_{1,\gamma}(i/n)}{\sqrt{n}}\sum_{k=n+1}^{k^{\star}+\ell}\tilde L_k\Bigr|.
\]
From Proposition \ref{pro:stepone_diff_mean_var}, we have $C_{\sigma} = \bigo(1)$.
By the triangle inequality, 
\[
    T_p(n,i) \geqslant \Phi_n(\ell) - R_n(\ell).
\]
We claim that (to be proven later), for $\ell=o(n)$, the following estimate holds:
\begin{align}\label{eq:mean_drift}
    \sum_{k=k^{\star}+1}^{k^{\star}+\ell} \E d_k = \biggl(I_1(f,\tau_1) \ell + \frac{(\tau_2-1) I_2(f)}{k_2^{\star}} \ell^2 \biggr) \bigl(1+\smallo(1)\bigr),
\end{align}
where $k_2^{\star} \deq k^{\star} - k_1$. 
It then follows that
\begin{align*}
    T_p(n,i) & \geqslant C_{\sigma} \rho_{1,\gamma}(i/n) \biggl|I_1(f,\tau_1) \ell + \frac{(\tau_2-1)I_2(f)}{k_2^{\star}}\ell^2\biggr|\bigl(1+\smallo(1)\bigr) - R_n(\ell). 
\end{align*}
Since $\{\tilde{L}_k, n+1 \leq k \leq k^{\star}+\ell\}$ is a martingale difference sequence, 
\[  
    R_n(\ell) = (k^{\star}-n+\ell)^{-\gamma}n^{\gamma-1/2}\bigop(\max\{\sqrt{\ell}, \sqrt{k^{\star}-n}\}).
\]

We first consider the late-change case, where $k^{\star} - n \asymp n$.
In this case,
\[ 
    \rho_{1,\gamma}\biggl(\frac{k^{\star}-n+\ell}{n}\biggr) = \bigo(1), 
    \quad 
    R_n(\ell) = \bigop(1).
\]  
Hence, for any $\varepsilon > 0$, there exists constants $K_{\varepsilon} > 0$ and $N_{\varepsilon}$ such that for all $n > N_{\varepsilon}$,
$ 
    \Prob\bigl(R_n(\ell) \leq K_{\varepsilon}\bigr) > 1 - \varepsilon.
$
Choose $\ell$ such that
\[
    C_{\sigma} \rho_{1,\gamma}\biggl(\frac{k^{\star}-n+\ell}{n}\biggr) \biggl|I_1(f,\tau_1) \ell + \frac{(\tau_2-1)I_2(f)}{k_2^{\star}}\ell^2\biggr| > c_{\alpha} + K_{\varepsilon}, 
\] 
then
\[
    \Prob(T_p(n,i) > c_{\alpha}) \geq \Prob\bigl(R_n(\ell) \leq K_{\varepsilon}\bigr) > 1 - \varepsilon.
\]
Taking $\varepsilon \to 0$ yields  
\[ 
    \Prob(T_p(n,i)>c_{\alpha}) \to 1.
\]
If $\tau_1 \neq 1$ and $I_1(f,\tau_1) \neq 0$, it suffices to choose $\ell \asymp \log(n)$.
If $\tau_1 = 1$ but $\tau_2 \neq 1$, and $I_2(f)\neq 0$, we can choose $\ell \asymp \sqrt{n}$ to obtain the same conclusion.

Then, we consider the early change case with $k^{\star} - n = \smallo(n)$. 
In this case, we have
\[ 
    \rho_{1,\gamma}\biggl(\frac{k^{\star}-n+\ell}{n}\biggr) \asymp n^{\gamma}(k^{\star}-n+\ell)^{-\gamma},
    \quad 
    R_n(\ell) = \smallop(1).
\]
Thus, the mean drift term is of order
\[ 
    \Phi_n(\ell) \asymp n^{\gamma}(k^{\star}-n+\ell)^{-\gamma} \biggl|I_1(f,\tau_1) \ell + \frac{(\tau_2-1)I_2(f)}{k_2^{\star}}\ell^2\biggr|.
\]
If $\tau_1 \neq 1$ and $I_1(f,\tau_1)\neq 0$, it suffices to choose $\ell\asymp \log(n)$. If $\tau_1 = 1$ but $\tau_2 \neq 1$ and $I_2(f)\neq 0$, we can choose $\ell \asymp n^{1/2-\delta}$ for some small $0<\delta<1/2$.


It remains to prove \eqref{eq:mean_drift}. WLOG, we assume that $\Sigma_0 = I_p$. Under $H_1$, there exist a change point $k^{\star}$ such that the covariance of $x_i$ changes from $\Sigma_0=I$ to some $\Sigma_1 \neq I_p$. 
To inspect the impact of the change of $\Sigma_1$ on the linear spectral statistics $L_k$ for $k > k^{\star}$, we examine the sequential $F$-matrix constructed after the change point, up to time $k^{\star}+\ell$. 
Under this setting, the sample covariance matrix $S_{2,k^{\star}+\ell}$ is defined as
\[
	S_{2,k^{\star}+\ell} = \frac{1}{k-k_1}\left( \sum_{i = k_1+1}^{k^{\star}}x_ix_i^{\top} + \sum_{i=k^{\star}+1}^{k^{\star}+\ell} \Sigma_1^{1/2}x_ix_i^{\top}\Sigma_1^{1/2} \right).
\]
By the Cauchy integral formula, it suffices to estimate the order of 
\[ 
    \E \Tr\bigl(G_{k^{\star}+\ell}(z)  - G_{k^{\star}}(z)\bigr).
\] 
For notational simplicity, let $k_2^{\star} \deq k^{\star}-k_1$, $F_{\star} \deq F_{k^{\star}}$, and $G_{\star}(z) \deq G_{k^{\star}}(z)$. 
Let 
\[ 
    Y \deq (y_{k^{\star}+1}, \ldots, y_{k^{\star}+\ell}), 
    \quad \text{where} \quad 
    y_i = \Sigma_1^{1/2}x_i, \quad i > k^{\star}.
\]
Using the Woodbury identity, we have
\begin{align}\label{eq:mixed_resolvent}
    G_{k^{\star}+\ell} 
    &= \left(\frac{k_2^{\star}}{k_2^{\star}+\ell} F_{\star} + \frac{1}{k_2^{\star}+\ell} S_1^{-1} YY^{\top} - z\right)^{-1} \notag \\
    &= \hat{G}_{\star} - \frac{1}{k_2^{\star}+\ell}\hat{G}_{\star}S_1^{-1}Y\left(I_{\ell}+\frac{1}{k_2^{\star}+\ell}Y^{\top}\hat{G}_{\star}S_1^{-1}Y\right)^{-1}Y^{\top}\hat{G}_{\star},
\end{align}
where $\hat{G}_{\star} \deq \bigl(\frac{k_2^{\star}}{k_2^{\star}+\ell}F_{\star}-z\bigr)^{-1}$. 
Applying the identity $A^{-1}-B^{-1}=A^{-1}(B-A)B^{-1}$ repeatedly, we get 
\begin{align}\label{exp}
    \hat{G}_{\star} 
    = \frac{k_2^{\star}+\ell}{k_2^{\star}} \biggl\{ G_{\star} + \frac{\ell z}{k_2^{\star}} G_{\star}^2 + \biggl(\frac{\ell z}{k_2^{\star}}\biggr)^2 G_{\star}^3 \biggr\} + R_1,
\end{align}
where $R_1$ is the remainder matrix with $\|R_1\|_2 = \bigop(\ell^3 p^{-3})$. For the inverse matrix involved in \eqref{eq:mixed_resolvent}, since $\|Y^{\top}\hat{G}_{\star}S_1^{-1}Y\|_2 \asymp \|Y^{\top}Y\|_2 = \bigop(\ell)$, by the Neumann series expansion, we can write 
\begin{align}\label{Neumann_exp}
    &\; \biggl(I_{\ell}+\frac{Y^{\top}\hat{G}_{\star}S_1^{-1}Y}{k_2^{\star}+\ell}\biggr)^{-1} 
    = \biggl\{ I_{\ell} + \frac{ \E \Tr(\hat{G}_{\star}S_1^{-1}\Sigma_1)}{k_2^{\star}+\ell}I_{\ell} + \frac{(1-\E)Y^{\top} \hat{G}_{\star}S_1^{-1}Y}{k_2^{\star}+\ell}\biggr\}^{-1} \notag \\
    =&\; \biggl(1+\frac{\E \Tr(\hat{G}_{\star}S_1^{-1}\Sigma_1)}{k_2^{\star}+\ell}\biggr)^{-1}I_{\ell} - \frac{\frac{1}{k_2^{\star}+\ell}(1-\E)Y^{\top}\hat{G}_{\star}S_1^{-1}Y}{\bigl\{1+\frac{1}{k_2^{\star}+\ell}\E\Tr(\hat{G}_{\star}S_1^{-1}\Sigma_1)\bigr\}^2}  + R_2,
\end{align}
where $R_2$ is the remainder matrix with $\|R_2\|_2 = \bigop(\ell^2 p^{-2})$. Plugging \eqref{exp} and \eqref{Neumann_exp} into \eqref{eq:mixed_resolvent}, we have 
\begin{align}
 	&\;\E \Tr \bigl(G_{k^{\star}+\ell} -  G_{\star} \bigr) \notag \\
    =& \; \frac{\ell}{k_2^{\star}} \E \Tr G_{\star} 
    + \frac{z\ell(k_2^{\star}+\ell)}{(k_2^{\star})^2} \E \Tr G_{\star}^2 
    + \frac{(z\ell)^2(k_2^{\star}+\ell)}{(k_2^{\star})^3} \E \Tr G_{\star}^3 
    - \frac{\frac{1}{k_2^{\star}+\ell}\E\Tr\bigl(Y^{\top}\hat{G}_{\star}^2S_1^{-1}Y\bigr)}{1+\frac{1}{k_2^{\star}+\ell}\E \Tr(\hat{G}_{\star}S_1^{-1}\Sigma_1)} \notag\\
    &\qquad +\frac{\frac{1}{(k_2^{\star}+\ell)^2}\Tr\bigl\{Y^{\top}\hat{G}_{\star}^2S_1^{-1}Y (1-\E )Y^{\top}\hat{G}_{\star}S_1^{-1}Y\bigr\}}{\bigl\{1+\frac{1}{k_2^{\star}}\E \Tr(G_{\star}S_1^{-1}\Sigma_1)\bigr\}^2} + \bigo(\ell^3 p^{-2}) \notag\\
 	= & \; \frac{\ell}{k_2^{\star}}\E \Tr G_{\star} 
    + \frac{z\ell}{k_2^{\star}} \E \Tr G_{\star}^2 
    + \frac{(z\ell)^2 }{(k_2^{\star})^2} \E \Tr G_{\star}^3 
    - \frac{\frac{\ell}{k_2^{\star}}\E \Tr(\tilde{\calG}\Sigma_1)}{1+\frac{1}{k_2^{\star}}\E \Tr( \calG \Sigma_1)} 
    -\frac{\ell^2}{k_2^{\star}}z\left(\frac{\frac{1}{k_2^{\star}}\E \Tr(\tilde{\mathcal{G}}\Sigma_1)}{1+\frac{1}{k_2^{\star}}\E \Tr(\mathcal{G}\Sigma_1)}\right)^2 \notag\\
 	&\qquad + \frac{\ell^2}{k_2^{\star}}\frac{\frac{1}{k_2^{\star}}\E \Tr(G^3S_1^{-1}\Sigma_1)}{1+\frac{1}{k^{\star}}\E \Tr(\mathcal{G}\Sigma_1)} + \frac{\frac{1}{(k_2^{\star})^2}\E \Tr\bigl\{(1-\E)Y^{\top}\tilde{\calG} Y (1-\E )Y^{\top} \calG_{\star} Y\bigr\}}{\bigl\{1+\frac{1}{k_2^{\star}}\E \Tr(\calG\Sigma_1)\bigr\}^2} 
    + \bigo(\ell^3 p^{-2}),\label{eq:tr_diff_G_Gstar}
\end{align}
where $\calG(z) \deq G_{\star}(z)S_1^{-1}$ and $\tilde{\calG}(z) \deq G_{\star}^2 (z)S_1^{-1}$. Note that $\tilde{\calG}(z) = \partial_z \calG(z)$.
It remains to estimate 
\[  
    \frac{1}{k_2^{\star}} \E \Tr(\calG\Sigma_1)
    \qquad \text{and} \qquad
    \frac{1}{k_2^{\star}} \E \Tr\Bigl\{(1-\E)Y^{\top} \tilde{\calG} Y (1-\E )Y^{\top}\calG Y\Bigr\}.
\] 
The limit of the first term is given in Lemma \ref{lem:trace_limit_Sigma} below. For the second term, we have
\begin{align*}
    \;& \sum_{i=1}^{\ell} \E \Bigl\{(1-\E )Y^{\top}\tilde{\calG} Y (1-\E )Y^{\top} \calG Y\Bigr\}_{ii} \notag \\
 	=&\; \sum_{i,j=1}^{\ell}\E  \Bigl\{y_i^{\top}\tilde{\calG} y_j - \Tr(\tilde{\calG} \Sigma_1)\one_{\{i=j\}}\Bigr\} \Bigl\{ y_j^{\top}\calG y_i - \Tr(\calG \Sigma_1)\one_{\{i=j\}}\Bigr\} \\
	=&\; (\nu_4-3)\ell \sum_{k=1}^p \E \bigl(\Sigma_1^{1/2}\tilde{\calG} \Sigma_1^{1/2}\bigr)_{kk}\bigl(\Sigma_1^{1/2}\calG \Sigma_1^{1/2}\bigr)_{kk} 
    + 2\ell \E \Tr\bigl(\tilde{\calG} \Sigma_1 \calG \Sigma_1\bigr) 
	+ \ell(\ell-1) \E  \bigl(y_1^{\top}\tilde{\calG} y_2 y_1^{\top}\calG y_2 \bigr) \\
	=&\; (\nu_4-3) \ell \sum_{k=1}^p \E \bigl(\Sigma_1^{1/2}\tilde{\calG} \Sigma_1^{1/2}\bigr)_{kk}\bigl(\Sigma_1^{1/2}\calG \Sigma_1^{1/2}\bigr)_{kk} 
    + \ell(\ell+1) \E \Tr\bigl(\tilde{\calG} \Sigma_1\calG \Sigma_1\bigr).
\end{align*}
Using the basic inequality $|A_{kk}| \leq \|A\|_2$, the first term in the above equation is of order $\bigo(\ell p)$. The second term is of larger order $\bigo(\ell^2 p)$. Therefore, the contribution of the first term is negligible compared to the second term, and we focus on the latter, whose limit is also given in Lemma \ref{lem:trace_limit_Sigma} below. 

\begin{lemma}\label{lem:trace_limit_Sigma}
Suppose Assumptions \ref{asm:data} and \ref{asm:sample_dimension} hold, then we have
\begin{align*}
    \frac{1}{k_2^{\star}} \E \Tr (\calG \Sigma_1) & = \tau_1 \mathtt{g}(z)+ \bigo(1/p),\\
	\frac{1}{k_2^{\star}} \E \Tr\bigl( \tilde{\calG}\Sigma_1 \calG \Sigma_1\bigr) 
	& = (\tau_2 - \tau_1^2 ) \mathtt{g}'(z)\mathtt{g}(z) +\tau_1^2 \mathtt{h}(z)
     + \bigo(1/p),
\end{align*}
where 
\begin{align*}
    \tau_1 \deq \lim_{p\to\infty} \frac{1}{p} \Tr(\Sigma_1), 
    \quad 
    \tau_2 \deq \lim_{p\to\infty} \frac{1}{p} \Tr(\Sigma_1^2),
    \quad 
    \mathtt{g}(z) \deq \lim_{p\to\infty} \frac{1}{k_2^{\star}} \E \Tr\calG,
    \quad 
    \mathtt{h}(z) \deq \lim_{p\to\infty} \frac{1}{k_2^{\star}} \E \Tr(\tilde{\calG}{\calG}).
\end{align*}
\end{lemma}
The proof of Lemma \ref{lem:trace_limit_Sigma} is postponed to Section \ref{sec:proof_lem_trace_limit_Sigma}.
Using this lemma and \eqref{eq:tr_diff_G_Gstar}, we have
\begin{align}
    \sum_{k=k^{\star}+1}^{k^{\star}+\ell} \E d_k 
    & = -\sum_{k=k^{\star}+1}^{k^{\star}+\ell} \frac{1}{2\pi \iu}\oint_{\Gamma} f(z) \E \Tr\bigl(G_{k^{\star}+\ell}-G_{\star}\bigr) \dif z + \sum_{k=k^{\star}+1}^{k^{\star}+\ell} \mu_k \\ 
    & = - \frac{\ell}{2\pi \iu} \oint_{\Gamma} f(z) \Biggl(\frac{\frac{1}{k_2^{\star}}\E \partial_z \{\Tr (\calG \Sigma_1)\}}{1+\frac{1}{k_2^{\star}}\E \Tr( \calG\Sigma_1)} - \frac{\frac{1}{k_2^{\star}}\E \partial_z(\Tr \calG)}{1+\frac{1}{k_2^{\star}}\E \Tr \calG}\Biggr) \dif z \notag \\ 
    & \qquad +\frac{\ell(\ell+1)}{2 k_2^{\star} \pi \iu} \oint_{\Gamma}  f(z) \Biggl(\frac{\frac{1}{k_2^{\star}}\E \Tr\bigl( \tilde{\calG}\Sigma_1 \calG \Sigma_1\bigr)}{\bigl\{1+\frac{1}{k_2^{\star}}\E \Tr(\calG\Sigma_1)\bigr\}^2} - \frac{\frac{1}{k_2^{\star}}\E \Tr\bigl( \tilde{\calG}\calG \bigr)}{\bigl\{1+\frac{1}{k_2^{\star}}\E \Tr \calG\bigr\}^2}  \Biggr)\dif z \notag \\
    &\qquad +\frac{\ell^2}{k_2^{\star}} \oint_{\Gamma} f(z) z\left[\left(\frac{\frac{1}{k_2^{\star}}\E \Tr(\tilde{\mathcal{G}}\Sigma_1)}{1+\frac{1}{k_2^{\star}}\E \Tr(\mathcal{G}\Sigma_1)}\right)^2 - \left(\frac{\frac{1}{k_2^{\star}}\E \Tr \tilde{\mathcal{G}} }{1+\frac{1}{k_2^{\star}}\E \Tr \mathcal{G}}\right)^2\right] \dif z \notag \\
    &\qquad - \frac{\ell^2}{k_2^{\star}}\oint_{\Gamma} f(z) \left(\frac{\frac{1}{k_2^{\star}}\E \Tr(G^3S_1^{-1}\Sigma_1)}{1+\frac{1}{k_2^{\star}}\E \Tr(\mathcal{G}\Sigma_1)} - \frac{\frac{1}{k_2^{\star}}\E \Tr(G^3S_1^{-1})}{1+\frac{1}{k_2^{\star}}\E \Tr \mathcal{G}}\right)\dif z + \bigo(\ell/p) \notag \\
    & = - \frac{\ell}{2\pi \iu} \oint_{\Gamma} f(z) \partial_z \log\biggl(\frac{1+\tau_1 \mathtt{g}(z)}{1+\mathtt{g}(z)}\biggr) \dif z \notag \\
    & \qquad + \frac{\ell(\ell+1)}{2 k_2^{\star} \pi \iu} \oint_{\Gamma} f(z) \Bigg(\frac{(\tau_2 - \tau_1^2 ) \mathtt{g}'(z)\mathtt{g}(z) +\tau_1^2 \mathtt{h}(z)}{\bigl\{1+\tau_1 \mathtt{g}(z)\bigr\}^2} - \frac{\mathtt{h}(z)}{\bigl\{1+\mathtt{g}(z)\bigr\}^2} \Bigg) \dif z \notag \\
    &\qquad +\frac{\ell^2}{k_2^{\star}} \oint_{\Gamma} f(z) z \left[\left(\frac{\tau_1\mathtt{g}(z)}{1+\tau_1\mathtt{g} (z)}\right)^2 - \left(\frac{\mathtt{g}(z)}{1+\mathtt{g}(z)}\right)^2\right] \dif z \notag \\
    &\qquad - \frac{\ell^2}{k_2^{\star}}\oint_{\Gamma} f(z) \left(\frac{\tau_1\mathtt{g}''(z)/2}{1+\tau_1\mathtt{g}(z)} - \frac{\mathtt{g}''(z)/2}{1+\mathtt{g}(z)}\right) \dif z
    + \bigo(\ell/p). \label{eq:power_mean_decomp}
\end{align}

First, we consider the case $\tau_1 \neq 1$.
From the integration by parts and $\mathtt{g}(z) = -\frac{1+z\underline{m}}{z\underline{m}}$, we have
\begin{align*}
    &\;\frac{1}{2\pi \iu}\oint_{\Gamma} f(z) \partial_z \log\biggl(\frac{1+\tau_1 \mathtt{g}(z)}{1+\mathtt{g}(z)}\biggr) \dif z 
    = - \frac{1}{2\pi \iu}\oint_{\Gamma} f'(z) \log\biggl(\frac{1+\tau_1 \mathtt{g}(z)}{1+\mathtt{g}(z)}\biggr) \dif z \\
    = &\; - \frac{1}{2\pi \iu}\oint_{\Gamma} f'(z) \log \bigl\{\tau_1 - (1-\tau_1)z\underline{m}(z)\bigr\} \dif z = - I_1(f,\tau_1).
\end{align*}
Under the condition $\tau_1\neq 1$, for selected test function $f$ such that $I_1(f,\tau_1) \neq 0$, the first term in \eqref{eq:power_mean_decomp} is of order $\bigo(\ell)$, and it dominates the remaining terms, which is of order $\bigo(\ell^2/p)$. 

Next, we consider the case $\tau_1 = 1$ and $\tau_2\neq 1$. 
In this case, the first, third and fourth terms in \eqref{eq:power_mean_decomp} vanishes, and we have
\begin{align}\label{eq:power_mean_decomp_tau1_eq1}
    \E \sum_{k=k^{\star}+1}^{k^{\star}+\ell} d_k 
    = \frac{\ell(\ell+1)(\tau_2 - 1)}{k_2^{\star}} 
    \frac{1}{2\pi \iu} \oint_{\Gamma} \frac{f(z)\mathtt{g}'(z)\mathtt{g}(z)}{\bigl\{1+\mathtt{g}(z)\bigr\}^2} 
    \dif z + \bigo(\ell/p).
\end{align}
Using $\mathtt{g}(z) = -\frac{1+z\underline{m}}{z\underline{m}}$, we have
\begin{align*} 
    &\; \frac{1}{2\pi \iu}\oint_{\Gamma} \frac{f(z)\mathtt{g}'(z)\mathtt{g}(z)}{\bigl\{1+\mathtt{g}(z)\bigr\}^2} 
    \dif z
    = 
    - \frac{1}{2\pi \iu}\oint_{\Gamma} f(z)(z\underline{m})' \biggl(\frac{1}{z\underline{m}(z)} + 1\biggr) \dif z \\
    =&\; 
    \frac{1}{2\pi \iu}\oint_{\Gamma} f'(z) \bigl\{\log(z\underline{m}) + z\underline{m}(z)\bigr\} \dif z
    = I_2(f). 
\end{align*}
For selected test function $f$ such that $I_2(f)\neq 0$, the first term in \eqref {eq:power_mean_decomp_tau1_eq1} is of order $\bigo(\ell^2/p)$, which dominates the remainder term $\bigo(\ell/p)$ provided that $\ell^2/p \gg 1$. 

In summary, we have 
\[ 
    T_p(n,i) \asymp \rho\biggl(\frac{i}{n}\biggr) \biggl| \ell I_1(f,\tau_1) + \frac{\ell^2}{k_2^{\star}} (\tau_2-1) I_2(f) \biggr|.
\]
Moreover, if one of the following conditions holds:
\begin{itemize}
    \item $\tau_1 \neq 1$ and $I_1(f,\tau_1) \neq 0$;
    \item $\tau_1 = 1$, $\tau_2\neq 1$, and $I_2(f) \neq 0$,
\end{itemize}
then for large enough $\ell$ satisfying $p^{1/2} \ll \ell \ll p$, we have $\sup_{i\geq 1} T_p(n,i) \to \infty$.
This complete the proof of Theorem \ref{thm:power}.

\section{Proof of Lemma \ref{lem:trace_limit_Sigma}}\label{sec:proof_lem_trace_limit_Sigma}
We first assume in addition that the entries of $x_i$ are i.i.d.\ standard Gaussian for all $i \in \intset{k}$. This allows us to exploit the orthogonal invariance of the model and to apply the \emph{Weingarten calculus method} \citep{Collins2006integration} to estimate the tracial quantities in Lemma \ref{lem:trace_limit_Sigma}. In the second step, the Gaussian assumption can be removed by using a standard \emph{Green function comparison} argument. In the following, we only provide the estimation for the second limit in Lemma \ref{lem:trace_limit_Sigma}, namely $\E \Tr\bigl(\tilde{\calG}\Sigma_1 \calG \Sigma_1\bigr)$, since the first one can be shown similarly.

\subsection{Weingarten calculus}
Let $Z_1 \in \mathbb{R}^{p\times k_1}, Z_2 \in \mathbb{R}^{p\times k_2^{\star}}$ have i.i.d.\ standard Gaussian entries. 
WLOG, we can assume that $\Sigma_1 = \diag(d_1, \ldots, d_p)$ is diagonal in the following analysis. 
We define the sample covariance matrices 
\[ 
    S_1^g \deq \frac{1}{k_1}Z_1Z_1^{\top}, 
    \quad 
    S_{2,k^\star}^g \deq \frac{1}{k_2^{\star}}Z_2Z_2^{\top}.
\]
We denote 
\[ 
    G^g(z) \deq \bigl\{ (S_1^{g})^{-1}S_2^g -z \bigr\}^{-1},
    \quad \calG^g(z) \deq G^gS_1^{-1}, 
    \quad \tilde{\calG}^g(z) \deq (G^g)^2 S_1^{-1}.
\]
We aim to estimate
\begin{align*}
    \E \Tr\bigl(\tilde{\calG}^g \Sigma_1 \calG^g \Sigma_1\bigr).
\end{align*}
In the following, we will drop the superscript $g$ for simplicity.
The matrix $\calG(z)$ has the unitary invariant property, namely, for any $p\times p$ Haar distributed orthogonal matrix $O$ independent of $Z_1, Z_2$,
\[ 
    O \calG(z) O^{\top} \overset{\mathtt{d}}{=} \calG(z).
\]
Using this property, we have
\begin{align*}
    &\; \E \Tr \left( \tilde{\calG} \Sigma_1 \calG  \Sigma_1\right) 
    = \E \Tr \left(O \tilde{\calG} O^{\top} \Sigma_1 O\calG O^{\top}  \Sigma_1 \right) \notag \\
    = &\; \sum_{\substack{i_1, i_2, i_3\\ j_1, j_2, j_3}} d_{j_2}d_{i_1} \E  \left(O_{i_1 j_1} \tilde{\calG}_{j_1 i_2} O_{j_2 i_2} O_{j_2 i_3} \calG_{i_3 j_3} O_{i_1 j_3} \right) \notag \\
    =&\; \sum_{\substack{i_1, i_2, i_3\\ j_1, j_2, j_3}} d_{j_2}d_{i_1} \E \bigl(\tilde{\calG}_{j_1 i_2} \calG_{i_3 j_3}\bigr) \E \bigl(O_{i_1 j_1} O_{i_1 j_3}  O_{j_2 i_2} O_{j_2 i_3} \bigr).
\end{align*}
Using the method of \emph{Weingarten calculus on orthogonal group} \citep[see e.g.][Lemma 6.2]{BHPZ17} we have 
\begin{align*}
	&\E  \left(O_{i_1 j_1} O_{i_1 j_3}  O_{j_2 i_2} O_{j_2 i_3} \right) = \frac{3}{p(p+2)} \one_{\{i_1=j_2\}} \one_{\{j_1=j_3=i_2=i_3\}} \notag \\
	&\qquad + \frac{1}{p(p+2)} \one_{\{i_1=j_2\}} \left(\one_{\{j_1=j_3\neq i_2=i_3\}} + \one_{\{j_1=i_2\neq j_3=i_3\}} + \one_{\{j_1=i_3\neq i_2=j_3\}}\right) \notag \\
	&\qquad + \frac{p+1}{p(p-1)(p+2)}\one_{\{i_1\neq j_2\}} \one_{\{j_1=j_3\neq i_2=i_3\}} \notag \\
	&\qquad - \frac{1}{p(p-1)(p+2)}\one_{\{i_1\neq j_2\}} \left(\one_{\{j_1=i_2\neq j_3=i_3\}} + \one_{\{j_1=i_3\neq i_2=j_3\}}\right).
\end{align*}
Hence, 
\begin{align*}
	&\;\E \Tr \bigl( \tilde{\calG} \Sigma_1 \calG \Sigma_1\bigr)   \\
    =&\; \frac{3 \Tr(\Sigma_1^2)}{p(p+2)}  \sum_{a} \E  (\tilde{\calG}_{ aa} \calG_{ aa}) 
    + \frac{\Tr(\Sigma_1^2)}{p(p+2)} \sum_{a\neq b} \E  (\tilde{\calG}_{ aa} \calG_{ bb} + \tilde{\calG}_{ ab} \calG_{ ba} + \tilde{\calG}_{ ab} \calG_{ ab}) \notag \\
    &\qquad + \frac{(p+1)\{(\Tr \Sigma_1)^2 - \Tr(\Sigma_1^2)\}}{p(p-1)(p+2)} \sum_{a\neq b} \E  (\tilde{\calG}_{ ab} \calG_{ ba}) \notag \\
    &\qquad - \frac{\{ (\Tr \Sigma_1)^2 - \Tr(\Sigma_1^2)\} }{p(p-1)(p+2)} \sum_{a\neq b} \E (\tilde{\calG}_{ aa} \calG_{ bb} + \tilde{\calG}_{ ab} \calG_{ ab})\\
	=&\;\frac{3 \Tr(\Sigma_1^2)}{p(p+2)} \sum_{a} \E  (\tilde{\calG}_{ aa} \calG_{ aa}) + \frac{\Tr(\Sigma_1^2)}{p(p+2)} \sum_{a\neq b} \E  (\tilde{\calG}_{ aa} \calG_{ bb}) \notag \\
    &\qquad + \frac{(\Tr \Sigma_1 )^2}{p(p+2)} \sum_{a\neq b} \E  (\tilde{\calG}_{ ab} \calG_{ ba})  - \frac{(\Tr \Sigma_1)^2}{p(p-1)(p+2)} \sum_{a\neq b} \E (\tilde{\calG}_{ aa} \calG_{ bb}) + \bigo(1) \\
    =&\; \frac{\Tr(\Sigma_1^2)}{p^2} \E  \Tr(\tilde{\calG}) \Tr(\calG) + \frac{(\Tr \Sigma_1 )^2}{p^2} \E  \Tr(\tilde{\calG}\calG) - \frac{(\Tr \Sigma_1 )^2}{p^3} \E  \Tr(\tilde{\calG}) \Tr(\calG)+\bigo(1) \notag \\
    =&\; \frac{\tau_2 - \tau_1^2}{p} \E  \Tr(\tilde{\calG}) \Tr(\calG) + \tau_1^2 \E  \Tr(\tilde{\calG}\calG)+ \bigo(1).
\end{align*}

\subsection{Green function comparison}\label{sec:green_comparison}
The goal of this section is to show the following equation using the Green function comparison method:
\begin{align}\label{eq:GFT}
    \E \Tr\bigl(\tilde{\calG} \Sigma_1 \calG \Sigma_1\bigr) 
    - \E \Tr\bigl( \tilde{\calG}^g \Sigma_1 \calG^g \Sigma_1\bigr) 
    = \bigo(1).
\end{align}
This equation allows us to extend the Weingarten calculus result to the general distribution setting.

We assume that $Z_1$ and $Z_2$ are independent of each other and independent of $X_1$ and $X_2$.
We construct a continuous interpolation of between the original data matrices and their Gaussian counterparts. 
For $ i = 1, 2$, define  
\begin{align*}
	X_i^t \deq \sqrt{t}Z_i+\sqrt{1-t}X_i, \qquad t \in [0,1].
\end{align*}
Let $S_i^t$ the corresponding sample covariance matrices formed by interpolated matrix $X_i^t$, and $\calG^t$ and $\tilde{\calG}^t$ be the resolvent matrices defined similarly as $\calG$ and $\tilde{\calG}$ but with $S_i$ replaced by $S_i^t$.

The LHS of \eqref{eq:GFT} can be rewritten as 
\begin{align*}
	\int_0^1 \frac{\dif }{\dif t} \E \Tr\bigl(\tilde{\calG}^t \Sigma_1 \calG^t \Sigma_1 \bigr) \dif t 
	=- \int_0^1 \E \Tr\Bigl\{\tilde{\calG}^t \bigl(\dot{S}_2^t - z_1\dot{S}_1^t\bigr) \tilde{\calG}^t\Sigma_1 \calG^t \Sigma_1 
    + \tilde{\calG}^t \Sigma_1 \calG^t \bigl(\dot{S}_2^t - z_2\dot{S}_1^t\bigr) \calG^t\Sigma_1 \Bigr\} \dif t, 
\end{align*}
where
\begin{align*}
    \dot{S}_1^t & \deq \frac{\dif S_1^t}{\dif t} 
    = \frac{1}{k_1} \Bigl\{Z_1Z_1^{\top} - X_1X_1^{\top} + c(t)\bigl(X_1Z_1^{\top} + Z_1X_1^{\top}\bigr) \Bigr\},
    \quad c(t) \deq \frac{1-2t}{2\sqrt{t(1-t)}},
\end{align*}
and $\dot{S_2^t}$ is defined similarly.
It suffices to estimate the order of each term in the integrand. All terms can be handled in a similar way, we only focus on the term 
\[ 
    \E \Tr \bigl(\dot{S_1^t} H_t\bigr),
    \quad \text{where} \quad
    H_t \deq \tilde{\calG}^t \Sigma_1 \calG^t \Sigma_1 \tilde{\calG}^t.
\] 
We expand the trace term as
\begin{align}\label{eq:expectation_trace_expansion}
    \E \Tr \bigl(\dot{S}_1^t H_t\bigr)
    & = \frac{1}{k_1} \E \sum_{\ell=1}^{k_1} \Bigl\{ 
        z_{\ell}^{\top} H_t z_{\ell} - x_{\ell}^{\top} H_t x_{\ell}  +  2c(t) x_{\ell}^{\top} H_t z_{\ell} \Bigr\},
\end{align}
where $x_{\ell} = (x_{\ell,1}, \ldots, x_{\ell,p})^{\top}$ and $z_{\ell}=(z_{\ell,1}, \ldots, z_{\ell,p})^{\top}$ are the $\ell$-th columns of $X_1$ and $Z_1$, respectively.
We denote $\partial_x M \deq \partial M /\partial x$.
Using the chain rule, it is straightforward to check that 
\begin{align*}
    & \partial_{z_{\ell,a}} H_{t,ij} = \sqrt{t} \partial_{x_{\ell,a}^t} H_{t,ij},
    \qquad 
    \partial_{x_{\ell,a}} H_{t,ij} = \sqrt{1-t} \partial_{x_{\ell,a}^t} H_{t,ij}.
\end{align*}
Using Lemma \ref{lem:cumulant_expansion}, we have
\begin{align*}
    \E x_{\ell}^{\top} H_t x_{\ell} & = \E \sum_{i,j} H_{t, ij} x_{\ell,i} x_{\ell,j}
    = \sum_{i,j} \E \partial_{x_{\ell,i}} \bigl(H_{t, ij} x_{\ell,j}\bigr) + \mathsf{Rem}_1 \\
    & = \sum_{i\neq j} \E \bigl(\partial_{x_{\ell,j}} \partial_{x_{\ell,i}} H_{t, ij}\bigr) + \E \Tr(H_t)  + \mathsf{Rem}_1 + \mathsf{Rem}_2\\ 
    & = (1-t)\sum_{i\neq j} \E \bigl(\partial_{x_{\ell,j}^t} \partial_{x_{\ell,i}^t} H_{t, ij}\bigr) + \E \Tr(H_t)  + \mathsf{Rem}_1 + \mathsf{Rem}_2, 
\end{align*}
where the remainder terms are bounded as
\[
    |\mathsf{Rem}_1| \lesssim \sum_{i\neq j} \E |\partial_{x_{\ell,i}}^2 H_{t}| = \bigo(1),
    \quad 
    |\mathsf{Rem}_2| \lesssim \sum_{i\neq j} \E |\partial_{x_{\ell,i}}^2 \partial_{x_{\ell,j}} H_{t}| = \bigo(1/p).
\]
Similarly, we obtain that
\begin{align*}
    \E z_{\ell}^{\top} H_t z_{\ell} 
    & = t \sum_{i\neq j} \E \bigl(\partial_{x_{\ell,j}^t} \partial_{x_{\ell,i}^t} H_{t, ij}\bigr) + \E \Tr(H_t),\\
    \E x_{\ell}^{\top} H_t z_{\ell} 
    & = \sqrt{t(1-t)} \sum_{i\neq j} \E \bigl(\partial_{x_{\ell,j}^t} \partial_{x_{\ell,i}^t} H_{t, ij}\bigr) + \bigo(1).
\end{align*}
These estimates and \eqref{eq:expectation_trace_expansion} yield
\[ 
    \E \Tr \bigl(\dot{S}_1^t H_t\bigr) = \bigo(1).
\]
Similarly, we can show that the other term in the integrand is also of order $\bigo(1)$. This completes the proof.

\bibliographystyle{abbrvnat}
\bibliography{reference}

\begin{appendix}
\section{Auxiliary lemmas}\label{sec:aux_lemmas}

\begin{lemma}[H\'{a}jek-R\'{e}nyi inequality]\label{lem:hajek-renyi}
Suppose $\{X_n, n\geq 1\}$ is a martingale difference sequence with $\sigma_n^2 \deq \E X_n^2 < \infty$. 
Let $S_n \deq \sum_{j=1}^n X_j$ and let $c_1\geq c_2\geq \cdots \geq c_n >0$ be given constants. 
Then, for any $m\in\intset{n}$ and any $x>0$, 
\[
    \Prob\Bigl( \max_{m \leq j\leq n} c_j |S_j| \geq x \Bigr) \leq \frac{1}{x^2}\biggl( c_m^2 \sum_{j=1}^m \sigma_j^2 + \sum_{j=n+1}^n c_j^2 \sigma_j^2 \biggr).
\]
\end{lemma}

\begin{lemma}[\citet{dharmadhikari1968bounds}]\label{lem:DFJ_inequality}
    Let $\{X_i, i\ge 1\}$ be a martingale difference sequence, and $S_n \deq \sum_{i=1}^n X_i$. Then, for any $r\geq 2$, 
    \[ 
        \E |S_n|^r \leq C_r n^{r/2 - 1} \sum_{i=1}^n \E |X_i|^r, 
    \]
    where $C_r = \{8(r-1)\max(1,2^{r-3})\}^r$.
\end{lemma}

\begin{lemma}[Corollary A.4 in \citet{dette2019determinants}]\label{lem:increment_control}
    Let $\{Z_t, 0\le t\le 1\}$ be a stochastic process with right-continuous sample paths. Assume that there exist constants $\gamma >0$, $\delta > 1$ such that for $r \le s \le t$ the inequality 
    \[ 
        \Prob\bigl( \min\{ |Z_t -Z_s|, |Z_r -Z_s| \} > \lambda \bigr) \le C \lambda^{-\gamma} |t-r|^{\delta}
    \]
    holds with a universal constant $C$. Further assume that there exist a function $\eta$ such that the inequality 
    \[ 
        \Prob\bigl( |Z_t - Z_s| > \alpha \bigr) \le \eta(|t-s|, \alpha)
    \] 
    is satisfied. Then, for any $\varepsilon >0$ and $r \in [0,1]$, we have
    \[ 
        \Prob\Bigl( \sup_{\substack{|s|, |t| <\varepsilon, r+s, r+t \in [0,1]}} |Z_{r+s} -Z_{r+t}| > 4 \lambda \Bigr) \le \varepsilon^{\delta} C C' \lambda^{-\gamma} + 4 \eta(\varepsilon, \lambda),
    \]
    where $C'$ is a constant depending on $\gamma$ and $\delta$ only.
\end{lemma}

\begin{lemma}[Cumulant expansion formula, Lemma 3.1 in \cite{10.1214/16-AAP1237}]\label{lem:cumulant_expansion}
    Let $h$ be a real-valued random variable with finite moments, and $f$ is a complex-valued smooth function on $\mathbb{R}$ with bounded derivatives. Let $c_k(h)$ be the $k$-th cumulant of $h$, given by $c_k(h)\deq(-\iu)^k\frac{\dif}{\dif t}\log \E e^{ith}|_{t=0}$. Then for any fixed $\ell \in \mathbb{N}$, we have
    \begin{equation*}
        \E [hf(h)]=\sum_{k=0}^{\ell} \frac{1}{k!}c_{k+1}(h) \E [f^{(k)}(h)]+R_{\ell+1}.
    \end{equation*}
    For any fixed cutoff $M>0$, the remainder term $R_{\ell+1}$ satisfies
    \begin{equation*}
        |R_{\ell+1}| = \bigo(1) \cdot \E |h|^{\ell+2}\sup_{|x|\leq M} \cdot\bigl|f^{(\ell+1)}(x)\bigr|+ \bigo(1) \cdot \E \bigl( |h|^{\ell+2}\one_{\{|h|>M\}} \bigr) \cdot \|f^{(\ell+1)}\|_{\infty}.
    \end{equation*}
\end{lemma}

\section{Proofs}  \label{sec:lemma_proof}

\subsection{Calculations of mean and variance}\label{sec:proof_explicit_mean_var}

In this section, we provide detailed calculations for the explicit forms of the mean and variance given in Section \ref{sec:simulation}.
Recall that the mean and variance involve contour integrals of the form
\begin{align*}
    \mu_k(f) 
    & = \frac{1}{2 \pi \iu} \oint_{\Gamma} \biggl(f'(z) z \underline{m} - f(z) \frac{\underline{m}'}{\underline{m}} \biggr) \dif z  + \frac{\nu_4-3}{4p \pi \iu} \oint_{\Gamma} f'(z) (1+z\underline{m})^2 \dif z\\
    & \qquad + \frac{1}{2 k_2 \pi \iu} \oint_{\Gamma} \biggl(-\frac{(zf)'' z \underline{m}}{2} + z f'(z) \frac{\underline{m}'}{\underline{m}} + \frac{f''(z)}{\underline{m}}\biggr) \dif z + \smallo\biggl(\frac{1}{p}\biggr),\\
    \sigma_k^2(f) & = \frac{\nu_4-3}{k_2 c_2} \biggl(\frac{1}{2\pi \iu}\oint_{\Gamma} f ' z \underline{m}\dif z\biggr)^2
    + \frac{1}{k_2 \pi \iu} \oint_{\Gamma} \frac{\{f'(z)\}^2}{\underline{m}} \dif z + \smallo\biggl(\frac{1}{p}\biggr).
\end{align*}   
We evaluate the contour integrals by applying the residue theorem to the region exterior to the contour $\Gamma$, which corresponds to the negative sum of the residues outside $\Gamma$. For polynomial functions $g(z) = z$ or $z^2$, the only non-vanishing contribution arises from the residue at infinity. In the case where $g(z) = \log(1+z)$, the residue at the pole $z=-1$ is also included. It is verified that the potential singularities at $z=0$ and $z=-c_2/c_1$ are removable, and thus their residues vanish. 
To illustrate the calculation procedure, we consider the integral 
\[  
    \oint_{\Gamma} \frac{\{f'(z)\}^2}{\underline{m}(z)} dz
\] as a representative example. 
We first establish the Laurent expansion of the reciprocal Stieltjes transform at $z=\infty$ up to the third order:
\begin{align*}
    \frac{1}{\underline{m}(z)} & = -z + M_1 - \frac{M_1^2 - M_2}{z} - \frac{M_1^3 - 2M_1 M_2 + M_3}{z^2} - \frac{C_3}{z^3} + \smallo\biggl(\frac{1}{z^3}\biggr),
\end{align*}
where $C_3 = M_1^4 - 3M_1^2 M_2 + 2M_1 M_3 + M_2^2 - M_4$, and the coefficient $M_i$ is the $i$-th moment of the probability distribution with Stieltjes transform $\underline{m}(z)$.
The explicit forms of $M_i$'s are given by
\begin{align*}
    & M_1 = \frac{c_2}{1-c_1},
    \quad 
    M_2 = \frac{c_2\bigl(1 + c_2 - c_1 c_2\bigr)}{(1-c_1)^3},
    \quad 
    M_3 = \frac{c_2\bigl(c_1^{2}c_2^{2} - 2c_1c_2^{2} - 3c_1c_2 + c_1 + c_2^{2} + 3c_2 + 1\bigr)}{(1-c_1)^5},\\
    & M_4 = \frac{c_2}{(1-c_1)^7}\Big( -c_1^{3}c_2^{3} + 3c_1^{2}c_2^{3} + 6c_1^{2}c_2^{2} - 4c_1^{2}c_2 + c_1^{2} - 3c_1c_2^{3} - 12c_1c_2^{2} - 2c_1c_2 + 3c_1  + c_2^{3} + 6c_2^{2} + 6c_2 + 1 \Big).
\end{align*}
For the case $g(z) = z$, the integrand simplifies to $1/\underline{m}(z)$, yielding a residue at infinity of $M_1^2 - M_2$. Hence, 
\[ 
    \frac{1}{2 \pi \iu}\oint_{\Gamma} \frac{1}{\underline{m}(z)} \dif z = -\mathrm{Res} \biggl(\frac{1}{\underline{m}(z)}, \infty\biggr) = - M_1^2 + M_2.
\]
For $g(z) = z^2$, the integrand becomes $4z^2/\underline{m}(z)$, and the residue at infinity is $4C_3$, leading to
\[ 
    \frac{1}{2 \pi \iu}\oint_{\Gamma} \frac{4z^2}{\underline{m}(z)} \dif z = -\mathrm{Res} \biggl(\frac{4z^2}{\underline{m}(z)}, \infty\biggr) = -4C_3.
\]  
In the case of $g(z) = \log(1+z)$, the integrand is $\frac{1}{(1+z)^2 \underline{m}(z)}$. This function possesses a pole at $z=-1$ with residue $-\frac{m'(-1)}{m(-1)^2}$. At infinity, the asymptotic behavior is derived as:
\[  
    \frac{1}{(1+z)^2 \underline{m}(z)} = \biggl(\frac{1}{z^2} - \frac{2}{z^3}\biggr)\bigl\{-z + M_1 + \smallo(1)\bigr\} = - \frac{1}{z} + \smallo\biggl(\frac{1}{z}\biggr),
\]
which corresponds to a residue of $1$. Therefore, we have
\[  
    \frac{1}{2 \pi \iu}\oint_{\Gamma} \frac{1}{(1+z)^2 \underline{m}(z)} \dif z = -\mathrm{Res} \biggl(\frac{1}{(1+z)^2 \underline{m}(z)}, -1\biggr) - \mathrm{Res} \biggl(\frac{1}{(1+z)^2 \underline{m}(z)}, \infty\biggr) = \frac{\underline{m}'(-1)}{\underline{m}(-1)^2} - 1.
\]
The other integrals in the expressions of $\mu_k(f)$ and $\sigma_k^2(f)$ can be evaluated using similar procedures and the results are summarized in Table \ref{tab:integral_summary}.
\begin{table}[htbp]
\centering
\renewcommand{\arraystretch}{1.4}
\setlength{\tabcolsep}{8pt}
\begin{tabular}{l|c|c|c}
\toprule
Integrals & $f(z) = z$ & $f(z) = z^2$ & $f(z) = \log(1+z)$ \\
\midrule
$\frac{1}{2\pi \iu} \oint f' z \underline{m} \dif z$ & $-M_1$ & $-2 M_2$ & $\underline{m}(-1) - 1$ \\
\hline
$\frac{1}{2\pi \iu} \oint f \frac{\underline{m}'}{\underline{m}} \dif z$ & $-M_1$ & $M_1^2 - 2 M_2$ & $\ln \underline{m}(-1)$ \\
\hline
$\frac{1}{2\pi \iu} \oint f' (1+z \underline{m})^2 \dif z$ & $0$ & $2 M_1^2$ & $-\{1 - \underline{m}(-1)\}^2$ \\
\hline
$\frac{1}{2\pi \iu} \oint (z f)'' z \underline{m}' \dif z$ & $ - 2 M_1$ & $ - 6 M_2$ & $ - 1 + \underline{m}'(-1)$ \\
\hline
$\frac{1}{2\pi \iu} \oint f' z \frac{\underline{m}'}{\underline{m}}\dif z$ & $-M_1$ & $2 M_1^2 - 4 M_2$ & $\frac{\underline{m}'(-1)}{\underline{m}(-1)} - 1$ \\
\hline
$\frac{1}{2\pi \iu} \oint \frac{f''}{\underline{m}}\dif z$ & $0$ & $-2 (M_1^2 - M_2)$ & $1 - \frac{\underline{m}'(-1)}{\underline{m}(-1)^2}$ \\
\hline
$\frac{1}{2\pi \iu} \oint \frac{(f')^2}{\underline{m}}\dif z$ & $- M_1^2 + M_2$ & $-4 C_3$ & $\frac{\underline{m}'(-1)}{\underline{m}(-1)^2} - 1$ \\
\bottomrule
\end{tabular}
\caption{Summary of contour integrals appearing in the evaluation of mean and variance.}
\label{tab:integral_summary}
\end{table}

\subsection{Proof of Proposition \ref{pro:stepone_diff_mean_var}}\label{sec:proof_stepone_diff_mean_var}
Recall the definition of $M_k(z)$ in \eqref{eq:Lk_Cauchy}. The following two lemmas provide the mean and covariance structure of the resolvent difference $M_k(z)$ and some useful trace limits, which are essential for proving Proposition \ref{pro:stepone_diff_mean_var}.
\begin{lemma}\label{lem:resolvent_mean_var}
Under Assumptions \ref{asm:data} and \ref{asm:sample_dimension}, for any $z_1, z_2 \in \mathbb{C}^+$, we have
\begin{align*}
    & \E M_k(z) = \biggl( 1 + \frac{1}{z\underline{m}}\biggr) (z\underline{m})' 
    + \frac{1}{k_2} \biggl\{\frac{z(z\underline{m})''}{2}
    +\biggl(\frac{z\underline{m}'}{\underline{m}}\biggr)'
    - \biggl(\frac{1}{\underline{m}}\biggr)'' \biggr\} 
    +\frac{\nu_4-3}{2p}\{(1+z\underline{m})^2\}'
    + \smallo(1/p),\\
    & \cov (M_k(z_1), M_k(z_2)) = \frac{\nu_4-3}{k_2c_2} \partial_{z_1} \bigl(z_1\underline{m}_1 \bigr) \partial_{z_2} \bigl(z_2\underline{m}_2 \bigr) - \frac{2}{k_2} \partial_{z_1}\partial_{z_2} \biggl( \frac{1/\underline{m}_1-1/\underline{m}_2}{z_1-z_2} \biggr) + \smallo(1/p).
\end{align*}
\end{lemma}

\begin{lemma}\label{lem:trace_limits}
    Under Assumptions \ref{asm:data} and \ref{asm:sample_dimension}, the following limits hold as $p \to \infty$:
    \begin{align}
        &\frac{\Tr G_{k-1}}{k_2}  \pto c_2 m(z) = \underline{m}(z) + \frac{1-c_2}{z}, 
        \quad 
        \frac{z \Tr G_{k-1}^2}{k_2} \pto z c_2 m'(z) = z\underline{m}'(z) - \frac{1-c_2}{z}, \label{eq:trace_G}\\
        &\frac{\Tr (G_{k-1}S_1^{-1})}{k_2}  \pto -  \frac{1+z \underline{m}(z)}{z \underline{m}(z)}, 
        \quad
        \frac{\Tr (G_{k-1}^2S_1^{-1})}{k_2} 
        \pto \frac{\underline{m}(z)+z\underline{m}'(z)}{\{z\underline{m}(z)\}^2}, \label{eq:trace_GS} \\  
        &\frac{1}{k_2} \Tr \bigl\{G_{k-1}(z_1)S_1^{-1}G_{k-1}(z_2)S_1^{-1}\bigr\} 
	    \pto -\frac{1}{z_1z_2\underline{m}_1\underline{m}_2} \biggl( 1+ \frac{1/\underline{m}_1 - 1/\underline{m}_2}{z_1-z_2} \biggr).\label{eq:trace_G1S1G2S1}     
    \end{align}
\end{lemma}
The proof of Lemma \ref{lem:resolvent_mean_var} and Lemma \ref{lem:trace_limits}  is given in Section \ref{sec:proof_lem_resolvent_mean_var} and Section \ref{sec:proof_lem_trace_limits} respectively. We first use the above two lemmas to prove Proposition \ref{pro:stepone_diff_mean_var}. From the definition of $L_k(f)$, Lemma \ref{lem:resolvent_mean_var} and the integration by part, we have
\begin{align*}
    \mu_k(f) & =-\frac{1}{2 \pi \iu} \oint_{\Gamma} f(z) \biggl\{(z\underline{m})' + \frac{\underline{m}'}{\underline{m}}\biggr\} \dif z - \frac{\nu_4-3}{4p \pi \iu} \oint_{\Gamma} f(z) \bigl\{(1+z\underline{m})^2\bigr\}' \dif z\\ 
    &\qquad - \frac{1}{2 k_2 \pi \iu} \oint_{\Gamma} f(z) \biggl\{ \frac{z(z\underline{m})''}{2} + \biggl(\frac{z\underline{m}'}{\underline{m}}\biggr)' - \biggl(\frac{1}{\underline{m}}\biggr)'' \biggr\} \dif z\\
    & = \frac{1}{2 \pi \iu} \oint_{\Gamma} \biggl(f'(z) z \underline{m} - f(z) \frac{\underline{m}'}{\underline{m}} \biggr) \dif z  + \frac{\nu_4-3}{4p \pi \iu} \oint_{\Gamma} f'(z) (1+z\underline{m})^2 \dif z\\
    & \qquad + \frac{1}{2 k_2 \pi \iu} \oint_{\Gamma} \biggl(-\frac{(zf)'' z \underline{m}}{2} + z f'(z) \frac{\underline{m}'}{\underline{m}} + \frac{f''(z)}{\underline{m}}\biggr) \dif z.
\end{align*}
We take $\Gamma_1 = \Gamma$ and choose $\Gamma_2$ to be a rectangle outside $\Gamma_1$, with each of its four sides at a distance $\varepsilon$ from $\Gamma_1$.
Using the integration by parts twice, first with respect to $z_1$ and then with respect to $z_2$, we have
\begin{align*}
    &\;\oint_{\Gamma_1\times \Gamma_2} f(z_1)f(z_2) \partial_{z_1}\partial_{z_2} \biggl(\frac{1/\underline{m}_1-1/\underline{m}_2}{z_1-z_2}\biggr) \dif z_1\dif z_2 
    =\;  \oint_{\Gamma_1} f'(z_1) \oint_{\Gamma_2}f'(z_2)\biggl(\frac{1/\underline{m}_1-1/\underline{m}_2}{z_1-z_2}\biggr) \dif z_2\dif z_1 \\
    =&\; -\oint_{\Gamma_2}\frac{f'(z_2)}{\underline{m}_2}\oint_{\Gamma_1} \frac{f'(z_1)}{z_1-z_2}\dif z_1 \dif z_2 
    =\; -2\pi\iu\oint_{\Gamma_2} \frac{\{f'(z)\}^2}{\underline{m}(z)} \dif z, 
\end{align*}
and thus
\begin{align*}
    \sigma_k^2(f) & = \biggl(\frac{-1}{2 \pi \iu}\biggr)^2 \oint_{\Gamma_1}\oint_{\Gamma_2} f(z_1) f(z_2) \cov(M_k(z_1), M_k(z_2)) \dif z_1 \dif z_2 \\
    & = \frac{\nu_4-3}{k_2c_2} \biggl(-\frac{1}{2 \pi \iu}\oint_{\Gamma_1} f(z) (z\underline{m})' \dif z\biggr)^2 
    + \frac{1}{2 k_2 \pi^2} \oint_{\Gamma_1\times \Gamma_2} f(z_1)f(z_2) \partial_{z_1}\partial_{z_2} \biggl(\frac{1/\underline{m}_1-1/\underline{m}_2}{z_1-z_2}\biggr) \dif z_1\dif z_2   \\
    & = \frac{\nu_4-3}{k_2 c_2} \biggl(\frac{1}{2\pi \iu}\oint_{\Gamma} f ' z \underline{m}\dif z\biggr)^2
    + \frac{1}{k_2 \pi \iu} \oint_{\Gamma} \frac{\{f'(z)\}^2}{\underline{m}} \dif z.
\end{align*}

\subsection{Proof of Remark \ref{rmk:integral}} \label{sec:proof_rmk_integral}

We denote $\underline{m}(x\pm \iu 0) = A \pm \iu B$, where \[
    A \equiv A(x) = -\frac{x(h^2+c_1)+c_2(1-c_2)}{2x(c_2+xc_1)},
    \quad 
    B \equiv B(x) = \frac{c_2(1-c_1)\sqrt{(x-b)(a-x)}}{2x(c_2+xc_1)}.
\]
We choose the counter $\Gamma$ to be a rectangle with sides parallel to the axes. The rectangle enclose the intervel $[a,b]$, and the horizontal sides are a distance $\varepsilon < 1$ away from the real axis. We let $\varepsilon \to 0$.
Using integration by parts, for any analytic function $f(z)$, we obtain
\begin{align}
   & \oint_{\Gamma} f(z) \{z\underline{m}(z)\}' \dif z  
    = -\oint_{\Gamma} f'(z)z\underline{m}(z) \dif z \notag\\
    & = - \int_b^a f'(x)x(A+\iu B) \dif x - \int_a^b f'(x)x(A-\iu B) \dif x =  2\iu \int_a^b f'(x) x B  \dif x,\label{eq:int_g_z_m}
\end{align}
and
\begin{align}
    &\oint_{\Gamma} \frac{f(z)\underline{m}'(z)}{\underline{m}(z)} \dif z 
    = - \oint_{\Gamma} f'(z) \log \{\underline{m}(z)\} \dif z \notag \\
    & = -\int_b^a f'(x) \log (A + \iu B) \dif x - \int_a^b f'(x) \log (A - \iu B) \dif x  
    = 2\iu \int_a^b f(x) \frac{A'B-AB'}{A^2+B^2} \dif x.\label{eq:int_g_mp_m}
\end{align}
Similarly, we can also get 
\begin{align}
    \oint_{\Gamma}\frac{f(z)}{\underline{m}(z)} \dif z 
    = 2\iu \int_{a}^b\frac{B f}{A^2+B^2} \dif x,\label{eq:int_g_over_m}
\end{align}
\begin{align*}
    \oint_{\Gamma} f(z) z (z\underline{m})'' \dif z = -\oint_{\Gamma} \{zf(z)\}'(z\underline{m})'\dif z \overset{\eqref{eq:int_g_z_m}}{=} -2\iu \int_a^b (xf''+2f') x B \dif x,
\end{align*}
\begin{align*}
   \oint_{\Gamma} \frac{f'(z)z\underline{m}'}{\underline{m}} \dif z \overset{\eqref{eq:int_g_mp_m}}{=} 2 \iu \int_a^b x f' \frac{A'B-AB'}{A^2+B^2} \dif x,
\end{align*}
\begin{align*}
    \oint_{\Gamma} f(z)\biggl(\frac{1}{\underline{m}}\biggr)''\dif z 
    = \oint_{\Gamma} \frac{f''(z)}{\underline{m}} \dif z 
    \overset{\eqref{eq:int_g_over_m}}{=} 2\iu  \int_a^b \frac{B f''}{A^2+B^2} \dif x,
\end{align*}
\begin{align*}
    \oint_{\Gamma} f(z) \{(1+z\underline{m})^2\}'\dif z 
    = - \oint_{\Gamma} f'(z) (1+z\underline{m})^2\dif z 
    = 4 \iu \int_a^b (1+xA) x B f' \dif x,
\end{align*}
\begin{align*} 
    \oint_{\Gamma_2} \frac{\{f'(z)\}^2}{\underline{m}(z)} \dif z 
    = 2\iu \int_a^b \frac{B(x)\{f'(x)\}^2}{A^2(x)+B^2(x)} \dif x.
\end{align*}
From these integrals and Lemma \ref{lem:resolvent_mean_var}, we have
\begin{align*}
    \E (L_k(f)) 
    &= -\frac{1}{\pi}\int_a^b\biggl( x B f'
    + \frac{A'B-AB'}{A^2+B^2}f\biggr) \dif x - \frac{\nu_4-3}{p\pi} \int_a^b (1+xA) x B f' \dif x \\
    &\qquad - \frac{1}{k_2 \pi} \int_a^b \biggl\{- \frac{(x f'' + 2f') x B}{2} 
    -\frac{(A'B-AB')x f' + B f''}{A^2+B^2}
    \biggr\} \dif x + \smallo\biggl(\frac{1}{p}\biggr),
\end{align*}
and
\begin{align*}
    \var (L_k(f)) & = \frac{\nu_4-3}{k_2c_2\pi^2}   \biggl(\int_a^b f'(x) x B(x)  \dif x \biggr)^2 
    +\frac{2}{k_2\pi} \int_a^b \frac{B(x)\{f'(x)\}^2}{A^2(x)+B^2(x)}\dif x + \smallo\biggl(\frac{1}{p}\biggr).
\end{align*}

\subsection{Proof of Lemma \ref{lem:resolvent_mean_var}}\label{sec:proof_lem_resolvent_mean_var}

Denote $k_2=k-1-k_1$, 
and we can write 
\[
    F_k = \frac{k_2}{k_2+1} F_{k-1} + \frac{S_1^{-1}x_kx_k^{\top}}{k_2+1}.
\]
This, together with the Sherman-Morrison identity, yields that
\[
    G_k = \hat{G}_{k-1} - \frac{\frac{1}{k_2+1}\hat{G}_{k-1}S_1^{-1}x_kx_k^{\top}\hat{G}_{k-1}}{1+\frac{1}{k_2+1}x_k^{\top}\hat{G}_{k-1}S_1^{-1}x_k},
\]
where $\hat{G} \equiv \hat{G}_{k-1}(z)\deq\bigl(\frac{k_2}{k_2+1} F_{k-1} -z \bigr)^{-1}$.
Using the identity $A^{-1}-B^{-1}=A^{-1}(B-A)B^{-1}$, we have 
\begin{align*}
    \hat{G}_{k-1} 
    = \frac{k_2+1}{k_2}G_{k-1} + \frac{z}{k_2}\hat{G}_{k-1}G_{k-1}.
\end{align*}
For brevity, we denote $G=G_{k-1}$. Plugging the above equation into the expression of $G_k$, we have 
\begin{align*}
    M_k(z) 
    & = \frac{1}{k_2} \Tr G + \frac{z}{k_2} \Tr \hat{G}G 
    - \frac{\frac{1}{k_2+1}x_k^{\top}\hat{G}^2 S_1^{-1}x_k}{1 + \frac{1}{k_2+1}x_k^{\top}\hat{G}S_1^{-1}x_k}\\ 
    & =\frac{1}{k_2}\Tr G + \frac{z(k_2+1)}{k_2^2}\Tr G^2 + \frac{z^2}{k_2^2} \Tr G^3 
    - \frac{\frac{k_2+1}{k_2^2}x_k^{\top}G^2 S_1^{-1}x_k + \frac{2z}{k_2^2}x_k^{\top}G^3S_1^{-1}x_k}{1 + \frac{1}{k_2}x_k^{\top}GS_1^{-1}x_k + \frac{z}{k_2^2} x_k^{\top}G^2 S_1^{-1}x_k} 
    + \smallop(1/p).
\end{align*}
Let $\Delta_k^{(\ell)}\deq\frac{x_k^{\top}G^{\ell}S_1^{-1}x_k - \Tr (G^{\ell}S_1^{-1})}{k_2}$, $\ell = 1, 2$. Using the Taylor expansion, we have
\begin{align*}
    &\;\frac{\frac{k_2+1}{k_2^2}x_k^{\top}G^2 S_1^{-1}x_k + \frac{2z}{k_2^2}x_k^{\top}G^3S_1^{-1}x_k}{1 + \frac{1}{k_2}x_k^{\top}GS_1^{-1}x_k + \frac{z}{k_2^2} x_k^{\top}G^2 S_1^{-1}x_k} \\
    =&\; \frac{\frac{1}{k_2}\Tr G^2 S_1^{-1} + \Delta_k^{(2)} +  \frac{1}{k_2^2} \Tr G^2 S_1^{-1} + \frac{2z}{k_2^2}\Tr G^3S_1^{-1}}{1 + \frac{1}{k_2} \Tr GS_1^{-1} + \Delta_k^{(1)} + \frac{z}{k_2^2} \Tr G^2 S_1^{-1}} + \bigop(p^{-3/2}) \\ 
    =&\; \alpha(z) \biggl(\frac{1}{k_2}\Tr G^2 S_1^{-1} + \frac{1}{k_2^2} \Tr G^2 S_1^{-1} + \frac{2z}{k_2^2}\Tr G^3S_1^{-1}\biggr) 
    - \alpha^2(z) \biggl(\frac{z(\Tr G^2 S_1^{-1})^2}{k_2^3} + \Delta_k^{(1)}\Delta_k^{(2)}\biggr)\\ 
    & \qquad + \alpha^3(z) \frac{\Tr G^2 S_1^{-1}}{k_2} (\Delta_k^{(1)})^2 + \frac{1}{k_2}(x_k^{\top} A_k x_k - \Tr A_k) + \bigop(p^{-3/2}),
\end{align*}
where 
\begin{align}\label{eq:Ak}
    A_k \equiv A_k(z) \deq \alpha(z) G^2 S_1^{-1} + \alpha'(z)  G S_1^{-1},
    \quad 
    \alpha(z) \deq \frac{ 1 }{1+\frac{1}{k_2} \Tr (GS_1^{-1}) }. 
\end{align}
From Lemma \ref{lem:trace_limits}, we have 
\begin{align}\label{eq:alpha}
    \alpha(z) &\pto -z\underline{m}(z),
\end{align}
\begin{align*}
    \E \bigl(\Delta_k^{(1)}\Delta_k^{(2)}\bigr)
    & = -\frac{(\nu_4-3)(1+z\underline{m})(\underline{m}+z\underline{m}')}{p(z\underline{m})^3} \\
    & \qquad + \frac{2c_2}{p(z\underline{m})^2} \biggl(\frac{1}{z} - \frac{\underline{m}'}{z \underline{m}^2} + \frac{\underline{m}'}{\underline{m}} - \frac{2(\underline{m}')^2}{\underline{m}^3} + \frac{\underline{m}''}{2\underline{m}^2}\biggr) + \smallo(1/p), \notag\\ 
    \E \bigl(\Delta_k^{(1)}\bigr)^2 & = \frac{\nu_4-3}{p} \biggl(\frac{1+z\underline{m}}{z\underline{m}}\biggr)^2 -\frac{2 c_2}{p(z\underline{m})^2} \biggl(1-\frac{\underline{m}'}{\underline{m}^2}\biggr) + \smallo(1/p),
\end{align*}
and thus
\begin{equation}\label{eq:Mk}
\begin{aligned}
    M_k(z) & =  \biggl( 1 + \frac{1}{z\underline{m}}\biggr) (z\underline{m})' 
 + \frac{1}{k_2} \biggl\{\frac{z(z\underline{m})''}{2}
    +\biggl(\frac{z\underline{m}'}{\underline{m}}\biggr)'
    - \biggl(\frac{1}{\underline{m}}\biggr)'' \biggr\}  \\ 
    & \qquad +\frac{\nu_4-3}{2p}\{(1+z\underline{m})^2\}' - \frac{1}{k_2}(x_k^{\top} A_k x_k - \Tr A_k) + \smallop(1/p).
\end{aligned}
\end{equation}
This implies the approximation of $\E M_k(z)$ in Lemma \ref{lem:resolvent_mean_var}, and it remains to calculate the covariance function $\cov(M_k(z_1), M_k(z_2))$. 
Using \eqref{eq:Mk}, we have
\begin{align*}
    \cov \bigl(M_k(z_1), M_k(z_2)\bigr)& = \frac{1}{k_2^2}\E \bigl\{\bigl(x_k^{\top} A_k(z_1)x_k - \Tr A_k(z_1)\bigr) \bigl(x_k^{\top} A_k(z_2)x_k - \Tr A_k(z_2)\bigr)\bigr\} + \smallo(1/p)\\
    & = \frac{\nu_4-3}{k_2^2} \E \sum_{i=1}^p [A_k(z_1)]_{ii}[A_k(z_2)]_{ii} + \frac{2}{k_2^2} \E \Tr \bigl\{ A_k(z_1)A_k(z_2) \bigr\} + \smallo(1/p).
\end{align*} 
Using the identity
\[
    A_k(z) = \partial_z \bigl\{\alpha(z) G(z) S_1^{-1}\bigr\},
\]
we have
\begin{align*}
    \E \Tr \bigl\{ A_k(z_1)A_k(z_2) \bigr\}
    & = \E \Tr \bigl[ \partial_{z_1}\{\alpha(z_1) G(z_1) S_1^{-1}\} \partial_{z_2}\{\alpha(z_2) G(z_2) S_1^{-1}\} \bigr] \notag\\
    & = \partial_{z_1}\partial_{z_2} \bigl[\alpha(z_1)\alpha(z_2) \E \Tr \{G(z_1) S_1^{-1}G(z_2) S_1^{-1}\} \bigr].\notag
\end{align*}
This, together with \eqref{eq:alpha} and Lemma \ref{lem:trace_limits}, implies that
\begin{align}\label{eq:trace_A1A2}
    \frac{1}{k_2} \E \Tr \bigl\{ A_k(z_1)A_k(z_2) \bigr\} 
    = -\partial_{z_1}\partial_{z_2} \biggl(\frac{1/\underline{m}_1-1/\underline{m}_2}{z_1-z_2}\biggr) + \smallo(1).
\end{align}
Using Lemma \ref{lem:trace_limits} again, we have
\begin{align*}
    \cov (M_k(z_1), M_k(z_2)) = \frac{\nu_4-3}{k_2c_2} \partial_{z_1} \bigl(z_1\underline{m}_1 \bigr) \partial_{z_2} \bigl(z_2\underline{m}_2 \bigr) - \frac{2}{k_2} \partial_{z_1}\partial_{z_2} \biggl( \frac{1/\underline{m}_1-1/\underline{m}_2}{z_1-z_2} \biggr) + \smallo(1/p).
\end{align*}

\subsection{Proof of Lemma \ref{lem:trace_limits}}\label{sec:proof_lem_trace_limits}

The estimate \eqref{eq:trace_G} follows from Lemma 6.2 in \citet{Zheng2012central}.
It remains to prove \eqref{eq:trace_GS} and \eqref{eq:trace_G1S1G2S1}.

\subsubsection*{Proof of \eqref{eq:trace_GS}} 
The second limit in \eqref{eq:trace_GS} can be derived from the first one by taking derivative with respect to $z$. Thus, it suffices to prove the first limit.
Recall that 
\[ 
    G_{k-1}(z) = \biggl(\frac{1}{k_2}S_1^{-1}XX^{\top} - z\biggr)^{-1},
    \quad 
    X = (x_{k_1+1}, \ldots, x_{k-1})_{p\times k_2}.
\]
We abuse the notation by denoting the $(i,j)$-th entry of $X$ by $x_{ij}$ and the $j$-th column of $X$ by $x_j$. From the definition of $G\equiv G_{k-1}(z)$, we have the following resolvent identity:
\begin{align}\label{eq:G_identity}
    G = \frac{1}{zk_2}S_1^{-1}XX^{\top}G - \frac{1}{z}I,
\end{align}
which implies that
\[
    \E \Tr G = \frac{1}{zk_2} \E \Tr (XX^{\top}GS_1^{-1}) - \frac{p}{z} = \frac{1}{zk_2} \sum_{i,j} \E \bigl\{ x_{ij} f_{ij}(X) \bigr\} - \frac{p}{z},
\]
where 
\[ 
    f_{ij}(X) \deq \{X^{\top} \calG(z)\}_{ji}, 
    \quad 
    \calG(z) \deq G(z) S_1^{-1} = (S_2 - z S_1)^{-1}.
\]
From Lemma \ref{lem:cumulant_expansion}, we have
\begin{align*}
    \E \bigl\{ x_{ij} f_{ij}(X) \bigr\}
    =
    \E \bigl\{\partial_{ij} f_{ij}(X)\bigr\}+ R_{1,ij},
\end{align*}
where the remainder term $R_{1,ij}$ satisfies that   
\[ 
    R_{1,ij} \lesssim \E\|\partial_{ij}^2 f_{ij}(X)\|_{\infty}.
\]
The first derivative is given by
\begin{align*}
    \partial_{ij} f_{ij}(X) 
    & = \{(\partial_{ij} X^{\top}) \calG\}_{ji} - \{X^{\top}  \calG (\partial_{ij} S_2) \calG\}_{ji} = \calG_{ii} - \frac{1}{k_2}(X^{\top} \calG)_{ji}(X^{\top} \calG)_{ji} - \frac{1}{k_2}(X^{\top} \calG X)_{jj}\calG_{ii}.
\end{align*}
Moreover, the second derivative is given by
\begin{align*}
    &\; \partial_{ij}^2 f_{ij}(X)  
    =\; \partial_{ij}  \calG_{ii} - \{\calG (\partial_{ij} S_2) \calG\}_{ii} 
    - \sum_a x_{aj} \Bigl\{ (\partial_{ij} \calG) (\partial_{ij} S_2) \calG + \calG (\partial_{ij}^2 S_2)\calG + \calG (\partial_{ij} S_2) (\partial_{ij} \calG)\Bigr\}_{ai}.
\end{align*}
Let $e_i$ be the $i$-th standard basis vector in $\bbR^p$.
Using the following norm estimates:
\begin{align*}
    &\|\partial_{ij} S_2\|_2 = \frac{1}{k_2}\|e_i x_j^{\top} + x_j e_i^{\top}\|_2 \lesssim \frac{\|x_j\|}{k_2},\quad 
    \|\partial_{ij}^2 S_2\|_2 = \frac{1}{k_2}\|e_i e_i^{\top}\|_2 = \frac{1}{k_2},\\ 
    &\|\partial_{ij} \calG\|_2 = \| - G (\partial_{ij} S_2) G S_1^{-1}\|_2 \le \|S_1^{-1}\|_2 \cdot \|G\|_2^2 \cdot \|\partial_{ij} S_2\|_2 \lesssim \frac{\|x_j\|}{k_2},
\end{align*}
we can easily get 
\begin{align*} 
    \|\partial_{ij}^2 f_{ij}(X)\|_{\infty}  
     \lesssim \frac{\|x_j\|}{k_2} + \frac{\|x_j\|^3}{k_2^2}.  
\end{align*} 
Thus, we obtain $R_{1,ij} = \bigo(1/\sqrt{p})$. 
Combining the above estimates, we have
\begin{align*}
    \E \Tr G
    & = \frac{1}{z} \E \Tr \calG - \frac{1}{zk_2} \E \Tr (\calG S_2 \calG) - \frac{1}{zk_2} \E\{ \Tr (\calG S_2) \Tr \calG\} -\frac{p}{z} + \bigo(\sqrt{p}).
\end{align*}
Using the identity $\calG S_2 = I + z G$ and the following covariance estimate
\begin{align}\label{eq:cov_traceG_traceGS}
    \E \{\Tr G \Tr \calG\} - \E (\Tr G) \E (\Tr \calG) = \bigo(1),
\end{align}
and \eqref{eq:trace_G}, 
we obtain that 
\begin{align*}
    \frac{1}{k_2} \E \Tr \calG & = \biggl(1 - \frac{z \E \Tr G + p}{k_2}\biggr)^{-1} \frac{z\E\Tr G + p}{k_2}  + \bigo(p^{-1/2})
    = - \frac{z \underline{m} + 1}{z \underline{m}} + \bigo(p^{-1/2}).
\end{align*}
This, together with the variance estimate
\begin{align}\label{eq:Var_traceG}
    \var(\Tr \calG) = \bigo(1),
\end{align}
completes the proof of \eqref{eq:trace_GS}.

Now we prove \eqref{eq:cov_traceG_traceGS} and \eqref{eq:Var_traceG}.
We denote
\begin{align*}
    S_2^{(j)} \deq S_2 - \frac{1}{k_2}x_jx_j^{\top}, 
    \quad 
    \calG^{(j)} \deq (S_2^{(j)} - z S_1)^{-1}, 
    \quad 
    j\in\intset{k_2}.
\end{align*}
Let $\calF_j$ be the $\sigma$-algebra generated by $\{x_{1}, \ldots, x_{j}\}$, and denote $\E_j \deq \E(\cdot | \calF_j)$, then we have
\begin{align}\label{eq:Var_traceG_decom}
    \var(\Tr \calG) & = \sum_{j=1}^{k_2} \E \bigl\{ (\E_j - \E_{j-1}) \bigl(\Tr \calG - \Tr \calG^{(j)}\bigr)\bigr\}^2 
    = \sum_{j=1}^{k_2}\E\biggl\{(\E_j - \E_{j-1}) \frac{\frac{1}{k_2}x_j^{\top}(\calG^{(j)})^2x_j}{1+\frac{1}{k_2}x_j^{\top}\calG^{(j)}x_j}\biggr\}^2.
\end{align}
Let $\Delta_j \deq \frac{1}{k_2}\Tr \calG^{(j)} - \frac{1}{k_2}\E \Tr \calG^{(j)}$. From the Taylor expansion, we have
\begin{align*}
     \frac{\frac{1}{k_2}x_j^{\top}(\calG^{(j)})^2x_j}{1+\frac{1}{k_2}x_j^{\top}\calG^{(j)}x_j} = \frac{\frac{1}{k_2} x_j^{\top}(\calG^{(j)})^2x_j}{1+\frac{1}{k_2}\E \Tr\calG^{(j)}} - \frac{\frac{1}{k_2}x_j^{\top}(\calG^{(j)})^2x_j}{(1+\frac{1}{k_2}\E \Tr\calG^{(j)})^2} \Delta_j + \bigo(\Delta_j^2).
\end{align*}
For the first term, we have
\begin{align*}
    \E \biggl| \frac{1}{k_2} (\E_j - \E_{j-1}) x_j^{\top}(\calG^{(j)})^2x_j \biggr|^2 & \lesssim \frac{1}{p^2} \E \Bigl| x_j^{\top} (\calG^{(j)})^2 x_j -\Tr (\calG^{(j)})^2 \Bigr|^2 = \bigo(1/p).
\end{align*}
Similarly, we can show that the expectations of the second and third terms are also $\bigo(1/p)$. These estimates and \eqref{eq:Var_traceG_decom} yield \eqref{eq:Var_traceG}.

Similarly, we can verify that 
\begin{align}\label{eq:cov_traceG_traceGS_proof}
    \var(\Tr G)=\var(\Tr (\calG S_1)) = \bigo(1).
\end{align}    
These estimates, together with the Cauchy-Schwarz inequality, imply \eqref{eq:cov_traceG_traceGS}. 

\subsubsection*{Proof of \eqref{eq:trace_G1S1G2S1}}

In the following, we denote $\calG_i=\calG(z_i)$ and $G_i=G(z_i)$ for $i=1,2$. 
From the identity \eqref{eq:G_identity}, we have
\[ 
    G_2 G_1 S_1^{-1} = \frac{1}{z_1 k_2} \calG_2 X X^{\top} \calG_1 - \frac{1}{z_1} \calG_2,
\]
and thus 
\begin{align*}
    \E \Tr (G_2 G_1 S_1^{-1}) & = \frac{1}{z_1 k_2} \E \Tr (X X^{\top} \calG_1 \calG_2) - \frac{1}{z_1} \E \Tr \calG_2 
    = \frac{1}{z_1 k_2} \sum_{i,j} \E \bigl\{ x_{ij} (X^{\top} \calG_1 \calG_2)_{ji} \bigr\} - \frac{1}{z_1} \E \Tr \calG_2.
\end{align*}
For brevity, we denote 
\[ 
    g_{ij}(X) = (X^{\top} \calG_1 \calG_2)_{ji}.
\]
By using the cumulant expansion, we have
\begin{align*}
    \E \{x_{ij} g_{ij}(X)\} = \E \{\partial_{ij} g_{ij}(X)\} + R_{2,ij},
\end{align*}
where the remainder term $R_{2,ij}$ satisfies that
\[ 
    R_{2,ij} \lesssim \E\|\partial_{ij}^2 g_{ij}(X)\|_{\infty}.
\]
The first derivative of $g_{ij}$ is given by
\begin{align*}
    \partial_{ij} g_{ij}(X) 
    & = (\calG_1 \calG_2)_{ii} 
    - \frac{1}{k_2} (X^{\top} \calG_1)_{ji} (X^{\top} \calG_1 \calG_2)_{ji} 
    - \frac{1}{k_2}(X^{\top} \calG_1 X)_{jj}(\calG_1\calG_2)_{ii} \\
    & \qquad - \frac{1}{k_2}(X^{\top} \calG_1 \calG_2)_{ji}(X^{\top} \calG_2)_{ji} 
    - \frac{1}{k_2}(X^{\top} \calG_1 \calG_2 X)_{jj}(\calG_2)_{ii}.
\end{align*}
The same argument as before shows that the remainder term $R_{2,ij} = \bigo(1/\sqrt{p}).$
Combining the above estimates, we have
\begin{align*}
    &\; \E \Tr (G_2 G_1 S_1^{-1}) \\ 
    =&\; \frac{1}{z_1 k_2} \E\bigl\{(k_2-p-z_1\Tr G_1)\Tr(\calG_1\calG_2) \bigr\}
    - \frac{1}{z_1k_2}\E\biggl\{\biggl(\frac{z_2\Tr \calG_2}{z_2 - z_1} - \frac{z_1\Tr \calG_1}{z_2-z_1}\biggr)\Tr \calG_2\biggr\} - \frac{1}{z_1} \E \Tr \calG_2 + \bigo(\sqrt{p})\\
    =&\; \frac{1}{z_1 k_2} \bigl(k_2 - p - z_1 \E \Tr G_1\bigr) \E \Tr (\calG_1 \calG_2)
    - \frac{1}{z_1 k_2} \biggl(\frac{z_2 \E \Tr \calG_2}{z_2 - z_1} - \frac{z_1 \E \Tr \calG_1}{z_2 - z_1}\biggr) \E \Tr \calG_2 - \frac{1}{z_1} \E \Tr \calG_2 + \bigo(\sqrt{p}),
\end{align*}
where the second equality following from \eqref{eq:Var_traceG}, \eqref{eq:cov_traceG_traceGS_proof}, and 
\begin{align}\label{eq:Var_Cov_traceG1G2}
    \var\bigl( \Tr(\calG_1\calG_2) \bigr) = \bigo(1),
    \quad
    \cov\bigl(\Tr \calG_1, \Tr \calG_2\bigr) = \bigo(1).
\end{align}
The proof of these two estimates is similar to that of \eqref{eq:Var_traceG} and is thus omitted.
Combining the above estimate with the identity $G_1 G_2 = (G_1-G_2)/(z_1-z_2)$, we have 
\begin{align*}
    \frac{z_2\E\Tr\calG_2-z_1\E\Tr\calG_1}{z_2-z_1} \biggl(1+\frac{1}{k_2}\E\Tr\calG_2\biggr) = \biggl(1-\frac{p}{k_2}-\frac{z_1}{k_2}\E\Tr G_1\biggr) \E \Tr\calG_1\calG_2 + \bigo(\sqrt{p}).
\end{align*}
By symmetry, we also have
\begin{align*}
    \frac{z_2\E\Tr\calG_2-z_1\E\Tr\calG_1}{z_2-z_1} \biggl(1+\frac{1}{k_2}\E\Tr\calG_1\biggr) = \biggl(1-\frac{p}{k_2}-\frac{z_2}{k_2}\E\Tr G_2\biggr) \E \Tr\calG_1\calG_2 + \bigo(\sqrt{p}).
\end{align*}
Combining the above two identities yields that
\begin{align*}
    &\; \frac{1}{k_2} \E \Tr (\calG_1 \calG_2) = \frac{\bigl(z_2\E\Tr \calG_2- z_1\E\Tr \calG_1\bigr) \bigl(\E \Tr \calG_2 - \E \Tr \calG_1\bigr)}{k_2 (z_2-z_1)\bigl(z_2 \E\Tr G_2 - z_1 \E \Tr G_1\bigr)}  + \bigo(1/\sqrt{p}).
\end{align*}
This, together with \eqref{eq:trace_G}, \eqref{eq:trace_GS}, and \eqref{eq:Var_Cov_traceG1G2}, completes the proof of \eqref{eq:trace_G1S1G2S1}.

\end{appendix}

\end{document}

%% file: table/table_H0_results.tex
\begin{table}[ht]
\centering
\begin{tabular}{llcccc}
\toprule
& & \multicolumn{4}{c}{Test Functions} \\
\cmidrule(lr){3-6}
Distribution & Weight & \texttt{linear} & \texttt{log} & \texttt{mix} & \texttt{square} \\
\midrule
\multirow{4}{*}{\texttt{Gaussian}}
 & $\rho_{1,0}$ & 0.0840 & 0.0630 & 0.0620 & 0.1315 \\
 & $\rho_{1,0.25}$ & 0.0855 & 0.0645 & 0.0680 & 0.1415 \\
 & $\rho_{1,0.45}$ & 0.0825 & 0.0725 & 0.0830 & 0.1660 \\
 & $\rho_{2}$ & 0.0135 & 0.0100 & 0.0090 & 0.0315 \\
\midrule
\multirow{4}{*}{\texttt{Uniform}}
 & $\rho_{1,0}$ & 0.0535 & 0.0450 & 0.0570 & 0.1235 \\
 & $\rho_{1,0.25}$ & 0.0600 & 0.0500 & 0.0670 & 0.1290 \\
 & $\rho_{1,0.45}$ & 0.0710 & 0.0570 & 0.0780 & 0.1410 \\
 & $\rho_{2}$ & 0.0055 & 0.0045 & 0.0070 & 0.0370 \\
\midrule
\multirow{4}{*}{\texttt{Student's $t(10)$}}
 & $\rho_{1,0}$ & 0.0755 & 0.0890 & 0.0715 & 0.1585 \\
 & $\rho_{1,0.25}$ & 0.0835 & 0.0885 & 0.0805 & 0.1645 \\
 & $\rho_{1,0.45}$ & 0.1040 & 0.0810 & 0.0905 & 0.1715 \\
 & $\rho_{2}$ & 0.0100 & 0.0190 & 0.0110 & 0.0455 \\
\bottomrule
\end{tabular}
\caption{Empirical size under different test function, weighted functions and data distributions. The nominal level is $\alpha=0.05$.}
\label{tab:empirical_size}
\end{table}

%% file: table/table_H1_results_gamma_0.tex
\begin{table}[ht]
\centering
\begin{tabular}{cclllll}
\toprule
& & \multicolumn{5}{c}{Change Magnitude ($\sigma^2$)} \\
\cmidrule(lr){3-7}
& $k^{\star}$ & \multicolumn{1}{c}{1.1} & \multicolumn{1}{c}{1.2} & \multicolumn{1}{c}{1.3} & \multicolumn{1}{c}{1.4} & \multicolumn{1}{c}{1.5} \\
\midrule
\multicolumn{2}{c}{}& \multicolumn{5}{l}{\texttt{Gaussian}} \\
\multirow{3}{*}{\texttt{linear}} & 350 & 114.84 (0.987) & 42.10 (1.0) & 25.53 (1.0) & 18.37 (1.0) & 14.50 (1.0) \\
& 450 & 133.16 (0.9925) & 52.22 (1.0) & 32.08 (1.0) & 23.69 (1.0) & 18.73 (1.0) \\
& 550 & 157.06 (0.9745) & 63.04 (1.0) & 39.44 (1.0) & 28.03 (1.0) & 22.36 (1.0) \\
\cmidrule(lr){2-7}
\multirow{3}{*}{\texttt{log}} & 350 & 53.96 (1.0) & 24.33 (1.0) & 16.10 (1.0) & 12.10 (1.0) & 9.87 (1.0) \\
& 450 & 63.87 (1.0) & 29.15 (1.0) & 18.83 (1.0) & 13.97 (1.0) & 11.45 (1.0) \\
& 550 & 75.65 (1.0) & 33.84 (1.0) & 22.18 (1.0) & 16.60 (1.0) & 13.27 (1.0) \\
\cmidrule(lr){2-7}
\multirow{3}{*}{\texttt{mix}} & 350 & 98.90 (0.9985) & 37.22 (1.0) & 23.32 (1.0) & 17.13 (1.0) & 13.74 (1.0) \\
& 450 & 118.16 (0.997) & 47.32 (1.0) & 29.61 (1.0) & 21.27 (1.0) & 16.93 (1.0) \\
& 550 & 140.61 (0.993) & 56.00 (1.0) & 35.39 (1.0) & 26.19 (1.0) & 20.52 (1.0) \\
\midrule
\multicolumn{2}{c}{}& \multicolumn{5}{l}{\texttt{Uniform}} \\
\multirow{3}{*}{\texttt{linear}} & 350 & 94.60 (0.995) & 36.52 (1.0) & 22.88 (1.0) & 16.84 (1.0) & 13.20 (1.0) \\
& 450 & 116.13 (0.9975) & 45.39 (1.0) & 28.79 (1.0) & 21.16 (1.0) & 16.74 (1.0) \\
& 550 & 138.91 (0.9915) & 55.67 (1.0) & 34.93 (1.0) & 25.63 (1.0) & 19.95 (1.0) \\
\cmidrule(lr){2-7}
\multirow{3}{*}{\texttt{log}} & 350 & 40.25 (1.0) & 19.04 (1.0) & 12.45 (1.0) & 9.53 (1.0) & 7.76 (1.0) \\
& 450 & 48.18 (1.0) & 22.34 (1.0) & 14.78 (1.0) & 11.02 (1.0) & 9.07 (1.0) \\
& 550 & 56.20 (1.0) & 26.21 (1.0) & 16.98 (1.0) & 12.88 (1.0) & 10.42 (1.0) \\
\cmidrule(lr){2-7}
\multirow{3}{*}{\texttt{mix}} & 350 & 82.20 (0.999) & 32.60 (1.0) & 20.57 (1.0) & 15.04 (1.0) & 12.04 (1.0) \\
& 450 & 100.18 (1.0) & 40.77 (1.0) & 26.03 (1.0) & 18.95 (1.0) & 14.92 (1.0) \\
& 550 & 116.64 (0.9985) & 49.29 (1.0) & 30.74 (1.0) & 22.70 (1.0) & 18.23 (1.0) \\
\midrule
\multicolumn{2}{c}{}& \multicolumn{5}{l}{\texttt{Student's $t(10)$}} \\
\multirow{3}{*}{\texttt{linear}} & 350 & 127.87 (0.9785) & 45.01 (1.0) & 27.91 (1.0) & 19.90 (1.0) & 15.84 (1.0) \\
& 450 & 151.23 (0.9805) & 57.80 (1.0) & 34.65 (1.0) & 25.93 (1.0) & 20.45 (1.0) \\
& 550 & 172.12 (0.958) & 68.46 (1.0) & 42.98 (1.0) & 30.66 (1.0) & 24.19 (1.0) \\
\cmidrule(lr){2-7}
\multirow{3}{*}{\texttt{log}} & 350 & 62.49 (1.0) & 28.16 (1.0) & 18.42 (1.0) & 13.83 (1.0) & 11.24 (1.0) \\
& 450 & 74.65 (1.0) & 33.34 (1.0) & 21.06 (1.0) & 15.78 (1.0) & 13.04 (1.0) \\
& 550 & 88.25 (1.0) & 39.10 (1.0) & 24.90 (1.0) & 18.51 (1.0) & 15.07 (1.0) \\
\cmidrule(lr){2-7}
\multirow{3}{*}{\texttt{mix}} & 350 & 112.59 (0.99) & 41.41 (1.0) & 25.37 (1.0) & 18.64 (1.0) & 14.71 (1.0) \\
& 450 & 136.76 (0.992) & 50.99 (1.0) & 31.97 (1.0) & 23.33 (1.0) & 18.25 (1.0) \\
& 550 & 156.53 (0.9795) & 62.42 (1.0) & 38.80 (1.0) & 28.44 (1.0) & 22.79 (1.0) \\
\bottomrule
\end{tabular}
\caption{EDD and Power (in parentheses) for \emph{homogeneous variance inflation} using weight $\rho_{1,\gamma}$ ($\gamma=0$).}
\label{tab:inflation_gamma_00}
\end{table}

\begin{table}[ht]
\centering
\begin{tabular}{cclllll}
\toprule
& & \multicolumn{5}{c}{Change Magnitude ($\rho$)} \\
\cmidrule(lr){3-7}
& $k^{\star}$ & \multicolumn{1}{c}{0.1} & \multicolumn{1}{c}{0.3} & \multicolumn{1}{c}{0.5} & \multicolumn{1}{c}{0.7} & \multicolumn{1}{c}{0.9} \\
\midrule
\multicolumn{2}{c}{}& \multicolumn{5}{l}{\texttt{Gaussian}} \\
\multirow{3}{*}{\texttt{linear}} & 350 & 7.33 (1.0) & 7.31 (1.0) & 7.33 (1.0) & 7.42 (1.0) & 7.36 (1.0) \\
& 450 & 9.05 (1.0) & 9.20 (1.0) & 9.24 (1.0) & 9.17 (1.0) & 9.27 (1.0) \\
& 550 & 11.15 (1.0) & 11.13 (1.0) & 10.93 (1.0) & 11.07 (1.0) & 11.26 (1.0) \\
\cmidrule(lr){2-7}
\multirow{3}{*}{\texttt{log}} & 350 & 5.37 (1.0) & 5.36 (1.0) & 5.31 (1.0) & 5.42 (1.0) & 5.57 (1.0) \\
& 450 & 6.03 (1.0) & 6.06 (1.0) & 6.09 (1.0) & 6.04 (1.0) & 6.30 (1.0) \\
& 550 & 6.99 (1.0) & 6.99 (1.0) & 6.93 (1.0) & 7.13 (1.0) & 7.19 (1.0) \\
\cmidrule(lr){2-7}
\multirow{3}{*}{\texttt{mix}} & 350 & 6.85 (1.0) & 6.90 (1.0) & 6.84 (1.0) & 6.93 (1.0) & 7.09 (1.0) \\
& 450 & 8.66 (1.0) & 8.52 (1.0) & 8.51 (1.0) & 8.56 (1.0) & 8.67 (1.0) \\
& 550 & 10.26 (1.0) & 10.43 (1.0) & 10.20 (1.0) & 10.28 (1.0) & 10.31 (1.0) \\
\midrule
\multicolumn{2}{c}{}& \multicolumn{5}{l}{\texttt{Uniform}} \\
\multirow{3}{*}{\texttt{linear}} & 350 & 6.67 (1.0) & 6.55 (1.0) & 6.67 (1.0) & 6.65 (1.0) & 6.70 (1.0) \\
& 450 & 8.43 (1.0) & 8.28 (1.0) & 8.37 (1.0) & 8.44 (1.0) & 8.45 (1.0) \\
& 550 & 10.08 (1.0) & 10.03 (1.0) & 10.08 (1.0) & 10.21 (1.0) & 10.19 (1.0) \\
\cmidrule(lr){2-7}
\multirow{3}{*}{\texttt{log}} & 350 & 4.25 (1.0) & 4.28 (1.0) & 4.30 (1.0) & 4.32 (1.0) & 4.50 (1.0) \\
& 450 & 4.83 (1.0) & 4.81 (1.0) & 4.76 (1.0) & 4.83 (1.0) & 5.00 (1.0) \\
& 550 & 5.58 (1.0) & 5.60 (1.0) & 5.55 (1.0) & 5.62 (1.0) & 5.76 (1.0) \\
\cmidrule(lr){2-7}
\multirow{3}{*}{\texttt{mix}} & 350 & 6.03 (1.0) & 6.17 (1.0) & 6.08 (1.0) & 6.11 (1.0) & 6.18 (1.0) \\
& 450 & 7.59 (1.0) & 7.66 (1.0) & 7.57 (1.0) & 7.64 (1.0) & 7.73 (1.0) \\
& 550 & 9.11 (1.0) & 8.92 (1.0) & 9.16 (1.0) & 9.22 (1.0) & 9.38 (1.0) \\
\midrule
\multicolumn{2}{c}{}& \multicolumn{5}{l}{\texttt{Student's $t(10)$}} \\
\multirow{3}{*}{\texttt{linear}} & 350 & 7.86 (1.0) & 7.95 (1.0) & 7.89 (1.0) & 7.92 (1.0) & 7.90 (1.0) \\
& 450 & 9.89 (1.0) & 10.01 (1.0) & 9.88 (1.0) & 9.94 (1.0) & 10.02 (1.0) \\
& 550 & 11.96 (1.0) & 11.76 (1.0) & 11.84 (1.0) & 12.14 (1.0) & 12.12 (1.0) \\
\cmidrule(lr){2-7}
\multirow{3}{*}{\texttt{log}} & 350 & 6.05 (1.0) & 6.05 (1.0) & 6.00 (1.0) & 6.13 (1.0) & 6.35 (1.0) \\
& 450 & 6.72 (1.0) & 6.74 (1.0) & 6.83 (1.0) & 6.97 (1.0) & 7.08 (1.0) \\
& 550 & 7.87 (1.0) & 7.70 (1.0) & 7.87 (1.0) & 7.99 (1.0) & 8.25 (1.0) \\
\cmidrule(lr){2-7}
\multirow{3}{*}{\texttt{mix}} & 350 & 7.46 (1.0) & 7.39 (1.0) & 7.45 (1.0) & 7.53 (1.0) & 7.51 (1.0) \\
& 450 & 9.25 (1.0) & 9.21 (1.0) & 9.29 (1.0) & 9.17 (1.0) & 9.38 (1.0) \\
& 550 & 11.04 (1.0) & 11.22 (1.0) & 11.13 (1.0) & 11.07 (1.0) & 11.38 (1.0) \\
\bottomrule
\end{tabular}
\caption{EDD and Power (in parentheses) for \emph{correlation structure change} using weight $\rho_{1,\gamma}$ ($\gamma=0$).}
\label{tab:ar1_gamma_00}
\end{table}

\begin{table}[ht]
\centering
\begin{tabular}{cclllll}
\toprule
& & \multicolumn{5}{c}{Change Magnitude ($\delta$)} \\
\cmidrule(lr){3-7}
& $k^{\star}$ & \multicolumn{1}{c}{2} & \multicolumn{1}{c}{2.5} & \multicolumn{1}{c}{3} & \multicolumn{1}{c}{3.5} & \multicolumn{1}{c}{4} \\
\midrule
\multicolumn{2}{c}{}& \multicolumn{5}{l}{\texttt{Gaussian}} \\
\multirow{3}{*}{\texttt{linear}} & 350 & 14.01 (1.0) & 13.34 (1.0) & 12.74 (1.0) & 12.17 (1.0) & 11.67 (1.0) \\
& 450 & 17.82 (1.0) & 16.74 (1.0) & 16.20 (1.0) & 15.25 (1.0) & 14.62 (1.0) \\
& 550 & 21.24 (1.0) & 20.25 (1.0) & 19.70 (1.0) & 18.86 (1.0) & 17.76 (1.0) \\
\cmidrule(lr){2-7}
\multirow{3}{*}{\texttt{log}} & 350 & 9.53 (1.0) & 9.05 (1.0) & 8.59 (1.0) & 8.33 (1.0) & 8.06 (1.0) \\
& 450 & 10.79 (1.0) & 10.19 (1.0) & 10.06 (1.0) & 9.68 (1.0) & 9.26 (1.0) \\
& 550 & 12.58 (1.0) & 12.07 (1.0) & 11.65 (1.0) & 11.32 (1.0) & 10.66 (1.0) \\
\cmidrule(lr){2-7}
\multirow{3}{*}{\texttt{mix}} & 350 & 13.01 (1.0) & 12.30 (1.0) & 11.72 (1.0) & 11.42 (1.0) & 10.95 (1.0) \\
& 450 & 16.08 (1.0) & 15.33 (1.0) & 14.77 (1.0) & 14.47 (1.0) & 13.63 (1.0) \\
& 550 & 19.44 (1.0) & 18.48 (1.0) & 17.94 (1.0) & 17.07 (1.0) & 16.47 (1.0) \\
\midrule
\multicolumn{2}{c}{}& \multicolumn{5}{l}{\texttt{Uniform}} \\
\multirow{3}{*}{\texttt{linear}} & 350 & 12.49 (1.0) & 11.90 (1.0) & 11.43 (1.0) & 10.96 (1.0) & 10.50 (1.0) \\
& 450 & 15.77 (1.0) & 14.98 (1.0) & 14.39 (1.0) & 13.84 (1.0) & 13.27 (1.0) \\
& 550 & 19.44 (1.0) & 18.20 (1.0) & 17.85 (1.0) & 16.71 (1.0) & 16.13 (1.0) \\
\cmidrule(lr){2-7}
\multirow{3}{*}{\texttt{log}} & 350 & 7.42 (1.0) & 7.14 (1.0) & 6.87 (1.0) & 6.61 (1.0) & 6.29 (1.0) \\
& 450 & 8.49 (1.0) & 8.19 (1.0) & 7.94 (1.0) & 7.66 (1.0) & 7.38 (1.0) \\
& 550 & 10.03 (1.0) & 9.50 (1.0) & 9.12 (1.0) & 8.70 (1.0) & 8.32 (1.0) \\
\cmidrule(lr){2-7}
\multirow{3}{*}{\texttt{mix}} & 350 & 11.41 (1.0) & 10.92 (1.0) & 10.53 (1.0) & 10.04 (1.0) & 9.49 (1.0) \\
& 450 & 14.41 (1.0) & 13.69 (1.0) & 12.97 (1.0) & 12.53 (1.0) & 12.03 (1.0) \\
& 550 & 17.28 (1.0) & 16.59 (1.0) & 15.64 (1.0) & 14.98 (1.0) & 14.56 (1.0) \\
\midrule
\multicolumn{2}{c}{}& \multicolumn{5}{l}{\texttt{Student's $t(10)$}} \\
\multirow{3}{*}{\texttt{linear}} & 350 & 15.13 (1.0) & 14.31 (1.0) & 13.57 (1.0) & 13.16 (1.0) & 12.57 (1.0) \\
& 450 & 19.01 (1.0) & 18.24 (1.0) & 17.22 (1.0) & 16.62 (1.0) & 15.96 (1.0) \\
& 550 & 23.00 (1.0) & 22.21 (1.0) & 20.94 (1.0) & 20.25 (1.0) & 19.31 (1.0) \\
\cmidrule(lr){2-7}
\multirow{3}{*}{\texttt{log}} & 350 & 10.71 (1.0) & 10.27 (1.0) & 9.96 (1.0) & 9.41 (1.0) & 9.14 (1.0) \\
& 450 & 12.03 (1.0) & 11.65 (1.0) & 11.26 (1.0) & 10.96 (1.0) & 10.43 (1.0) \\
& 550 & 14.29 (1.0) & 14.12 (1.0) & 13.10 (1.0) & 12.55 (1.0) & 12.07 (1.0) \\
\cmidrule(lr){2-7}
\multirow{3}{*}{\texttt{mix}} & 350 & 13.90 (1.0) & 13.33 (1.0) & 12.83 (1.0) & 12.32 (1.0) & 11.77 (1.0) \\
& 450 & 17.66 (1.0) & 16.98 (1.0) & 16.06 (1.0) & 15.26 (1.0) & 14.78 (1.0) \\
& 550 & 21.46 (1.0) & 20.36 (1.0) & 18.96 (1.0) & 18.38 (1.0) & 18.03 (1.0) \\
\bottomrule
\end{tabular}
\caption{EDD and Power (in parentheses) for \emph{heterogeneous variance inflation} using weight $\rho_{1,\gamma}$ ($\gamma=0$).}
\label{tab:spike_gamma_00}
\end{table}

%% file: table/table_comparison.tex
\begin{table}[htbp]
\centering
\begin{tabular}{lccccc}
\toprule
\multicolumn{6}{c}{Homogeneous Variance Inflation (Change Magnitude $\sigma^2$)} \\
\cmidrule(lr){2-6}
Method & 1.1 & 1.2 & 1.3 & 1.4 & 1.5 \\
\midrule
\textsf{Our} & 82.03 (0.992) & 30.79 (0.992) & 18.29 (0.992) & 13.65 (0.998) & 11.01 (0.998) \\
\textsf{A19} & 126.65 (0.794) & 59.25 (0.798) & 36.45 (0.822) & 27.00 (0.838) & 22.16 (0.830) \\
\textsf{LL23} & 402.33 (0.332) & 343.44 (0.784) & 222.37 (0.956) & 126.45 (0.988) & 78.70 (0.996) \\
\midrule
\multicolumn{6}{c}{Correlation Structure Change (Change Magnitude $\rho$)} \\
\cmidrule(lr){2-6}
Method & 0.1 & 0.3 & 0.5 & 0.7 & 0.9 \\
\midrule
\textsf{Our} & 5.75 (0.998) & 5.85 (0.992) & 5.78 (0.992) & 5.89 (0.998) & 6.23 (0.998) \\
\textsf{A19} & 10.72 (0.792) & 10.65 (0.792) & 10.43 (0.818) & 9.81 (0.834) & 9.23 (0.842) \\
\textsf{LL23} & 17.68 (0.990) & 16.19 (0.988) & 13.54 (0.990) & 9.79 (0.988) & 6.35 (0.996) \\
\midrule
\multicolumn{6}{c}{Heterogeneous Variance Inflation (Change Magnitude $\delta$)} \\
\cmidrule(lr){2-6}
Method & 2 & 2.5 & 3 & 3.5 & 4 \\
\midrule
\textsf{Our} & 9.95 (0.998) & 9.43 (0.992) & 8.47 (0.992) & 8.00 (0.998) & 7.51 (0.998) \\
\textsf{A19} & 18.38 (0.792) & 15.10 (0.794) & 11.97 (0.828) & 9.99 (0.836) & 8.42 (0.848) \\
\textsf{LL23} & 57.93 (0.990) & 38.11 (0.988) & 25.80 (0.990) & 20.84 (0.988) & 16.90 (0.996) \\
\bottomrule
\end{tabular}
\caption{Comparison of EDD and Power (in parentheses) with existing methods (\textsf{A19}: \citet{Avanesov2019structural}, \textsf{LL23}: \citet{Li2023online}).}
\label{tab:comparison}
\end{table}